\documentclass[11pt,a4paper,twoside]{article}
\usepackage{a4,subfigure}
\usepackage{color}
\usepackage{bm, amsmath, amssymb, amsthm} 
\usepackage{graphicx}

\def\calf{{\cal F}}

\def\<{\langle}
\def\>{\rangle}
\def\eps{\epsilon}
\def\NN{\mathbb{N}}

\def\RR{\mathbb{R}}

\def\tr{\operatorname{Tr}}
\def\id{\operatorname{id\,}}
\def\Div{\operatorname{div}}

\newcommand{\dtau}{\operatorname{d}\!\tau}

\def\Ric{\operatorname{Ric}}

\def\eq{\hspace*{-1.5mm}&=&\hspace*{-1.5mm}}
\def\plus{\hspace*{-1.5mm}&+&\hspace*{-1.5mm}}
\def\minus{\hspace*{-1.5mm}&-&\hspace*{-1.5mm}}
\def\dt{\partial_t}

\newtheorem{corollary}{Corollary}
\newtheorem{definition}{Definition}
\newtheorem{example}{Example}
\newtheorem{remark}{Remark}
\newtheorem{lemma}{Lemma}
\newtheorem{proposition}{Proposition}
\newtheorem{theorem}{Theorem}

\author{Vladimir Rovenski\thanks{
        E-mail: rovenski@math.haifa.ac.il,\quad zelenko@math.haifa.ac.il.
        }
        \ and \
        Leonid Zelenko
        \\
        {\small Mathematical Department, University of Haifa, Mount Carmel, Haifa, 31905, Israel}
}

\title{Prescribing the mixed scalar curvature \\ of a foliated Riemann-Cartan manifold}

\begin{document}

\date{}

\maketitle

\begin{abstract}
The mixed scalar curvature is one of the simplest curvature invariants of a foliated Riemannian manifold.
We~explore the problem of prescribing the mixed scalar curvature of a foliated Riemann-Cartan manifold
by conformal change of the structure in tangent and normal to the leaves directions.
Under certain geometrical assumptions and in two special cases: along a compact leaf and for a closed fibred manifold,
we reduce the problem to solution of a leafwise elliptic equation, which has three stable solutions
-- only one of them corresponds to the case of a foliated Riemannian manifold.



\end{abstract}

\tableofcontents

\section*{Introduction}
\label{sec:intro}

Geometrical problems of prescribing curvature invariants
of Riemannian manifolds using conformal change of metric are popular for a long time,
i.e., the study of constancy of the scalar curvature
was began by Yamabe in 1960 and completed by Trudinger, Aubin and Schoen in 1986, see~\cite{aub}.

 The metrically-affine geometry was founded by E.\,Cartan in 1923--1925, who suggested
using an asymmetric connection $\bar\nabla$ instead of Levi-Civita connection $\nabla$ of $g$;
in~extended theory of gravity the torsion of $\bar\nabla$ is represented by the spin tensor of matter.
Notice that $\bar\nabla$ and $\nabla$ are projectively equivalent (have the same geodesics) if and only if
the difference $\mathfrak{T}:=\bar\nabla-\nabla$, called the {contorsion tensor}, is antisymmetric.
 \textit{Riemann-Cartan} ({RC}) spaces, i.e., with metric connection: $\bar\nabla g =0$, appear in such topics as
homogeneous and almost Hermitian spaces \cite{gps}, and geometric flows~\cite{ah}.

Foliations, i.e., partitions of a manifold into collection of submanifolds, called leaves,
arise in topology and have applications in differential geometry, analysis
and theoretical physics, where many models are foliated.
One of the simplest curvature invariants of a foliated Riemannian manifold is
the mixed scalar curvature ${\rm S}_{\rm mix}$, i.e., an averaged sectional curvature over all planes that contain vectors from both
 -- tangent and normal -- distributions, see~\cite{rov-m}.
The
prescribing of
${\rm S}_{\rm mix}$
by conformal change of the metric in normal to the leaves directions
and certain Yamabe type problem have been studied in~\cite{rz2013,rz1}.
 In the paper, we examine the problem of prescribing the mixed scalar curvature $\bar{\rm S}_{\,\rm mix}$ of a foliated RC manifold
by conformal change of the structure in tangent and normal to the leaves directions.
In particular, we explore the following {Yamabe type~problem}:

\textit{Given foliated {RC} manifold $(M,g,\bar\nabla)$
find a $({\cal D},{\cal D}^\perp)$-conformal {RC} structure, i.e.,
\begin{equation}\label{E-conf-RC}
 g'=g^\top+u^2\,g^\bot,\qquad
 {\mathfrak{T}'}=u^2\,\mathfrak{T}^\top+\mathfrak{T}^\bot,
\end{equation}
with (leafwise) constant~mixed scalar~curvature.} Here $u\in C^\infty(M)$ is positive and
\begin{equation*}
 g^\top(X,Y):=g(X^\top,Y^\top),\quad
 g^\perp(X,Y):=g(X^\perp,Y^\perp),\quad
 \mathfrak{T}^\top_XY:=(\mathfrak{T}_X Y)^\top,\quad \mathfrak{T}^\bot_XY:=(\mathfrak{T}_XY)^\bot.
\end{equation*}
We show that under certain geometric assumptions, including $\nabla$-harmonicy of $\calf$ and $g_{\,|\,T\calf}>0$,
the conformal factor in \eqref{E-conf-RC} obeys leafwise elliptic equation
\begin{equation}\label{E-Yam1-init}
 -\Delta^\top u -\beta(x)\,u =  \Psi_1(x)\,u^{-1} -\Psi_2(x)\,u^{-3} +\Psi_3(x)\,u^{3}
\end{equation}
with smooth functions $\beta$ and $\Psi_i\ (i=1,2,3)$ described in Section~\ref{subsec:prel}.
The case of $g_{\,|T\calf}<0$ reduces to the above by change $g\rightarrowtail-g$.
Notice that $\Psi_3$ represents the novel \textit{mixed scalar $\mathfrak{T}^\top$-curvature}, see Section~\ref{sec:Smix},
and the stable solution of \eqref{E-Yam1-init} in the case of $\Psi_3=0$ has been found in \cite{rz1}.
 By stable solution of elliptic equation we mean a stable stationary solution of its parabolic counterpart.

Using
spectral parameters of the Schr\"{o}dinger operator along compact leaves,
\begin{equation}\label{E-Schred1}
 \mathcal{H}: u \mapsto -\Delta^\top u -\beta(x)\,u,
\end{equation}
we prove that \eqref{E-Yam1-init} has three stable solutions,
one of them ($\Psi_3=0$) corresponds to the Riemannian~case.

 Since the topology of the leaf through a point can change dramatically with the point,
there are difficulties in studying leafwise elliptic equations.
 Thus, we examine two
 formulations of the problem:

\noindent\ \
{1}. $\bar{\rm S}_{\,\rm mix}$ \underline{is prescribed on a compact leaf} $F$.
Under some geometric assumptions we get \eqref{E-Yam1-init},
whose solutions $u_*{\in} C^\infty(F)$ form a compact set in $C(F)$ and can be extended smoothly onto~$M$.

\noindent\ \
{2}. $\bar{\rm S}_{\,\rm mix}$ \underline{is prescribed on a closed manifold} $M$.
Under certain geometric assumptions we get \eqref{E-Yam1-init} on any~$F$,
whose unique solution $u_*\in C^\infty(F)$ on any leaf $F$   belongs to $C^\infty(M)$ when
\begin{equation}\label{E-2conditions}
 {\calf} \ \ \textrm{is defined by an orientable fiber bundle}\ \ {\pi: M\to B.}
\end{equation}
The main results of the paper are Theorems~\ref{T-A}--\,\ref{T-mainB1-R} (and their corollaries) about foliations of arbitrary (co)dimension,
similar results for codimension-one foliations and flows are omitted.

The paper is organized as follows.
Section~\ref{sec:prelim} contains geometrical results of our paper.
Section~\ref{sec:Smix} gives preliminaries for foliated {RC} manifolds.
Section~\ref{subsec:prel} derives the transformation law for $\bar{\rm S}_{\,\rm mix}$
under $({\cal D},{\cal D}^\bot)$-conformal change of RC structure;
this yields, under certain geometrical assumptions,
elliptic equation \eqref{E-Yam1-init} for the conformal factor.
The results in Section~\ref{sec:RC-space} are separated into three cases
according the sign of the mixed scalar $\mathfrak{T}^\top$-curvature represented by $\Psi_3$.
To prescribe $\bar{\rm S}_{\,\rm mix}$ on a closed leaf (Theorem~\ref{T-A}) we use the existence of a solution to \eqref{E-Yam1-init},
and to~prescribe $\bar{\rm S}_{\,\rm mix}$ on a closed fibred manifold (Theorem~\ref{T-B}) we
use the existence and uniqueness of a solution to \eqref{E-Yam1-init}, see Section~\ref{sec:app2},
where we also prove that \eqref{E-Yam1-init} has three stable solutions,
 which are expressed in terms of spectral parameters of operator~\eqref{E-Schred1}.

\section{Foliated Riemann-Cartan manifolds}
\label{sec:prelim}

\subsection{The mixed scalar curvature}
\label{sec:Smix}

A pseudo-Riemannian metric of index $q$ on manifold $M$ is an element $g\in{\rm Sym}^2(M)$
(of the space of symmetric $(0,2)$-tensor fields)
 such that each $g_x\ (x\in M)$ is a {non-degenerate bilinear form of index} $q$ on the tangent space $T_xM$.
When~$q=0$, $g$ is a Riemannian metric (resp. a Lorentz metric~when $q=1$).
Let $\mathfrak{X}_M$ be the module over $C^\infty(M)$ of all vector fields on~$M$.

 The~Levi-Civita connection $\nabla: \mathfrak{X}_M\times \mathfrak{X}_M\to \mathfrak{X}_M$ of $g$, represented using the Lie~bracket,
\begin{eqnarray}\label{eqlevicivita}
\nonumber
 &&\quad 2\,g(\nabla_X\,Y, Z) = X g(Y,Z) + Y g(X,Z) - Z\,g(X,Y) \\
     && +\,g([X, Y], Z) - g([X, Z], Y) - g([Y, Z], X)\quad (X,Y,Z\in\mathfrak{X}_M),
\end{eqnarray}
is \textit{metric compatible}, $\nabla g=0$,  and has zero torsion.

A subbundle ${\cal D}\subset TM$ (called a distribution) is non-degenerate,
if ${\cal D}_x$ is a non-degenerate subspace of $(T_x M,\, g_x)$ for
$x\in M$;
in this case, its orthogonal distribution ${\cal D}^\bot\subset TM$
is also non-degenerate.
Thus, we consider a connected manifold $M^{n+p}$ with a pseudo-Riemannian metric $g$ and
complementary orthogonal non-degene\-rate distributions ${\cal D}$ and ${\cal D}^\bot$ of ranks
$\dim_{\,\RR}{\cal D}_x=p\ge1$ and $\dim_{\,\RR}{\cal D}^\bot_x=n\ge1$ for every $x\in M$
(called an \textit{almost-product structure} on $M$), see~\cite{bf}.

Let $X^\top$ be the ${\cal D}$-component of $X\in\mathfrak{X}_M$ (resp., $X^\perp$ the ${\cal D}^\bot$-component of $X$),
and $\mathfrak{X}^\top_M$ (resp. $\mathfrak{X}^\bot_M$) the module over $C^\infty(M)$
of all vector fields in ${\cal D}$ (resp. ${\cal D}^\bot$).
In the paper, ${\cal D}$ is integrable and tangent to a foliation $\calf$.
The~integrability tensor of ${\cal D}^\bot$ is defined by $T^\bot(X,Y) =\frac12\,[X,\,Y]^\top\ (X,Y\in\mathfrak{X}^\bot_M)$,
the second fundamental forms of ${\cal D}$ and~${\cal D}^\bot$
are given~by
\begin{eqnarray*}
 h^\top(X,Y) = (\nabla_{X} Y)^\bot\ \ (X,Y\in\mathfrak{X}^\top_M),\quad
 h^\bot(X,Y) = (\nabla_{X} Y +\nabla_{Y} X)^\top/2\ \ (X,Y\in\mathfrak{X}^\bot_M),
\end{eqnarray*}
and mean curvature vectors are
 $H^\top=\tr_g h^\top$
 and
 $H^\bot=\tr_g h^\bot$.
 We call ${\cal D}$
 \textit{totally umbilical}, \textit{harmonic}, or \textit{totally geodesic}, if
 $h^\top=\frac1p\,H^\top g^\top,\ H^\top =0$, or $h^\top=0$, resp.
Examples of harmonic foliations are parallel
circles or winding lines on a flat~torus
 and a Hopf field
 of great circles
 on $S^3$.

\smallskip

Recall that a linear connection $\bar\nabla$ on $M$ is a map
 $\bar\nabla:\mathfrak{X}_M\times\mathfrak{X}_M\to\mathfrak{X}_M$
 with the properties:
 \[
 \bar\nabla_{f X_1+X_2} Y =f\bar\nabla_{X_1} Y+\bar\nabla_{X_2} Y,\quad
 \bar\nabla_X(fY +Z) =X(f)Y+f\bar\nabla_{X}Y +\bar\nabla_{X}Z,
\]
where $f\in C^\infty(M)$.
Thus, linear connections over $M$ form an affine space, and the difference of two connections is a $(1,2)$-tensor.

Computing terms in the definition $\bar R_{X,Y}=[\bar\nabla_Y,\bar\nabla_X]+\bar\nabla_{[X,Y]}$
of the curvature tensor of $\bar\nabla=\nabla +\mathfrak{T}$,
and comparing with similar formula for $R_{X,Y}$, we find the following relation:
\begin{equation}\label{E-RC-2}
 \bar R_{X,Y} = R_{X,Y} +(\nabla_Y \mathfrak{T})_X -(\nabla_X \mathfrak{T})_Y +[\mathfrak{T}_Y,\mathfrak{T}_X].
\end{equation}

 Let $\{E_a,\,{\cal E}_i\}_{a\le p,\,i\le n}$ be a local orthonormal frame on $TM$ such that
 $\{E_a\}\subset{\cal D}$ and $\{{\cal E}_i\}\subset{\cal D}^\bot$
 and $\eps_a=g(E_a,E_a),\ \eps_i=g({\cal E}_{i},{\cal E}_{i})$.
 We use the following convention for various  tensors:
 $\mathfrak{T}_i=\mathfrak{T}_{{\cal E}_i}$ etc.
 The following function on a metric-affine manifold $(M,g,\bar\nabla)$:
\begin{equation}\label{eq-wal2}
 \bar{\rm S}_{\rm mix}
 =\frac{1}{2}\sum\nolimits_{a,i} \eps_a\eps_i\big( g({\bar R}_{E_a, {\cal E}_i} E_a, {\cal E}_i)
 + g({\bar R}_{{\cal E}_i, E_a} {\cal E}_i, E_a)\big)
\end{equation}
is well-defined and is called the \textit{mixed scalar curvature} of $(\widetilde{\cal D},{\cal D})$.
This definition does not depend on the order of distributions and on the choice of a local frame. Moreover, see \eqref{E-RC-2},
\begin{eqnarray}\label{E-RC-5}
\nonumber
 &&\hskip-9mm \bar{\rm S}_{\rm mix} = {\rm S}_{\rm mix} + Q,\ \ {\rm where}  \\
\nonumber
 &&\hskip-9mm Q = \frac12\sum\nolimits_{a,i} \eps_a \eps_i
 \big( g( (\nabla_i \mathfrak{T})_a E_a , {\cal E}_i ) - g( (\nabla_a \mathfrak{T})_i E_a, {\cal E}_i )
 + g( [\mathfrak{T}_i,\, \mathfrak{T}_a] E_a , {\cal E}_i ) \\
 && +\, g( (\nabla_a \mathfrak{T})_i {\cal E}_i, E_a ) - g( (\nabla_i \mathfrak{T})_a {\cal E}_i , E_a )
 + g( [\mathfrak{T}_a,\, \mathfrak{T}_i]\, {\cal E}_i , E_a ) \big)\,,
\end{eqnarray}
and ${\rm S}_{\rm mix}$ is the mixed scalar curvature of $\nabla$, see~\cite{rov-m,wa1}.
Recall the formula:
\begin{equation}\label{E-PW}
 {\rm S}_{\,\rm mix}=g(H^\bot,H^\bot)-\<h^\bot,h^\bot\>+\<T^\bot,T^\bot\>+g(H^\top,H^\top)-\<h^\top,h^\top\>
  +\Div(H^\bot + H^\top).
\end{equation}
For a vector field
on $M$ and for the gradient and Laplacian of a function $f\in C^\infty(M)$ we~have
\begin{eqnarray*}
 &&\hskip-6mm \Div X = \Div^\bot X +\Div^\top X,\quad
  \Div^\bot X=\sum\nolimits_{\,i}\eps_i\,g(\nabla_{{\cal E}_i} X,\,{\cal E}_i),\quad
 \Div^\top X=\sum\nolimits_{\,a}\eps_a\,g(\nabla_{E_a}X,\,E_a),\\
 &&\hskip-6mm g(\nabla f, X) = X(f),\qquad  \Delta\,f =\Div(\,\nabla f).
\end{eqnarray*}
We also use notations for traces of $\mathfrak{T}$:
 $\tr^\bot\mathfrak{T}:=\sum\nolimits_{\,i}\eps_i\,\mathfrak{T}_i\,{\cal E}_i$
 and
 $\tr^\top\mathfrak{T}:=\sum\nolimits_{\,a}\eps_a\,\mathfrak{T}_a E_a$.

\noindent\
Among all {metric-affine spaces} $(M,g,\bar\nabla)$,
\textit{{RC} spaces} have metric compatible connection, i.e.,
\begin{equation}\label{E-nabla-g}
 g(\mathfrak{T}_XY,Z) = -g(\mathfrak{T}_XZ,Y)\quad (X,Y,Z\in\mathfrak{X}_M).
\end{equation}
 The leaves of a foliation $\calf$ on $(M,g,\bar\nabla)$ are submanifolds with induced metric $g^\top$ and
metric connection $\bar\nabla^\top_X Y:=(\bar\nabla_X Y)^\top\ (X,Y\in \mathfrak{X}^\top_M)$.
Since, see~\eqref{E-nabla-g},
\[
 g^\top(\mathfrak{T}^\top_XY,Z) +g^\top(\mathfrak{T}^\top_XZ,Y) =
 g(\mathfrak{T}_XY,Z) +g(\mathfrak{T}_XZ,Y) =0\quad (X,Y,Z\in\mathfrak{X}^\top_M),
\]
the leaves (equipped with the metric $g^\top$ and connection $\bar\nabla^\top$) are themselves {RC} manifolds.
For {RC} spaces, the curvature tensor $\bar R$ has some symmetry properties, e.g.
\begin{equation}\label{E-Rsym}
 g(\bar R_{X,Y}Z,U) = -g(\bar R_{X,Y}U,Z),\quad g(\bar R_{X,Y}Z,U) = -g(\bar R_{Y,X}Z,U).
\end{equation}
The sectional curvature $\bar K(X\wedge Y)=g(\bar R_{X,Y}X,Y)/\,[g(X,X)g(Y,Y)-g(X,Y)^2]$
of {RC} spaces doesn't depend on the choice of a basis
in a non-degenerate plane $X\wedge Y$.
In this case, \eqref{E-RC-5}~reads
\begin{eqnarray}\label{E-RC-5b}
\nonumber
 &&\hskip-9mm \bar{\rm S}_{\rm mix} = {\rm S}_{\rm mix} + Q,\ \ {\rm where}  \\
\nonumber
 &&\hskip-9mm Q =
 \sum\nolimits_{\,i,a}\eps_i\,\eps_a\big[\,g((\nabla_{i} \mathfrak{T})_a E_a,{\cal E}_i)
  + g((\nabla_{a}\mathfrak{T})_i\,{\cal E}_i, E_a) + g(\mathfrak{T}_a{\cal E}_i, \mathfrak{T}_i E_a)
 - g(\mathfrak{T}_i\,{\cal E}_i, \mathfrak{T}_a E_a)\,\big]\,.
\end{eqnarray}
To show this, we use \eqref{E-RC-2}, \eqref{E-Rsym} and the equality
$g((\nabla_{X}\mathfrak{T})_YX,Y)=-g((\nabla_{X}\mathfrak{T})_YY,X)$, see \eqref{E-nabla-g}.

\begin{example}[\bf {RC} products]\label{Ex-2twist}\rm
The \textit{doubly-twisted product} of {RC} manifolds $(B,g_B,\mathfrak{T}_B)$ and $(F, g_F,\mathfrak{T}_F)$
is a manifold $M=B\times F$ with the metric $g = g^\top + g^\perp$ and the contorsion tensor
$\mathfrak{T} = \mathfrak{T}^\top + \mathfrak{T}^\bot$,~where
\begin{eqnarray*}
 g^\top(X,Y) \eq v^2 g_B(X^\top,Y^\top),\quad g^\bot(X,Y)=u^2 g_F(X^\bot,Y^\bot),\\
 \mathfrak{T}^\top_XY \eq u^2(\mathfrak{T}_B)_{X^\top}Y^\top,\quad \mathfrak{T}^\bot_XY=v^2(\mathfrak{T}_F)_{X^\bot}Y^\bot,
\end{eqnarray*}
and the warping functions $u,v\in C^\infty(M)$ are positive.
For $v=1$ we have the \textit{twisted product}, if, in addition, $u\in C^\infty(B)$
then this is a \textit{warped product},
and for $u=v=1$ -- the product of {RC} manifolds.
Let $g_B$ be positive definite.
 One may show that $\bar\nabla g=0$, see \eqref{E-nabla-g}, for the new connection $\bar\nabla=\nabla + \mathfrak{T}$:
\begin{eqnarray*}
  -\bar{\nabla}_X\,g(Y,Z) \eq (g^\top \!+ g^\perp)(\mathfrak{T}^\top_XY + \mathfrak{T}^\bot_XY,\,Z)
  + (g^\top \!+ g^\perp)( \mathfrak{T}^\top_XZ + \mathfrak{T}^\perp_XZ,\,Y)\\
 \eq u^{2} v^{2}[\,g_B((\mathfrak{T}_B)_{X^\top}Y^\top,Z^\top) + g_B((\mathfrak{T}_B)_{X^\top}Z^\top,Y^\top) \\
  \plus g_F((\mathfrak{T}_F)_{X^\bot}Y^\bot,Z^\bot) + g_F((\mathfrak{T}_F)_{X^\bot}Z^\bot,Y^\bot)\,] = 0\,.
\end{eqnarray*}
Hence, $(M,g,\bar\nabla)$ is a {RC} space, which will be denoted by $B\times_{(v,u)} F$.
Its second fundamental forms (w.r.t. $\nabla$) are, see \cite{pr},
 $h^\bot = -\nabla^\top(\log u)\,g^\bot$
 and
 $h^\top=-\nabla^\bot(\log v)\,g^\top$.
By~the above, $H^\bot=-n\,\nabla^\top(\log u)$ and $H^\top=-p\,\nabla^\bot(\log v)$.
 Hence, the \textit{leaves} $B\times\{y\}$ and the \textit{fibers} $\{x\}\times F$
of a {RC} doubly-twisted product $B\times_{(v,u)} F$ are totally umbilical
w.r.t. $\bar\nabla$ and $\nabla$. Since
\begin{eqnarray*}
 \Div\,H^\bot \eq -n\,(\Delta^\top u)/u -(n^2-n)\,g(\nabla^\top u,\nabla^\top u)/u^2,\\
 g(H^\bot,H^\bot)-\<h^\bot,h^\bot\> \eq (n^2-n)\,g(\nabla^\top u,\nabla^\top u)/u^2,\\
 \Div\,H^\top \eq -n\,(\Delta^\bot v)/v -(p^2-p)\,g(\nabla^\bot v,\nabla^\bot v)/v^2,\\
 g(H^\top,H^\top)-\<h^\top,h^\top\> \eq (p^2-p)\,g(\nabla^\bot v,\nabla^\bot v)/v^2,
\end{eqnarray*}
where $\Delta^\top$ is the leafwise Laplacian and $\Delta^\bot$ is the ${\cal D}^\bot$-Laplacian,
the formula \eqref{E-PW} reduces~to
 ${\rm S}_{\,\rm mix} = -n\,(\Delta^\top u)/u -p\,(\Delta^\bot\,v)/v$\,.
We have
 $Q = n\,u(\tr\mathfrak{T}^\top)(u) + p\,v(\tr\mathfrak{T}^\bot)(v)$, see \eqref{E-RC-5};
hence,
\begin{equation*}
 \bar{\rm S}_{\,\rm mix}= -n\,(\Delta^\top u)/u +n\,u(\tr\mathfrak{T}^\top)(u)
 -p\,(\Delta^\bot v)/v +p\,v(\tr\mathfrak{T}^\bot)(v).
\end{equation*}
The last formula is the linear PDE (with given $\bar{\rm S}_{\,\rm mix}$) along a leaf for unknown function $u$,
\begin{equation}\label{E-Smix_2twisted}
 -(\Delta^\top u) -\beta\,u + u^2(\tr\mathfrak{T}^\top)(u) = (\bar{\rm S}_{\,\rm mix}/n)\,u\,,
\end{equation}
where $\beta = \frac pn\,(v^{-1}\Delta^\perp\,v - v(\tr\mathfrak{T}^\bot)(v))$.
Let $B$ be a closed manifold, with $g_B>0$ and $\tr\mathfrak{T}_B=0$.
Thus, $\tr\mathfrak{T}^\top\!=0$, and \eqref{E-Smix_2twisted} becomes the eigenvalue problem.
 Thus, the product $B\times_{(v,e_0)} F$ has leafwise constant
$\bar{\rm S}_{\rm mix}$ equal to $n\lambda_0$, see \eqref{E-Schred1}.
For~$\mathfrak{T}_B=0=\mathfrak{T}_F$ we obtain Riemannian doubly-twisted products of leafwise constant ${\rm S}_{\,\rm mix}$,
see~\cite{rz1}.
\end{example}

In \cite{op2016}, the ${\cal K}$-sectional curvature of a symmetric $(1,2)$-tensor ${\cal K}$ is defined.
On this way, we introduce the following scalar invariant of a foliation.
For  a $(1,2)$-tensor ${\cal K}$, the \textit{mixed scalar ${\cal K}$-curvature} is defined by
\begin{equation*}
  {\rm S}_{{\rm mix},{\cal K}} := \frac12\sum\nolimits_{a,i} \eps_a\eps_i\big(
  g(\,[{\cal K}_i,\,{\cal K}_a]\,{\cal E}_i, E_a) + g(\,[{\cal K}_a,\,{\cal K}_i]\,E_a, {\cal E}_i)\big).
\end{equation*}
Note that the {mixed scalar $\mathfrak{T}$-curvature} in RC case is
  ${\rm S}_{{\rm mix},\mathfrak{T}} =  \sum\nolimits_{a,i} \eps_a\eps_i\,g(\,[\mathfrak{T}_i,\,\mathfrak{T}_a]\,{\cal E}_i, E_a)$.
Both tensors $\mathfrak{T}^\top$ and $\mathfrak{T}^\bot$ obey \eqref{E-nabla-g}.
For example, the \textit{mixed scalar $\mathfrak{T}^\top$-curvature} in RC case~is
\begin{equation}\label{E-SK2}
 {\rm S}^\top_{{\rm mix},\mathfrak{T}} :=
 \sum\nolimits_{a,i} \eps_a\eps_i\,g(\,[\mathfrak{T}^\top_i,\,\mathfrak{T}^\top_a]\,{\cal E}_i, E_a).
\end{equation}

\subsection{Transformation of the mixed scalar curvature}
\label{subsec:prel}

Let $\calf$ be a foliation on a {RC} space $(M,g,\bar\nabla)$ with $\bar\nabla=\nabla+\mathfrak{T}$ and $g_{\,|T\calf}>0$.
Obviously, $({\cal D},{\cal D}^\bot)$-conformal structures \eqref{E-conf-RC}
preserve the decomposition $TM = {\cal D} + {\cal D}^\bot$.
 From \eqref{E-nabla-g} we~get
\begin{eqnarray*}
 {g'}({\mathfrak{T}'}_XY,\,Z) + {g'}({\mathfrak{T}'}_XZ,\,Y)
 \eq u^{2}\,[\,g(\mathfrak{T}_XY,\,Z) + g(\mathfrak{T}_XZ,\,Y)\,] = 0\,.
\end{eqnarray*}
Hence, ${g'}$ is parallel w.r.t.
${\nabla'}+{\mathfrak{T}'}$,
where ${\nabla'}$ is the Levi-Civita connection of ${g'}$.
 Put
\begin{eqnarray}\label{E-ab-tau}
 a_\mathfrak{T} =
 {\rm S}^\top_{{\rm mix},\mathfrak{T}},\quad
 b_\mathfrak{T} = -\sum\nolimits_{i,a}\eps_i\,\eps_a\,
 g(T^{\bot}(\mathfrak{T}_i E_a +\mathfrak{T}_a {\cal E}_i,\,{\cal E}_i),\,E_a)\,.
\end{eqnarray}
Note that $b_\mathfrak{T}=0$ when either ${\cal D}^\bot$ is integrable or
$\bar\nabla$ and $\nabla$ are projectively equivalent.
For a $({\cal D},{\cal D}^\bot)$-conformal structure \eqref{E-conf-RC}
 we have
 $a'_\mathfrak{T} = u^2 a_\mathfrak{T}$
 and
 $b'_\mathfrak{T} = u^{-2} b_\mathfrak{T}$\,.

\begin{proposition}\label{P-main2}
After transformation \eqref{E-conf-RC},
the mixed scalar curvature $\bar{\rm S}\,'_{\rm mix}$ of the {RC} manifold $(M,{g'},{\nabla'}+{\mathfrak{T}'})$
along any $\nabla$-minimal leaf $F$ is
\begin{eqnarray}\label{E-Yam1-u2}
\nonumber
 &&\hskip-6mm \bar{\rm S}\,'_{\rm mix} = \bar{\rm S}_{\rm mix}
 +n(\tr^\top\mathfrak{T})^\bot(u)\,u^{-1}
 +(\tr^\bot\mathfrak{T})^\bot(u)\,u^{-1}
 +n u\,(\tr^\top\mathfrak{T})^\top(u)
  -(u^2{-}1)\,g(\tr^\top\mathfrak{T},\,H^\bot) \\
 &&\hskip-7mm -\,n u^{-1}\Delta^\top u + 2\,u^{-1}H^\bot(u)
   + (u^{-4}{-}1)\<T^\bot,T^\bot\>_{g} - (u^{-2}{-}1)(\<h^\top,h^\top\>_{g} - b_\mathfrak{T}) -(u^{2}{-}1)\,a_\mathfrak{T}\,.
\end{eqnarray}
 If $u=c$ is leafwise constant then
\begin{equation*}
 \bar{\rm S}\,'_{\rm mix} = \bar{\rm S}_{\rm mix}
  + (c^{-4}-1)\<T^\bot,T^\bot\>_{g} - (c^{-2}-1)(\<h^\top,h^\top\>_{g} - b_\mathfrak{T}) -(c^{2}-1)\,a_\mathfrak{T}\,.
\end{equation*}
\end{proposition}

\noindent\textbf{Proof}.
Since ${\cal E}'_i=u^{-1}{\cal E}_i$ is a
${g'}$-orthonormal frame of ${\cal D}^\bot$,
terms of $Q$ in \eqref{E-RC-5} are transformed~as
\begin{eqnarray*}
 &&\hskip-6mm \sum\nolimits_{a,i}\eps_i\eps_a\,{g'}(({\nabla'}_{{\cal E}'_i}{\mathfrak{T}'})_a E_a,{\cal E}'_i)
 =\sum\nolimits_{a,i}\eps_i\eps_a\,g(\nabla_{{\cal E}_i}\mathfrak{T}_aE_a,{\cal E}_i)\\
 && +\,n\,u^{-1}(\tr^\top\mathfrak{T})^\bot(u) +u^{-1}(\tr^\bot\!\mathfrak{T})^\bot(u) +n\,u(\tr^\top\!\mathfrak{T})^\top(u)
 -(u^2{-}1)\,g(\tr^\top\mathfrak{T},H^\bot)\\
 && +\,(u^{-2}{-}1)\big[
 \sum\nolimits_{i,j}\eps_i\eps_j\,g(\mathfrak{T}_j\,{\cal E}_i, T^{\bot}({\cal E}_i,{\cal E}_j))
 -\sum\nolimits_{a,i}\eps_i\eps_a\,g(T^{\bot}({\cal E}_i,\mathfrak{T}_a{\cal E}_i),\,E_a) \big],\\
 &&\hskip-6mm \sum\nolimits_{a,i}\eps_i\eps_a\,{g'}(({\nabla'}_{E_a}\mathfrak{T})_{i'}{\cal E}'_i,E_a)
 = \sum\nolimits_{a,i}\eps_i\eps_a\,g((\nabla_{E_a}\mathfrak{T})_i\,{\cal E}_i,E_a)\\
 &&+\,(u^{-2}-1)\big[
 \sum\nolimits_{i,j}\eps_i\eps_j\,g(\mathfrak{T}_j{\cal E}_i, T^{\bot}({\cal E}_i,{\cal E}_j))
 -\sum\nolimits_{i,a}\eps_i\eps_a\,g(T^{\bot}({\cal E}_i,\mathfrak{T}_i E_a),\,E_a) \big]\,,\\
 &&\hskip-6mm \sum\nolimits_{a,i}\eps_i\eps_a\,{g'}({\cal T'}_a{\cal E}'_i, {\cal T'}_{i'} E_a)
 = \sum\nolimits_{a,i}\eps_i\eps_a\,g(\mathfrak{T}_a{\cal E}_i, \mathfrak{T}_i E_a)
 +(u^{2}-1)\sum\nolimits_{a,i}\eps_i\,\eps_a\,g^\top(\,\mathfrak{T}_a{\cal E}_i, \mathfrak{T}_i E_a)\,,\\
 &&\hskip-6mm \sum\nolimits_{a,i}\eps_i\eps_a\,{g'}(\mathfrak{T}_{i'} {\cal E}'_i, \mathfrak{T}'_a E_a)
 = \sum\nolimits_{a,i}\eps_i\eps_a\,g(\mathfrak{T}_i{\cal E}_i,\mathfrak{T}_a E_a)
 +\,(u^{2}-1)\sum\nolimits_{a,i}\eps_i\eps_a\,g^\top(\,\mathfrak{T}_i{\cal E}_i, \mathfrak{T}_a E_a)\,,
\end{eqnarray*}
where $X(u)=g(X,\,\nabla^\top u)$. Here we used the consequences of \eqref{eqlevicivita},
\begin{eqnarray*}
 {\nabla'}_{{\cal E}_i}\,E_a - \nabla_{{\cal E}_i}\,E_a \eq
 {\nabla'}_{E_a}\,{\cal E}_i -\nabla_{E_a}\,{\cal E}_i
 = E_a(\log u)\,{\cal E}_i - T^{\bot}_i({\cal E}_i) ,\\
 {\nabla'}_{{\cal E}_i}\,{\cal E}_i -\nabla_{{\cal E}_i}\,{\cal E}_i
 \eq -n\,u\nabla^\top u -(n-2)\,u^{-1}\nabla^\bot u +(u^2-1) H^\bot,\\
 {\nabla'}_{E_a}\,E_a - \nabla_{E_a}\,E_a \eq (u^{-2}-1)\,H^\top = \nabla_{E_a}\,E_a\,,
\end{eqnarray*}
where $g(T^{\bot}_Z(X),\,Y) = g(T^\bot(X,Y),\ Z)$ for all $X,Y\in\mathfrak{X}^\bot_M$ and $Z\in\mathfrak{X}^\top_M$.
Thus,
\begin{eqnarray*}
 Q' \eq Q +n\,u^{-1}(\tr^\top\mathfrak{T})^\bot(u) +u^{-1}(\tr^\bot\mathfrak{T})^\bot(u) +n\,u\,(\tr^\top\mathfrak{T})^\top(u) \\
 \minus (u^2-1)\,g(\tr^\top\mathfrak{T},\,H^\bot) -(u^{2}-1)\,a_\mathfrak{T} +(u^{-2}-1)\,b_\mathfrak{T}\,.
\end{eqnarray*}
 From the above, equalities $Q' = \bar{\rm S}\,'_{\rm mix} - {\rm S}\,'_{\rm mix}$,
$Q = \bar{\rm S}_{\rm mix} -{\rm S}_{\rm mix}$ and Lemma~\ref{C-norms1} we obtain \eqref{E-Yam1-u2} on~$F$.
The result for leafwise constant $u$ follows from \eqref{E-Yam1-u2}.
\qed

\begin{lemma}
\label{C-norms1}
Let $\calf$ be a foliation of a pseudo-Riemannian manifold $(M,g)$.
Then, after transformation \eqref{E-conf-RC}$_1$,
the mixed scalar curvature along any minimal leaf $F$ is
\begin{equation}\label{E-Kmixphi-h2}
 {\rm S}\,'_{\rm mix} = {\rm S}_{\rm mix} -n\,u^{-1}\Delta^\top u
 + 2\,u^{-1}H^\bot(u) +(u^{-4}-1)\<T^\bot,T^\bot\>_{g} -(u^{-2}-1)\<h^\top,h^\top\>_{g}\,.
\end{equation}
\end{lemma}

\noindent\textbf{Proof}. This is similar to the proof of Proposition~2.10 in \cite{rz1} for $g_{\,|T\calf}>0$.
Notice that a ${\cal D}^\bot$-conformal change of pseudo-Riemannian metrics preserves
total umbilicity, harmonicity, and total geodesy of foliations, and preserves total umbilicity of $\,{\cal D}^\bot$.
\qed

\smallskip

 One may rewrite \eqref{E-Yam1-u2}~as the second order PDE for the function $u>0$,
\begin{eqnarray}\label{E-Yam1-u2-PDE}
\nonumber
 && -\Delta^\top u +(2/n)\,H^\bot(u) -(\beta^\top+\Phi)\,u  = \Psi_1(x)\,u^{-1} -\Psi_2(x)\,u^{-3} +\Psi_3(x)\,u^{3}\\
 && -(\tr^\top\mathfrak{T})^\bot(u) -(1/n)(\tr^\bot\mathfrak{T})^\bot(u) -u^2\,(\tr^\top\mathfrak{T})^\top(u)
 +((u^3-u)/n)\,g(\tr^\top\mathfrak{T},\,H^\bot),
\end{eqnarray}
where $n\Phi={\rm S}\,'_{\rm mix}$ is the mixed scalar curvature after transformation \eqref{E-conf-RC} and
$\Psi_i$ and $\beta^\top$ are given by
\begin{equation}\label{E-Psi-geom}
 \beta^\top=\Psi_2-\Psi_1-\Psi_3-\frac 1n\,\bar{\rm S}_{\rm mix},\quad
 \Psi_1 = \frac 1n\,(\<h^\top,h^\top\>-b_\mathfrak{T}),\quad
 \Psi_2=\frac 1n\,\<T^\bot,T^\bot\>,\quad
 \Psi_3=\frac 1n\,a_\mathfrak{T}\,.
\end{equation}
Remark that \eqref{E-Yam1-u2-PDE} reduces itself to \eqref{E-Yam1-init} under certain geometric assumptions.

\begin{example}[\bf Flows]\label{Ex-1dim-a}\rm
 If ${\cal D}$
 is spanned by a nonsingular vector field $N$ then $N$ defines a flow (a~one-dimensional foliation).
An example is provided by a circle action $S^1 \times M \to M$ without fixed points.
 A flow of $N$ is \textit{geodesic} if the orbits are geodesics, (i.e., $H^\top=0$),  and a flow is \textit{Riemannian}
 if the metric is bundle-like (i.e., $h^\bot=0$).
 Let $g(N,N)=\eps_{N}\in\{-1,1\}$. In~this case, ${\rm S}_{\rm mix}=\eps_N\Ric_{N}$ and
$\bar{\rm S}_{\,\rm mix}=\eps_N\,\overline\Ric_{N}$ (the Ricci curvature in the $N$-direction).
Thus, for RC case we obtain, see \eqref{E-RC-5},
\begin{equation*}
 \overline\Ric_{N}=\Ric_{N} + \sum\nolimits_{\,i}\eps_i\,\big[g((\nabla_{N}\mathfrak{T})_i{\cal E}_i,N)
 + g((\nabla_{{\cal E}_i}\mathfrak{T})_NN,{\cal E}_i)
 +g(\mathfrak{T}_i N, \mathfrak{T}_N {\cal E}_i) - g(\mathfrak{T}_N N, \mathfrak{T}_i{\cal E}_i)\big].
\end{equation*}
We have $h^\bot=h_{sc}^\bot\,N$, where $h_{sc}^\bot=\eps_N\<h^\bot,N\>$ is the scalar second fundamental form
of~${\cal D}^\bot$. Define the functions $\tau^\bot_i=\tr\,((A^\bot)^{\,i})\ (i\ge0)$,
where the shape operator $A^\bot:{\cal D}^\bot\to {\cal D}^\bot$ obeys $g(A^\bot(X),Y)=h_{sc}^\bot(X,Y)$.
 An easy computation shows that
\begin{eqnarray*}
\nonumber
 H^\top\eq \eps_{N}\nabla_N\,N,\quad h^\top=H^\top\tilde g, \\
 H^\bot\eq \eps_{N}\,\tau^\bot_1 N,\quad \tau^\bot_1= \eps_{N} \tr_{g} h^\bot_{sc},\quad
 \< h^\bot , h^\bot \> = \eps_{N}\,\tau^\bot_2\,.
\end{eqnarray*}
 Let $\{{\cal E}_i\}_{a\le n}$ be a local orthonormal frame on ${\cal D}^\bot$.
Using equalities, see \cite{rz2013},
\begin{eqnarray*}
 \Div N \eq \sum\nolimits_{\,i}\eps_i\,g(\nabla_{{\cal E}_i} N, {\cal E}_i)
 = -g(N,\sum\nolimits_{\,i}\eps_i\,\nabla_{{\cal E}_i}\,{\cal E}_i) = -g(N, H^\bot) = -\tau^\bot_1,\\
 \Div(\tau^\bot_1 N) \eq N(\tau^\bot_1) +\tau^\bot_1\Div N
 = N(\tau^\bot_1)-(\tau^\bot_1)^2,
\end{eqnarray*}
we reduce \eqref{E-PW} to the following:
\begin{equation}\label{E-RicNs1aa-main}
 \Ric_{N} = \Div(\nabla_NN) +(N(\tau^\bot_1) -\tau^\bot_2) +\eps_N\<T^\bot,T^\bot\>.
\end{equation}
Consider transformation \eqref{E-conf-RC} of a {RC} structure
and assume $H^\top=0$ along a compact leaf $F$ (a closed geodesic).
Then, along $F$, the Ricci curvature of $\nabla$ in the $N$-direction is transformed~as
\begin{equation*}
 \Ric'_{N} = \Ric_{N} -n\,u^{-1}N(N(u)) + 2\,u^{-1}\tau^\bot_1 N(u) +(u^{-4}-1)\<T^\bot,T^\bot\>_{g}\,,
\end{equation*}
see \eqref{E-Kmixphi-h2}.
Note that the vector field $\mathfrak{T}_XN$ belongs to $\mathfrak{X}^\bot_M$ for any $X\in\mathfrak{X}_M$.
Hence,
\[
 a_\mathfrak{T}=0,\quad
 b_\mathfrak{T}=-\eps_N\sum\nolimits_{i}\eps_i\,g(T^\bot(\,\mathfrak{T}_iN+\mathfrak{T}_N{\cal E}_i,\,{\cal E}_i),\,N)\,,
\]
where, as usual, $\{{\cal E}_i\}_{a\le n}$ is a local orthonormal frame on ${\cal D}^\bot$.
By the above, see also \eqref{E-Yam1-u2},
\begin{eqnarray}\label{E-Yam1-dim1}
 \overline{\Ric}\,'_{N} \eq \overline\Ric_{N}
 +n\,(\mathfrak{T}_NN)(u)\,u^{-1} +(\tr^\bot\mathfrak{T})^\bot(u)\,u^{-1} - n\,u^{-1}N(N(u))\\
\nonumber
 \plus 2\,u^{-1}\tau^\bot_1 N(u) + (u^{-2}-1)\,b_\mathfrak{T} + (u^{-4}-1)\<T^\bot,T^\bot\>_{g}\,.
\end{eqnarray}
Assuming $\nabla^\bot u=0$ along a compact leaf $F$, we
reduce \eqref{E-Yam1-dim1} along $F$ to a shorter form
\begin{equation*}
 -N(N(u)) +\frac2n\,\tau^\bot_1 N(u) -(\beta^\top+\Phi)\,u  = \Psi_1\,u^{-1} -\Psi_2\,u^{-3}\,,
\end{equation*}
where
 $\beta^\top=\Psi_2-\Psi_1-\frac 1n\,\overline\Ric_{N},\
 \Psi_1 = -\frac 1n\,b_\mathfrak{T},\
 \Psi_2=\frac 1n\,\<T^\bot,T^\bot\>$
 and
 $\Phi=\frac 1n\,\overline{\Ric}\,'_{N}$.
\end{example}

\subsection{Main results}
\label{sec:RC-space}

As promised in the introduction we present two types of solutions to the problem of prescribing $\bar{\rm S}_{\rm mix}$:
 1)~along a compact leaf $F$;  2)~on a closed $M$ under fiber bundle assumption~\eqref{E-2conditions}.

 We~will use spectral parameters (see Section~\ref{sec:app2}) of the elliptic operator \eqref{E-Schred1}, where

\bigskip
\centerline{either $\beta=\beta^\top+\Phi$ when $\Phi\ne{\rm const}$,\ or $\beta=\beta^\top$ when $\Phi={\rm const}$.}

\bigskip\noindent
Here $n\Phi={\rm S}\,'_{\rm mix}$ is the mixed scalar curvature after transformation \eqref{E-conf-RC}.
We can add a real constant to $\beta$ to provide $\beta<0$; then $\mathcal{H}$ is invertible in $L_2(F)$
and $\mathcal{H}^{-1}$ becomes bounded in $L_2(F)$.
If $\Phi={\rm const}$ on $F$ then the ground state $e_0$ does not depend on $\Phi$, see Section~\ref{sec:app2}.

For a positive function $f\in C(F)$ define the quantity
$\delta(f):=(\min\nolimits_{\,F}f)/(\max\nolimits_{\,F}f)\in(0,\,1]$.

\noindent
In case of (\ref{E-2conditions}), the leafwise constants $\lambda_j$ and functions $\{e_j\}$  on~$M$ are smooth.
The least eigenvalue, $\lambda_0$, is simple and obeys on any compact leaf~$F$ the inequalities
\begin{equation}\label{E-beta-lambda}
 -\max\nolimits_{\,F}\beta \le \lambda_0 \le -\min\nolimits_{\,F}\beta\,;
\end{equation}
its eigenfunction $e_0$ (called the ground state) may be chosen positive, see~\cite{rz2013,rz1}.
According to Section~\ref{R-burgers-heat}, we consider three cases:
$\Psi_3>0$, $\Psi_3<0$ and $\Psi_3\equiv0$.
For each case we find $u_*\in C^\infty(F)$, which solves \eqref{E-Yam1-init} under some geometric assumptions.
Assuming $\nabla^\bot u_*=0$ on $F$, we extend the function smoothly onto $M$ and get a required {RC} structure
 $({g'}=g^\top+u_*^2 g^\perp, {\mathfrak{T}'}=u_*^2\,\mathfrak{T}^\top + \mathfrak{T}^\bot)$ on~$M$.
 The~following condition on a leaf $F$ is helpful for the case of $\Psi_3>0$:
\begin{eqnarray}\label{E-Psi-2B}
 27\max\nolimits_{\,F}\Psi_{2}^{2}\cdot\max\nolimits_{\,F}|\Psi_{3}|/\min\nolimits_{\,F}|\Psi_{1}|^{3}
 < \delta^{\,8}(e_0)\,.
\end{eqnarray}
 The~following condition is helpful for the case of $\Psi_3<0$:
\begin{eqnarray}\label{E-Psi-1B}
 \delta(|\Psi_{3}|)\,\delta^{\,2}(e_0)>1/3\,.
\end{eqnarray}
 We also introduce the quantities
\begin{equation}\label{E-ineqL0-2a}
 K_1 = \frac{\psi^{+}_3
 \max\{ 18\psi^{+}_1\psi^{+}_2,\ 4(\psi^{+}_1)^3 \!+27(\psi^{+}_2)^2\psi^{+}_3\}}{4 \psi^{-}_2},\quad
 \psi^{+}_i
 =\max_{\,F}|\Psi_i|,\quad
 \psi^{-}_i
 =\min_{\,F}|\Psi_i|.
\end{equation}
The case of $\Psi_3\equiv0$ has applications for pseudo-Riemannian manifolds, see Theorem~\ref{T-mainB1-R}.

\begin{theorem}\label{T-A}
Let $(M, g, \bar\nabla)$ be a foliated {RC} manifold with the following conditions
along a compact leaf $F$: $g_{\,|TF}>0$, nowhere integrable normal distribution and
\begin{equation}\label{E-RC-lwise}
 H^\bot=0=H^\top,\quad
 (\tr^\top\mathfrak{T})^\top = 0\,.
\end{equation}
Suppose that any of conditions holds on $F$:

\smallskip

1)~$\Psi_3>0,\ \Psi_1>0$ and \eqref{E-Psi-2B};\quad
2)~$\Psi_3<0$ and $\Psi_1<0$;\quad
3)~$\Psi_3=0$ and $\Psi_1>0$.

\smallskip
\noindent
Then for any $\Phi$ obeying, respectively,

\smallskip
1)~$\Phi<-\beta^\top$,\quad
2)~$\Phi>-\beta^\top +1+\delta^{\,-4}(e_0)\sqrt{K_1}$,\quad
3)~$-\beta^\top - \delta^{\,4}(e_0)\,\frac{\min\nolimits_{\,F}\Psi_1^2}{4\max\nolimits_F\Psi_2} < \Phi < -\beta^\top$

\smallskip
\noindent
there exists $u_*\in C^\infty(M)$ such that $M$ with
$({g'}=g^\top+u_*^2 g^\perp, {\mathfrak{T}'}=u_*^2\,\mathfrak{T}^\top + \mathfrak{T}^\bot)$
has $\bar{\rm S}\,'_{\rm mix}=n\Phi$ along $F$;
moreover, $y_2^- \le u_*/e_0 \le y_2^+$ and the set $\{u_{*|F}\}$ of all solutions is compact~in~$C(F)$.
\end{theorem}

\noindent\textbf{Proof}.

1)~By conditions, $\lambda_0>0$; hence, each of bicubic polynomials $P^-(y)$ and $P^+(y)$,
see Section~\ref{sec:attr2}, has three positive roots: $y_3^-<y_2^-<y_1^-$ and $y_3^+<y_2^+<y_1^+$
(which can be expressed by Cardano or trigonometric formulas).
Since \eqref{E-beta-lambda} and \eqref{E-Psi-2B} yield \eqref{E-Psi-2}, we apply Theorem~\ref{thattract2}(i)
 of Section~\ref{sec:attr2}.
Hence, there exists $u_*\in C^\infty(M)$, which solves \eqref{E-Yam1-init} along $F$,
and $y_2^- \le u_*/e_0 \le y_2^+$ holds.

2)~By conditions,
each of bicubic polynomials $P^-(y)$ and $P^+(y)$
has two positive roots: $y_1^-<y_2^-$ and $y_1^+<y_2^+$, see Section~\ref{sec:attr}.
Since \eqref{E-ineqL0-2a} provides \eqref{E-ineqL0},
we apply Theorem~\ref{thattract1}.

3)~For $\Psi_3\equiv0$, the problem amounts to finding a positive solution of the elliptic equation
\begin{equation}\label{E-Yam1-init-3}
 -\Delta^\top u
 -(\beta^\top+\Phi)\,u  = \Psi_1\,u^{-1} -\Psi_2\,u^{-3}\,,
\end{equation}
see \cite{rz2013,rz1}, where
 $\beta^\top=\Psi_2-\Psi_1-\frac1n\,{\rm S}_{\rm mix},\
 \Psi_1 = \frac1n\,(\<h^\top,h^\top\>-b_\mathfrak{T}),\
 \Psi_2=\frac1n\,\<T^\bot,T^\bot\>$.
 For $\Psi_1>0$ and $\Psi_2\ne0$ each of biquadratic polynomials $P^-$ and $P^+$ has two positive roots
$y_1^- < y_2^-$ and $y_1^+ < y_2^+$, see Section~\ref{sec:attr3}.
Conditions and \eqref{E-beta-lambda} yield \eqref{E-lambda0-cond1};
thus, we apply Theorem~\ref{thattract3}.\,\qed

\smallskip
In next corollary we assume that ${\cal D}^\bot$ is integrable.

\begin{corollary}\label{C-A}
Let $(M, g, \bar\nabla)$ be a foliated {RC} manifold with the following conditions
on a compact leaf $F$: $g_{\,|TF}>0$,
integrable normal distribution
and \eqref{E-RC-lwise}.
Suppose that any of conditions holds on $F$:

\smallskip
1)~$\Psi_3>0$
and \eqref{E-Psi-2B};\quad
2)~$\Psi_3<0$
and \eqref{E-Psi-1B};\quad
3)~$a_\mathfrak{T}=0$ and $h^\top=0$.

\smallskip
\noindent
Then for any $\Phi$ obeying, respectively,
\[
 1)~\Phi<-\beta^\top
 -\delta^{\,-2}(e_0)\,\frac{(\max_{\,F}\Psi_1)^{1/2}(3\,\min_{\,F}\Psi_3 +\max_{\,F}\Psi_3)}{\sqrt2\,(3\min_{\,F}\Psi_3-\max_{\,F}\Psi_3)^{1/2}},\quad
 2)~\Phi>-\beta^\top,\quad
 3)~\Phi<-\beta^\top
\]
there exists $u_*\in C^\infty(M)$ unique in $\mathcal{U}_{\,1}$ such that $M$ with
$({g'}=g^\top+u_*^2 g^\perp, {\mathfrak{T}'}=u_*^2\,\mathfrak{T}^\top + \mathfrak{T}^\bot)$
has $\bar{\rm S}\,'_{\rm mix}=n\Phi$ on $F$; moreover, $y_2^- \le u_*/e_0 \le y_2^+$.
\end{corollary}

\noindent\textbf{Proof}. By assumptions,  $\Psi_2=0$ and $b_\mathfrak{T}=0$.

 1) if $\Psi_3>0$ then, in addition, assume $\Psi_1>0$.
 Each of $P^-$ and $P^+$ are reduced to biquadratic polynomials, having two positive roots
$y_{2}^+ <y_{1}^+$ and $y_{2}^- < y_{1}^-$, given in Remark~\ref{R-04} of Section~\ref{sec:attr2}.

 2) if $\Psi_3<0$ then, in addition, assume $h^\top=0$, hence, $\Psi_1=0$.
Each of polynomials $P^-$ and $P^+$ has only one positive root, $y_{2}^-$ and $y_{2}^+$,
given in Remark~\ref{R-03} of Section~\ref{sec:attr}.

3) if $\Psi_3\equiv0$ then, in addition, assume $\Psi_1>0$.
Each of polynomials $P^-$ and $P^+$ has one positive~root,
$y_{2}^-$,
see Remark~\ref{R-05} of Section~\ref{sec:attr3}.
\qed

\begin{remark}\rm
Under stronger geometric conditions along $F$,
see \eqref{E-Psi-1}, \eqref{E-case1-P} and \eqref{E-case1-mu} in Section~\ref{sec:attr2},
the solution $u_{*|F}$, obtained in case 1) of Theorem~\ref{T-A},
is unique in the set $\{\tilde u\in C(F):\,y_3^- < \tilde u/e_0 < y_1^+\}$.
In case 2) of Theorem~\ref{T-A}, if $\Phi>-\beta^\top +1+\delta^{\,-4}(e_0)\sqrt{K_2}$
then the solution $u_{*|F}$ is unique in $\mathcal{U}_{\,1}=\{\tilde u\in C(F):\,\tilde u/e_0 > y_1^-\}$.
In case 3) of Theorem~\ref{T-A}, the solution $u_{*|F}$ is unique in $\mathcal{U}_{\,1}$.
Here
\begin{eqnarray}
\nonumber
 K_2 \eq \frac{\max\big\{36\psi^{+}_1 \psi^{+}_2 \psi^{-}_3 (\psi^{-}_3 +\psi^{+}_3),\
  27 \psi^{+}_3 (\psi^{+}_2)^2 (\psi^{+}_3)^2 + 3 (\psi^{-}_3)^2 +(\psi^{+}_1)^3(\psi^{+}_3 +3 \psi^{-}_3)^2\big\}}
 {8\psi^{+}_2(3 \psi^{-}_3 -\psi^{+}_3)} .
\end{eqnarray}
 In Corollary~\ref{C-A}:
case 2) without assumption \eqref{E-Psi-1B}, and case 1) under weaker assumption
\[
 \Phi<-\beta^\top-2\,\delta^{\,-1}(e_0)(\max\nolimits_{\,F}\Psi_1\cdot\max\nolimits_{\,F}\Psi_3)^{1/2},
\]
we obtain only existence of $u_*\in C^\infty(M)$, but the set $\{u_{*|F}\}$ of all solutions is compact in~$C(F)$.
\end{remark}

\begin{example}\label{R-main0-folB}\rm
Let $a_\mathfrak{T}=0$, $H^\bot=0$, and $T^\bot=0$ hold along a $\nabla$-totally geodesic $F$.
Then $\Psi_1=\Psi_2=\Psi_3=0$, and \eqref{E-Yam1-init} becomes the linear elliptic equation
 $-\Delta^\top u -(\beta^\top+\Phi)\,u =0$\,,
on $F$, where $\beta^\top=-\frac 1n\,\bar{\rm S}_{\rm mix}$.
Suppose that $\bar{\rm S}_{\rm mix}\ne{\rm const}$ and $\Phi={\rm const}$.
Then
${\cal H}\,u_*=\Phi\,u_*$, where $u_*=e_0$ (the~ground state) and $\Phi=\lambda_0$ (the ground level)
for ${\cal H}: u\mapsto -\Delta^\top u
-\beta^\top u $.
\end{example}

In the following theorem, we consider two cases: $\Psi_3<0$ and $\Psi_3\equiv0$,
concerning the sign of the mixed scalar $\mathfrak{T}^\top$-curvature, introduced in \eqref{E-SK2}.
For $\Psi_3>0$, explicit conditions for uniqueness of a solution are difficult; hence, we omit this case.
 In corollary, for integrable normal bundle, we present explicit conditions for uniqueness of a solution
for three cases.

\begin{theorem}\label{T-B}
Let $\calf$ be a
foliation of a closed {RC} manifold $(M, g, \bar\nabla)$
with $g_{\,|T\calf}>0$,
nowhere integrable normal distribution and conditions \eqref{E-2conditions} and:
\begin{equation}\label{E-RC-lwise2}
 H^\bot=0=H^\top,\quad
 \tr^\top\mathfrak{T} = 0\,,\quad
 (\tr^\bot\mathfrak{T})^\bot = 0\,.
\end{equation}
Suppose that any of conditions holds:

\smallskip

1)~$\Psi_3<0,\ \Psi_1<0$ and \eqref{E-Psi-1B};\quad 2)~$\Psi_3=0$, $\Psi_1>0$.

\smallskip
\noindent
Then for any $\Phi$ obeying, respectively,

\smallskip
1)~$\Phi>1-\beta^\top +\delta^{\,-4}(e_0)\sqrt{K_2}$,\quad
2)~$-\beta^\top - \delta^{\,4}(e_0)\min\nolimits_F\Psi_1^2/(4\max\nolimits_F\Psi_2) < \Phi < -\beta^\top$,

\smallskip
\noindent
there exists a leafwise smooth $u_*\in C(M)$, unique in $\mathcal{U}_{\,1}=\{\tilde u\in C(M):\,\tilde u/e_0 > y_1^-\}$,
such that $M$ with
$({g'}=g^\top+u_*^2 g^\perp, {\mathfrak{T}'}=u_*^2\,\mathfrak{T}^\top + \mathfrak{T}^\bot)$
has $\bar{\rm S}\,'_{\rm mix}=n\Phi$; moreover, $y_2^- \le u_*/e_0 \le y_2^+$.
\end{theorem}

\noindent\textbf{Proof}. 1)~As in the proof of case 2) of Theorem~\ref{T-A},
we apply Theorem~\ref{thattract1} of Section~\ref{sec:attr}, and then Theorem~\ref{contdeponqeigfunct}.
2)~We apply Theorem~\ref{thattract3} of Section~\ref{sec:attr3}, and then Theorem~\ref{contdeponqeigfunct}.
\qed

\begin{corollary}\label{C-B}
Let $\calf$ be a
foliation of a closed {RC} manifold $(M, g, \bar\nabla)$
with conditions $g_{\,|T\calf}>0$,
\eqref{E-2conditions} and \eqref{E-RC-lwise2}.
Suppose that any of conditions holds:

\smallskip
1)~$\Psi_3>0$, $h^\top=0$ and \eqref{E-Psi-2B};\quad
2)~$\Psi_3<0$, $h^\top\ne0$ and \eqref{E-Psi-1B};\quad
3)~$\Psi_3=0$ and $h^\top\ne0$.

\smallskip
\noindent
Then for any $\Phi$ obeying, respectively,

\[
 1)~\Phi<-\beta^\top -\delta^{\,-2}(e_0)\,\frac{(\max_{\,F}\Psi_1)^{1/2}(3\,\min_{\,F}\Psi_3 +\max_{\,F}\Psi_3)}{\sqrt2\,(3\min_{\,F}\Psi_3-\max_{\,F}\Psi_3)^{1/2}},\quad
 2)~\Phi>-\beta^\top,\quad
 3)~\Phi<-\beta^\top
\]

\noindent
there exists a leafwise smooth $u_*\in C^\infty(M)$, unique in $\mathcal{U}_{\,1}$, such that $M$ with
$({g'}=g^\top+u_*^2 g^\perp, {\mathfrak{T}'}=u_*^2\,\mathfrak{T}^\top + \mathfrak{T}^\bot)$
has leafwise constant $\bar{\rm S}\,'_{\rm mix}=n\Phi$; moreover, $y_2^- \le u_*/e_0 \le y_2^+$.
\end{corollary}

\begin{remark}\rm
If $H^\top=0$ and $h^\top\ne0$ on $M$ (see Corollary~\ref{C-A} cases 1 and 3,
then the foliation is harmonic and nowhere totally geodesic.
There exist foliations of any codimension $>1$ with harmonic,
nowhere totally geode\-sic leaves on (compact) Lie groups with left-invariant
metrics, see \cite{ty84}; furthermore, the metric can be chosen to be bundle-like.
Such foliations have leafwise constant mixed scalar~curvature.
\end{remark}

The above has consequences for foliated pseudo-Riemannian manifolds.

\begin{theorem}\label{T-mainB1-R}
Let $(M, g)$ be a foliated pseudo-Riemannian manifold with conditions $T^\bot\ne0$, $H^\bot=0=H^\top$ and $h^\top\ne0$
along a~compact leaf $F$ with $g_{\,|TF}>0$. Then for any $\Phi$ obeying
\begin{equation}\label{E-lambdaPhi-R}
 -\beta^\top - \delta^{\,4}(e_0)\min\nolimits_{\,F}\<h^\top,h^\top\>^2/(4\,n\max\nolimits_{\,F}
 \<T^\bot,T^\bot\>) < \Phi < -\beta^\top,
\end{equation}
there exists a leafwise smooth $u_*\in C(M)$, unique in~${\cal U}_{\,1}=\{\tilde u\in C(F):\ \tilde u/e_0 > y_1^-\}$
such that $(M,\,{g'}=g^\top+u_*^2 g^\perp)$ has ${\rm S}\,'_{\rm mix}=n\Phi$ on $F$;
moreover, $y_2^- \le u_*/e_0 \le y_2^+$.
\end{theorem}

\noindent\textbf{Proof}.
The problem means to find a positive solution $u_*\in C^\infty(M)$ of elliptic equation on~$F$:
\begin{equation}\label{E-Yam1-init-3R}
 -\Delta^\top u -(\beta^\top+\Phi)\,u = \Psi_1\,u^{-1} -\Psi_2\,u^{-3}\,,
\end{equation}
see~\eqref{E-Yam1-init-3}, where
 $\beta^\top=\Psi_2-\Psi_1-\frac1n\,{\rm S}_{\rm mix},\
 \Psi_1 = \frac1n\,\<h^\top,h^\top\>,\
 \Psi_2=\frac1n\,\<T^\bot,T^\bot\>$.
\newline
In conditions, each of biquadratic polynomials $P^-$, $P^+$
has two positive roots $y_2^- < y_1^-$ and $y_2^+ < y_1^+$, see Section~\ref{sec:attr3}.
The case of $h^\top_{\,|F}=0$ is not applicable, see paragraph~(c$_1$) in Section~\ref{R-burgers-heat}.
Thus, the mixed scalar curvature of the metric ${g'}=g^\top+u_*^2 g^\perp$ along $F$ is~$n\Phi$.
\qed

\smallskip

If $T^\bot=0$ on $F$ then each polynomial $P^-$ and $P^+$ has one positive~root,
$y_{1}^-$ and $y_{1}^-$, see Remark~\ref{R-05} in Section~\ref{sec:attr3}.
If~${\rm S}_{\rm mix}>0$ then there are no compact $\nabla$-harmonic leaves, see \cite[Theorem~2]{wa1}.

\begin{corollary}\label{T-main0-Riemfol-R}
Let $T^\bot=0$, $H^\bot=0=H^\top$, and $h^\top\ne0$ on $F$.
Then for any $\Phi$ obeying $\Phi<-\beta^\top$ there exists a leafwise smooth $u_*\in C(M)$, unique in~${\cal U}_{\,1}$,
such that $(M,\,{g'}=g^\top+u_*^2 g^\perp)$ has ${\rm S}\,'_{\rm mix}=n\Phi$ on $F$; moreover, $y_2^- \le u_*/e_0 \le y_2^+$.
\end{corollary}

Similar results (to Theorem~\ref{T-mainB1-R}) for a closed manifold $M$ with \eqref{E-2conditions}, extend our results in \cite{rz1}.

\section{The nonlinear heat equation}
\label{sec:app2}

Let $(F, g)$ be a closed $p$-dimensional Riemannian manifold,
$H^l(F)$ the Hilbert space of differentiable by Sobolev real functions on $F$
with the inner product $(\,\cdot,\cdot\,)_{l}$ and the norm $\|\cdot\|_l$,
e.g. $H^0(F)=L_2(F)$.
If $B$ and $C$ are Banach spaces with norms $\|\cdot\|_B$ and $\|\cdot\|_C$, denote by $\mathcal{B}^r(B,C)$ the Banach space
of all bounded $r$-linear operators $A:\,\prod_{i=1}^r B\rightarrow C$ with the norm
 $\|A\|_{\mathcal{B}^r(B,C)}=\sup_{v_1,\dots, v_r\in B\setminus{0}}\frac{\|A(v_1,\dots,v_r)\|_C}{\|v_1\|_B\cdot\ldots\cdot\|v_r\|_B}$.
If~$r=1$, we shall write $\mathcal{B}(B,C)$ and  $A(\cdot)$,
and if $B=C$ we shall write $\mathcal{B}^r(B)$ and $\mathcal{B}(B)$, respectively.
Denote by $\|\cdot\|_{C^k}$ the norm in the Banach space $C^k(F)\ (k\ge1)$,
and $\|\cdot\|_C$ for $k=0$.
In~coordinates $(x_1,\dots, x_p)$ on~$F$, we have
$\|f\|_{C^k}=\max\nolimits_{\,F}\max\nolimits_{\,|\alpha|\le k}|d^\alpha f|$, where $\alpha\ge0$ is the multi-index of order
$|\alpha|=\sum_{i}\alpha_i$ and $d^\alpha$ is the partial derivative.
Consider the nonlinear elliptic equation, see \eqref{E-Yam1-init},
\begin{equation}\label{Cauchy-stat1}
 -\Delta u -\beta\,u = \Psi_1(x)\,u^{-1} -\Psi_2(x)\,u^{-3} +\Psi_3(x)\,u^{3}\,,
\end{equation}
where $\Psi_i$ and $\beta$ are arbitrary smooth functions on $F$, and $\Psi_2\ge0$.
If $\Psi_i\ (i=1,2,3)$ are real constants then
\eqref{E-Yam1-init} belongs to reaction-diffusion equations, which are well understood.
The lhs of \eqref{Cauchy-stat1} is the Schr\"{o}dinger operator $\mathcal{H}:=-\Delta -\beta\id$ with domain of definition~$H^2(F)$.
 One can add a real constant to $\beta$ such that $\mathcal{H}$ becomes invertible in
$L_2$ (e.g., $\lambda_0>0$) and $\mathcal{H}^{-1}$ is bounded in $L_2(F)$.
 Recall the Elliptic regularity Theorem, see \cite{aub}:

\smallskip
\textit{If $\,0\notin\sigma(\mathcal{H})$ then
$\mathcal{H}^{-1}: H^k(F)\rightarrow H^{k+2}(F)$
for any integer~$k\ge0$}.

\smallskip
\noindent
For $k=0$, we have $\mathcal{H}^{-1}: L_2(F)\rightarrow H^2(F)$, and the embedding of
$H^2(F)$ into $L_2(F)$ is continuous and compact;
hence, the operator $\mathcal{H}^{-1}:\,L_2(F)\rightarrow L_2(F)$ is compact.
Thus, the spectrum of $\mathcal{H}$ is discrete,
the least eigenvalue $\lambda_0$ of $\mathcal{H}$
is simple, its eigenfunction $e_0(x)$ (called the \textit{ground state})
can be chosen positive, see~\cite{rz1}.
 Since $(\beta(x)u,\,u)_0\ge \beta^-(u,\,u)_0$, where $\beta^-=\min_{\,F}\beta$, we have
\begin{equation*}
 (\mathcal{H}u,\,u)_0=\!\int_F(|\nabla u(x)|^2-\beta(x)|u(x)|^2)\,{\rm d}x
 \le\!\int_F(|\nabla u(x)|^2-\beta^-|u(x)|^2)\,{\rm d}x=(-\Delta u-\beta^- u, u)_0
\end{equation*}
for any $u\in\rm{Dom}(\mathcal{H})$.
Since $\beta^-$ is the maximal eigenvalue of the linear operator $\Delta+\beta^-\id$,
by the variational principle for eigenvalues, we obtain $\lambda_0\le-\beta^-$, see~\eqref{E-beta-lambda}.
Similarly, $\lambda_0\ge-\max_{\,F}\beta$.
 To solve \eqref{Cauchy-stat1}, we look for attractor of the Cauchy's problem for the
 {heat equation},
\begin{eqnarray}\label{Cauchy}
 \partial_t u = \Delta u +\beta\,u +\Psi_1(x)\,u^{-1}-\Psi_2(x)\,u^{-3} +\Psi_3(x)\,u^{3},\quad
 \quad u(x,0) = u_0(x)>0\,.
\end{eqnarray}

Let $\mathcal{C}_{t}=F\times[0,t),\ (0<t\le\infty)$, be cylinder with the base $F$.
By~\cite[Theorem~4.51]{aub}, (\ref{Cauchy}) has a~unique smooth solution in $\mathcal{C}_{t_0}$ for some $t_0>0$.
Substituting $u=e_0 w$ into \eqref{Cauchy} and using
$\Delta(e_0 w)=e_0\Delta w+w\Delta e_0+\<2\,\nabla e_0,\nabla w\>$
and
$\Delta e_0 +\beta e_0=-\lambda_0e_0$, yields the Cauchy's problem
\begin{equation}\label{Cauchw1R-C}
\partial_t w =\Delta w+\,\<2\,\nabla\log e_0,\,\nabla w\> +f(w,\,\cdot\,),\quad
 w(\cdot\,,0) ={u_0}/{e_0}>0
\end{equation}
 for $w(x,t)$, where
 $f(w,\,\cdot\,)=-\lambda_0\,w+(\Psi_1 e_0^{-2})\,{w^{-1}}-(\Psi_2 e_0^{-4})\,w^{-3} +(\Psi_3 e_0^2)\,{w^{3}}$.
 From (\ref{Cauchw1R-C}) we obtain the differential inequalities
\begin{equation}\label{difineqR-C}
 \phi_-(w) \le \partial_t w - \Delta w - \<2\nabla\log e_0,\nabla w\> \le \phi_+(w),
\end{equation}
where the functions $\phi_-$ and $\phi_-$ are defined for each case separately.

Define the parallelepiped $\mathcal{P}=[\Psi_1^-,\Psi_1^+]\times[\Psi_2^-,\Psi_2^+]\times[\Psi_3^-,\Psi_3^+]\subset\RR^3_+$,
where
\begin{eqnarray*}
 \Psi_k^+ \eq\max\nolimits_{\,F}\,(|\Psi_k|\,e_0^{-2k}),\quad
 \Psi_k^-=\min\nolimits_{\,F}\,(|\Psi_k|\,e_0^{-2k})\quad (k=1,2),\\
 \Psi_3^+ \eq\max\nolimits_{\,F}\,(|\Psi_3|\,e_0^{2}),\qquad
 \Psi_3^-=\min\nolimits_{\,F}\,(|\Psi_3|\,e_0^{2}).
\end{eqnarray*}
Then $\mathcal{P}_0=\{(\Psi_1(x),\Psi_2(x),\Psi_3(x)):\,x\in F\}$ is a closed subset of $\mathcal{P}$.
We shall use the following.

\begin{proposition}[Scalar maximum principle, see
\cite{ah}]\label{P-weak-max}
Let $X_t$ and $g_t$ be smooth families of vector fields
and metrics on a closed manifold $F$, and $f\in C^\infty(\RR\times[0,T))$.
Suppose that $u:F\times[0,T)\to\RR$ is a $C^\infty$ supersolution to
 $\dt\,u \ge\Delta_t\,u - X_t(u) +f(u,t)$,
and $y:[0,T]\to\RR$ solves the Cauchy's problem for ODEs: $y\,' = f(y(t),t),\ y(0)=c$.
If~$u(\cdot\,, 0)\ge c$ then $u(\cdot\,,t)\ge y(t)$ for~$t\in[0,T)$.
\end{proposition}


Let $F\times\RR^n$ be the product with a compact leaf $F$, and
$g(\cdot,q)$ a leafwise Riemannian metric (i.e., on
$F_q=F\times\{q\}$ for $q\in\RR^n$) such that the volume form of the
leaves $d\,{\rm vol}_F=|g|^{1/2}\,{\rm d} x$ depends on $x\in F$
only (e.g., the leaves are minimal submanifolds, see~Section~\ref{sec:RC-space}).
This assumption simplifies arguments used in the proof of Theorem~\ref{contdeponqeigfunct} below
(we consider products $\mathbb{B}=L_2\times\RR^n$ and
$\mathbb{B}_k=H^k\times\RR^n$ instead of infinite-dimensional vector
bundles over~$\RR^n$), on the other hand, it is sufficient for proving the geometric results.
 The leafwise Laplacian in a local chart $(U,x)$ on $(F,g)$ is written as
$\Delta\,u = \nabla_i(g^{ij}\,\nabla_j\,u) =
|g|^{-1/2}\partial_{i}(|g|^{1/2} g^{ij} \partial_j\,u)$, see \cite{aub}.
This defines a self-adjoint elliptic operator
$-\Delta_q$, where $q\in\RR^n$ is a parameter and $\Delta_0=\Delta$,
\begin{equation*}
 \Delta_q = g^{ij}(x,q) \partial^2_{ij} + b^j(x,q)\partial_j\,.
\end{equation*}
Here $b^j=|g|^{-\frac12}\partial_{i}(|g|^{\frac12} g^{ij})$ are
smooth functions in $U\times\RR^n$. The Schr\"odinger~operator
\begin{equation}\label{E-Hq}
 \mathcal{H}_q=-\Delta_q-\beta(x,q)\id,\quad q\in\RR^n
\end{equation}
acts  in the Hilbert space $L_2$ with the domain $H^2$ and it is self-adjoint.

\begin{theorem}[see \cite{rz1}]\label{contdeponqeigfunct}
Let $\lambda(q)$ be the least eigenvalue  of $\mathcal{H}_q\ (q\in\RR^n)$.
If $\beta\in C^\infty(F\times\RR^n)$ then $\lambda\in C^\infty(\RR^n)$ and
there exists a unique function $e\in C^\infty(F\times\RR^n)$ such that
$e(\cdot,q)$ is a positive eigenfunction of $\mathcal{H}_q$
related to $\lambda(q)$ with $\Vert e(\cdot,q)\Vert_{L_2}=1$.
\end{theorem}

\begin{theorem}[see \cite{rz1}]\label{prsmoothparstat}
Let $f\in C^\infty(D\times\RR^n)$ and $u_*(x)\in\mathrm{Int}\,G$ be a smooth solution of
\begin{equation}\label{nonlinstatpar1a}
 \Delta_q\,u + f(u,x,q)=0,
\end{equation}
with $q=0$ such that $\lambda=0$ is not an eigenvalue of
${\mathcal H}=-\Delta-\partial_u f(u_*(x),x,0)$ on $L_2$ with domain in $H^2$.
Then for any integers $k\ge 0$ and $l\ge 1$ and $\alpha\in(0,1)$ there are open neighborhoods
$U_*\subseteq C^{k+2,\alpha}$ of $u_*$ and $V_0\subseteq\RR^n$ of $\,0$
such that for any $q\in V_0$
(\ref{nonlinstatpar1a}) has in $U_*$ a unique solution $u(x,q)$,
in particular, $u_*(x){=}u(x,0)$
such that
$q\rightarrow u(\cdot,q)$ belongs to~$C^l(V_0,U_*)$.
\end{theorem}

\subsection{Comparison ODE with constant coefficients}
\label{R-burgers-heat}

Let $\beta$ and  $\Psi_i\ (i=1,2,3)$ be real constants with $\Psi_2>0$ (the case of $\Psi_2<0$ is studied similarly).
Then reaction-diffusion equation \eqref{Cauchy} can be compared with the ordinary differential equation with constant coefficients,
whose solutions can be written explicitly and easily investigated.
Namely, leafwise constant solutions of \eqref{Cauchy} obey the Cauchy's problem for~ODE:
\begin{equation}\label{zermod}
 y\,'=f(y) = P(y^2)/y^3,\qquad y(0)=y_0>0
\end{equation}
with the polynomial $P(z) = \Psi_3\,z^{3} +\beta\,z^{2} +\Psi_1\,z -\Psi_2$.
Recall that $P(z)$ (when $\Psi_3\ne0$)
has three different real roots if and only if the discriminant $D_P=-{\rm Res}(P,P^{\,\prime})/\Psi_3$ is positive,
where ${\rm Res}(P,P^{\,\prime})$ is the resultant of two polynomials.
Consequently, $P(z)$ has one real root if and only if $D_P<0$.
Remark that $D_P = 4\Psi_2\beta^3+\Psi_1^2\beta^2-18\,\Psi_1\Psi_2\Psi_3\,\beta - (4\Psi_1^3 + 27\,\Psi_2^2\Psi_3)\Psi_3$
is a cubic polynomial in~$\beta$,
which is positive when $\beta\to\infty$. By Maclaurin method in what follows, one may take
 $\beta > 1+\big(\max\{18\Psi_2|\Psi_1\Psi_3|,\ |\,4\Psi_1^3\Psi_3 + 27\,\Psi_2^2\Psi_3^2\,|\,\}
 /(4\Psi_2)\big)^{1/2}$.

\noindent
\textbf{Maclaurin method}. Suppose that the first $m$ leading coefficients of the real polynomial
$P_n(t)=a_0t^n+a_1t^{n-1}+\ldots + a_{n-1}\,t+a_n$ are nonnegative,
i.e., $a_0 > 0,\,a_{1}\ge0,\ldots,\, a_{m-1}\ge0$,
and the next coefficient is negative, $a_{m} < 0$. Then $1 + (B/a_0)^{1/m}$
is an upper bound for the positive roots of this polynomial, where $B$ is the largest of the
absolute values of negative coefficients of $P_n(t)$.
 Note that $P_n(t)>0$ for all $t\in[0,1]$ (so, $a_n>0$) if
\begin{equation}\label{E-a3}
 a_n > \sum\nolimits_{\,0\le i<n}|a_i|,\quad{\rm for\ all}\quad a_i<0.
\end{equation}

 We look for stable stationary solutions of \eqref{zermod}, i.e., roots of $P(y^2)$.
 If there exists a real root $\tilde y>0$ such that $f'(\tilde y)<0$ then
 $y=\tilde y$ is a one-point attractor for the semigroup associated to (\ref{zermod}).
 The basin of attractor is determined by other two positive roots of which surround~$\tilde y$.

\begin{figure}
\begin{tabular*}{\textwidth}{@{}l@{}l@{}r@{}}
\hskip15mm\subfigure
{\includegraphics[width=0.31\textwidth]{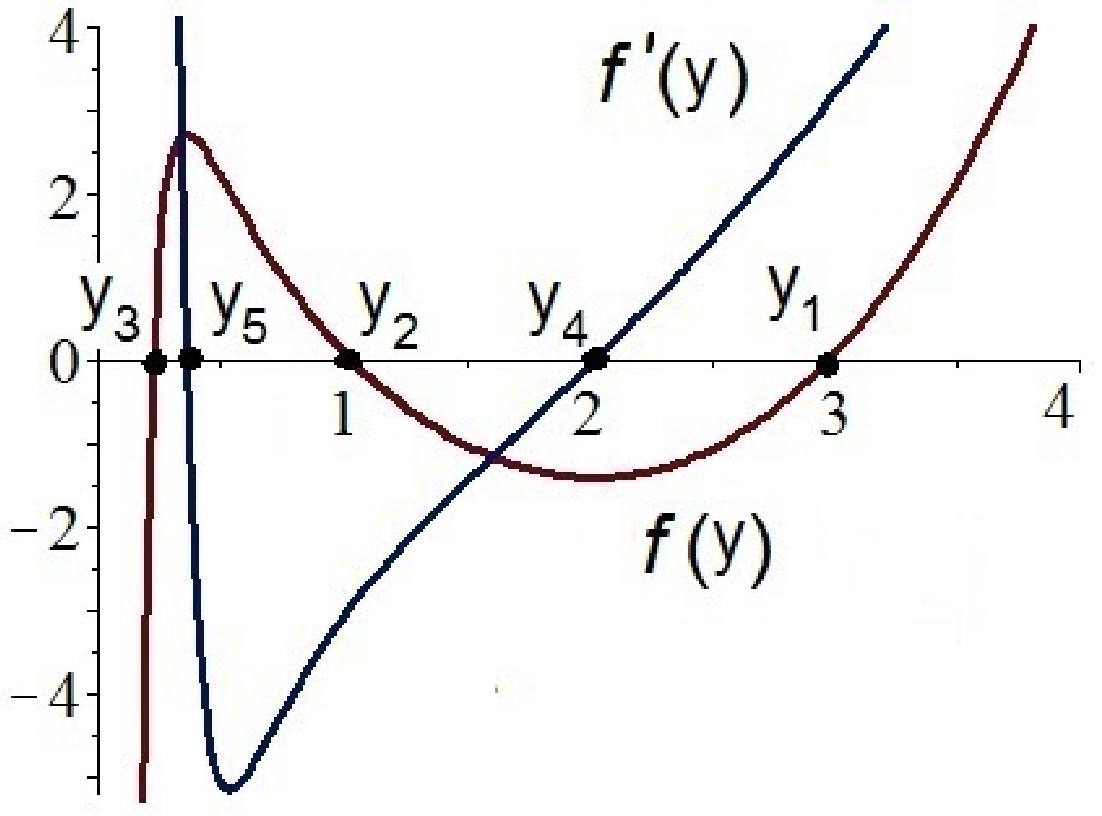}
}
\hskip15mm\subfigure
{\includegraphics[width=0.3\textwidth]{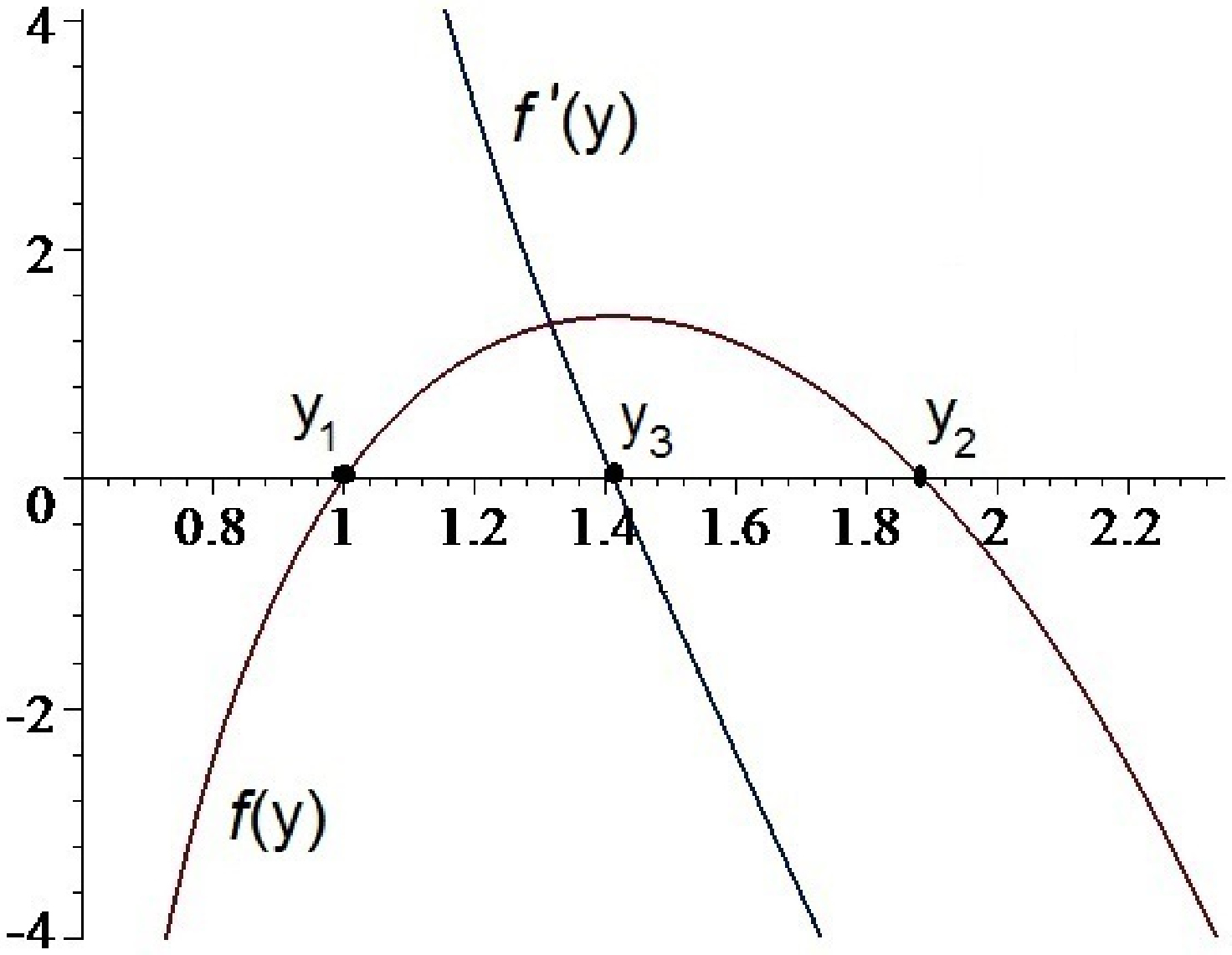}}
\end{tabular*}
\vskip-2mm
\caption{\small  Section~\ref{R-burgers-heat}, cases (a) and (b): graphs of functions $f$ and $f^\prime$.
\newline
(a)~$\Psi_1>0,\ \Psi_2>0,\ \Psi_3>0,\ D_P>0
$;\quad
(b)~$\Psi_1<0,\ \Psi_2>0,\ \Psi_3<0$\,.}
\label{fig:1}
\end{figure}

\textbf{(a)} Let $\Psi_3>0$. Thus, $P(z)$ has the properties:
 $P(0)=-\Psi_2<0,\ P(\infty)=\infty$ and $P(-\infty)=-\infty$.
The condition $D_P>0$ and the fact that both roots of the quadratic polynomial
$P^\prime(z)$ are
positive imply that all three roots of $P(z)$ are
positive, $z_3<z_2<z_1$.
Indeed, $P(z)$ increases in the semi-axis $(-\infty,\,0]$;
hence, in view of $P(0)<0$, it has no
negative roots.
Note that if $\beta^2-3\Psi_1\Psi_3>0$, $\beta<0$ and $\Psi_1>0$,
then both roots $z_4>z_5$ of $P^{\,\prime}(z)$ are
positive. Thus, conditions
\begin{equation*}
 \Psi_1>0,\quad \Psi_2>0,\quad\Psi_3>0,\quad\beta<0,\quad D_P>0
\end{equation*}
guarantee existence of a stable stationary solution $y_2=z_2^2>0$
(and unstable solutions $y_1=z_1^2>0$ and $y_3=z_3^2>0$)
of \eqref{zermod}, see Fig.~\ref{fig:1}(a).
Hence, $f'(y)$ has two
positive roots, $y_5<y_4$.
We~conclude that the basin of a single-point attractor $y=y_2$ for the semigroup of operators of \eqref{zermod}
is the (invariant) set of continuous functions $y(t)$, whose values belong to $(y_3,\, y_1)$.

\textbf{(b)} Let $\Psi_3<0$.
The cubic polynomial $P(z)$ has the properties: $P(0)=-\Psi_2<0,\ P(\infty)=-\infty$.
Its maximal real root $z_2$ is an attractor for
the heat equation. Note that the condition $D_P>0$ and the fact that
the maximal root $z_0$ of the derivative  $P^\prime$ is
positive imply that $z_2>0$ (and $z_1>0$ is the minimal positive root of $P$).
Indeed, otherwise all the roots of $P$ are
negative, hence both roots of $P_3^\prime$
are negative in contradiction with the assumption.
If $\beta>0$, $\Psi_1<0$ and $\beta^2-3\,\Psi_1\Psi_3>0$ (the discriminant of $P'$ is positive)
then both roots of
$P'  (z)=3\,\Psi_3 z^{2} + 2\,\beta z + \Psi_1$ are real and the maximal root
$z_0=\frac{(\beta^2-3\Psi_1\Psi_3)^{1/2}-\beta}{3\Psi_3}$ is positive.
In~view of
\[
 27\,\Psi_3^2\cdot D_P = 4(\beta^2-3\,\Psi_1\Psi_3)^{3} - (27\,\Psi_2\Psi_3^2+9\,\beta\,\Psi_1\Psi_3-2\,\beta^3)^2,
\]
the~condition $D_P>0$ implies the inequality $\beta^2-3\,\Psi_1\Psi_3>0$. Thus, the conditions
\begin{equation*}
 \Psi_1<0,\quad \Psi_2>0,\quad\Psi_3<0,\quad\beta>0,\quad D_P>0
\end{equation*}
guarantee existence of a stable stationary solution $y_2=z_2^2>0$
(and existence of unstable stationary solution $y_1=z_1^2>0$) of \eqref{zermod}, see Fig.~\ref{fig:1}(b).
Note that $f''(y)=(6\,\Psi_3\,y^{6} +2\,\Psi_1\,y^{2} -12\,\Psi_2)/y^{5}$ is negative for $y>0$.
Hence, $f(y)$ is concave for $y>0$, and $f'(y)$ is monotone decreasing
(with $f'(0+)=\infty$ and $f'(\infty)=-\infty$) and has one
positive root.
We~conclude that the basin of a single-point attractor $y=y_2$ for the semigroup of \eqref{zermod}
is the (invariant) set of continuous functions~$y(t)>y_1$.

\textbf{(c)} Let $\Psi_3=0$.
Then $P(z) = \beta\,z^{2} +\Psi_1\,z -\Psi_2$, see \eqref{zermod}.
A positive root $\tilde z$ of $P(z)$ corresponds to a stationary solution $\tilde y=\sqrt{\tilde z}$ of (\ref{zermod});
moreover, if $P\,'(\tilde z)<0$ then $\tilde y$ is a single-point attractor.

\begin{figure}
\begin{tabular*}{\textwidth}{@{}l@{}l@{}r@{}}
\hskip15mm\subfigure
{\includegraphics[width=0.36\textwidth]{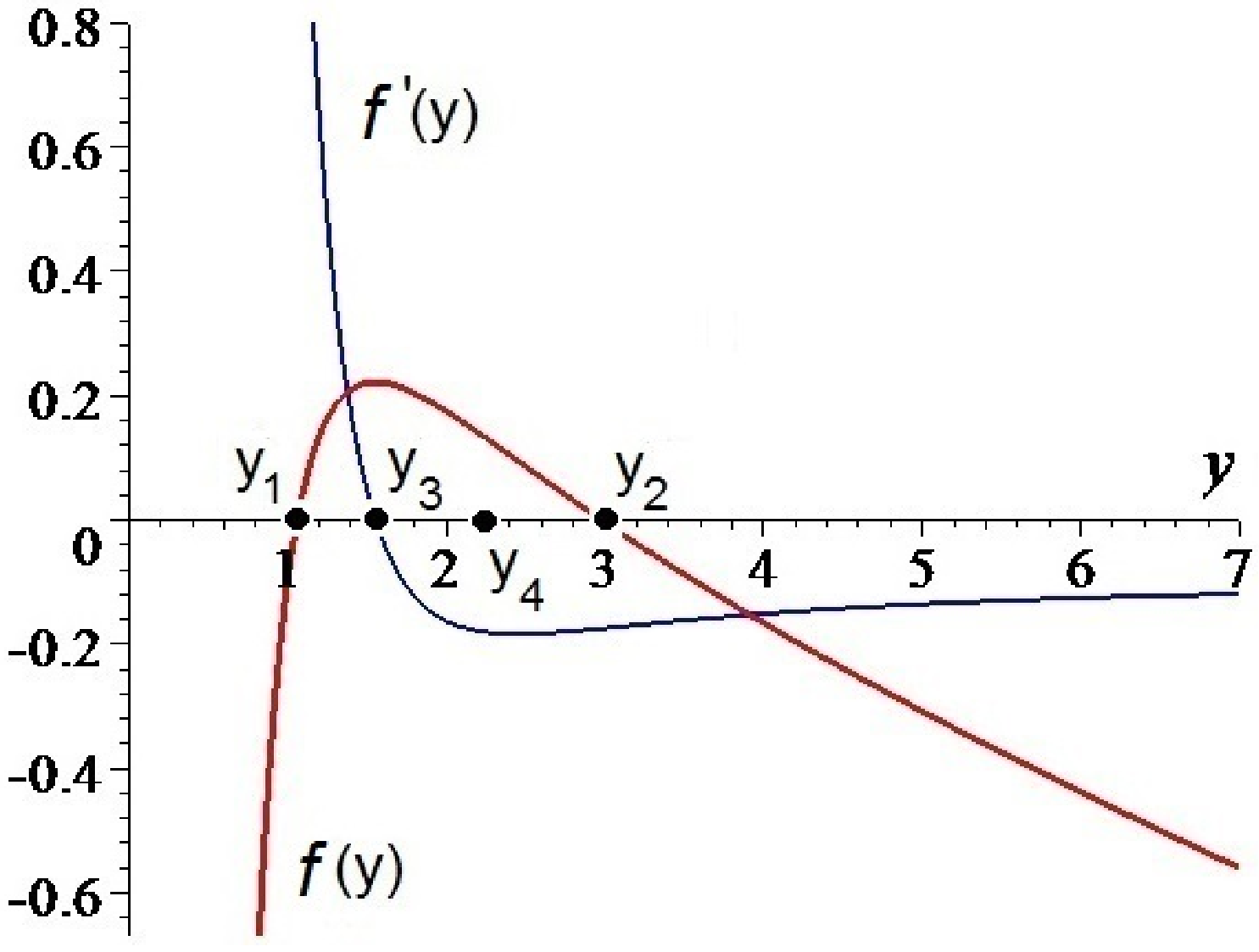}}
\hskip10mm\subfigure
{\includegraphics[width=0.28\textwidth]{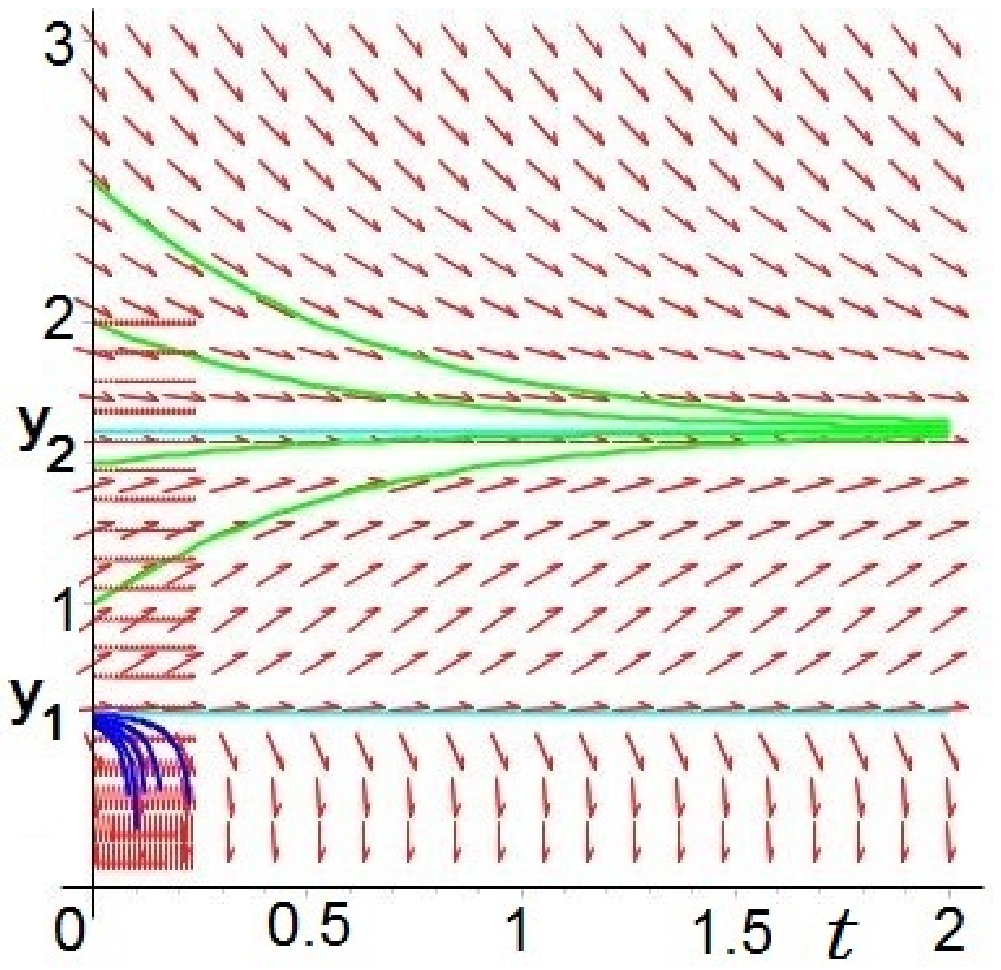}}
\end{tabular*}
\vskip-2mm
\caption{\small Section~\ref{R-burgers-heat}, paragraph~(c$_1$): $P(y^2)=\beta\,y^4 +\Psi_1\,y^{2} -\Psi_2$ with $\beta<0$ and $4\,|\beta|\,\Psi_2<\Psi_1^2$.
\newline
(a)~Graphs of $f$ and $f^\prime$ for $\Psi_3=0$, $\Psi_1>0$, $\Psi_2>0$.
\qquad
(b)~$y_1$~is unstable, and $y_2$ is stable.}
\label{fig:2}
\end{figure}

(c$_1$) Assume $\beta<0$. We have $P(0)=-\Psi_2<0$ and $P(\infty)=-\infty$.
Thus, $P(z)$ has real roots if and only if $P(z_0)>0$, where $z_0=-\Psi_1/\beta$ is a root of $P'(z)=0$.
In our case, the inequality $P(z_0)>0$ is valid when $-(\Psi_1)^2/(4\Psi_2)<\beta<0$.
Maximal root $y_{2}=\big(\frac{\Psi_1+(\Psi_1^2-4\,|\beta|\,\Psi_2)^{1/2}}{2 |\beta|}\big)^{1/2}$
of $f(y)=0$ is asymptotically stable, but the second (minimal) root $y_1$ is unstable;
moreover, $f'(y)$ has a unique positive root $y_3$, and $f'(y)$ takes minimum at $y_4$, see Fig.~\ref{fig:2}.
If $-4\,\beta\Psi_2=\Psi_1^2$ then (\ref{zermod}) has one positive stationary solution,
and has no stationary soluti\-ons if $-4\,\beta\Psi_2>\Psi_1^2$.

(c$_2$) Let $\beta>0$.  We have $P(0)=-\Psi_2<0$ and $P(\infty)=\infty$.
Thus, $P(z)$ has one positive root $z_2$, which corresponds to unstable stationary solution
of (\ref{zermod}), because $P'(z_2)>0$. One may show that for $\beta=0$, (\ref{zermod})
has a unique positive stationary solution, which is unstable.

(c$_3$) Let $\Psi_2=0$, then $f(y)=\beta\,y+\Psi_1 y^{-1}$.
If $\beta\ge 0$ then there are no positive stationary solutions.
If $\beta<0$ and $\Psi_1>0$ then $f(y)=0$ has one positive root $y_2=({\Psi_1/|\beta|})^{1/2}$.
The solution $y_1$ is stable (attractor) because $f^\prime(y_2)<0$.

\subsection{A fixed point of a one-parametric semigroup}

\begin{definition}\label{dfeqicont}\rm
Let $(X,\,d)$ and $(Y,\,d^{\,\prime})$ be metric spaces.
A family of mappings $\{f_\alpha:\,X\rightarrow Y\}_{\alpha\in\mathcal{A}}$
is called \textit{equicontinous},
if for any $\varepsilon>0$ there is $\delta>0$ such that for any pair of
points $x_1,x_2\in X$ satisfying the condition $d(x_1,x_2)<\delta$
the inequality $d^{\,\prime}(f_\alpha(x_1),f_\alpha(x_2))<\varepsilon$ holds for any $\alpha\in\mathcal{A}$.
A~family of mappings $\{f_t:\,X\rightarrow Y\}_{t\in[a,\,b]}$
is called \textit{continuous by $t$ uniformly with respect to} $x\in X$, if the family of
mappings $\{\phi_x(t):=f_t(x):\,[a,\,b]\rightarrow Y\}_{x\in X}$ is equicontinuous.
\end{definition}

The following lemma extends the Arzela-Ascoli Theorem.

\begin{lemma}\label{genArz}
Let $\mathcal{U}=\{u_\alpha\}_{\alpha\in\mathcal{A}}$ be a family
of functions defined in a closed interval $[a,\,b]$, with values
in a Banach space $E$ and having the properties:

$(a)$ for any $t\in[a,\,b]$ the set
$U_t=\{u_\alpha(t)\}_{\alpha\in\mathcal{A}}$ is precompact in $E$;

$(b)$ the family $\mathcal{U}$ is equicontinuous.

\noindent
Then the family $\mathcal{U}$ is precompact in $C([a,\,b],\,E)$.
\end{lemma}

\noindent\textbf{Proof}.
For any $n\in\NN$ consider on $[a,\,b]$ the grid
$t_k=a+\frac{k}{n}(b-a)\;(k=0,1,\dots,n)$ and the set
$\mathcal{U}_n$ of all functions $u\in C([a,\,b],\,E)$ having the properties:

 -- for any $k\in\{0,1,\dots,n\}\;$ $u(t_k)\in U_{t_k}$;

 -- $u(t)$ is linear in each interval $[t_k,\,t_{k+1}]$, i.e.,
$u(t)=n(t_{k+1}-t)u(t_k)+n(t-t_k)u(t_{k+1})$.

\noindent
It is easy to see that each set $\mathcal{U}_n$ is homeomorphic to
the product $U_{t_0}\times U_{t_1}\times\dots\times U_{t_n}$;
hence, and in view of (a), it is precompact in $C([a,\,b],\,E)$. On
the other hand, (b) easily implies that for any
$\varepsilon>0$ it is possible to choose $n\in\NN$ such that
$\Vert u-u_n\Vert_{C([a,\,b],\,E)}<\varepsilon$ for any $u\in\mathcal{U}\;$,
where $u_n$ is a function from $\mathcal{U}_n$ such that
$u_n(t_k)=u(t_k)\;(k=0,1,\dots,n)$. So, for any $\varepsilon>0$, the set
$\mathcal{U}$ has a precompact $\varepsilon$-net in $C([a,\,b],\,E)\;$;
hence, it is precompact in $C([a,\,b],\,E)$.
\qed

\begin{lemma}\label{lmconvperfunc1}
Let $\{u_n(t)\}_{n=1}^\infty$ be a sequence of continuous functions,
defined in a closed interval $[a,\,b]$, with values in a Banach space
$E$ and having the following properties:

\noindent\
$(a)$ there exists a sequence $\tau_n>0$ such that $\lim_{\,n\rightarrow\infty}\tau_n=0$
and $u_n(t+\tau_n)=u_n(t)$ if $\,t,\,t+\tau_n\in[a,\,b]$;

\noindent\
$(b)$ the sequence $\{u_n(t)\}_{n=1}^\infty$ converges uniformly in $[a,\,b]$ to a function $u:[a,\,b]\to Y$.

\noindent
Then $u(t)$ is constant.
\end{lemma}

\noindent\textbf{Proof}.
By (b), the sequence $\{u_n(t)\}_{n=1}^\infty$ is
equicontinuous, i.e., for any $\varepsilon>0$ there exists $\delta>0$
such that
$|t-s|<\delta\,(t,\,s\in[a,\,b])$ implies
$\Vert u_n(t)-u_n(s)\Vert_E<\varepsilon/3$ for any $n\in\NN$.
By~conditions, we can choose $N>0$ such that
$\tau_n\le\delta$ and $|u(t)-u_n(t)|<\varepsilon/3$ for any $n\ge N$ and $t\in[a,\,b]$.
Let us take $a\le t_1<t_2\le b$.
In view of (a), $u_n(t_1)=u_n(t_1+k_n\tau_n)$, where
$k_n=\Big[\frac{t_2-t_1}{\tau_n}\Big]$. Observe that $0\le
t_2-t_1-k_n\tau_n<\tau_n$. The above arguments imply
\begin{equation*}
\Vert u(t_2)-u(t_1)\Vert_E\le \Vert u(t_2)-u_n(t_2)\Vert_E+\Vert
u(t_1)-u_n(t_1)\Vert_E+\Vert
u_n(t_2)-u_n(t_1+k_n\tau_n)\Vert_E<\varepsilon
\end{equation*}
for $n\ge N$.
Since $\varepsilon>0$ is arbitrary, this proves the claim.
\qed

\smallskip

Now we turn to the main result of this section.
A mapping $f:X\to Y$ between metric spaces is called \textit{compact}
if it is continuous and maps each bounded set in $X$ onto a precompact set in~$Y$.

\begin{theorem}\label{thcomfixpoint}
Let $G$ be a closed bounded convex subset of a Banach space $E$.
Suppose that a one-parametric semigroup of mappings $\{S_t:G\to G\,\}_{t\ge 0}$
has the properties:

$(a)$ for any $t>0$ the mapping $S_t$ is compact;

$(b)$ for any $\;0<a<b$ the family of mappings $S_t$ is continuous by $t$  in the segment $[a,\,b]$

\quad
 uniformly w.r.t. $u\in G$.

\noindent
Then the semigroup $\{S_t\}_{t\ge 0}$ has in $G$ a common fixed point.
\end{theorem}

\noindent\textbf{Proof}.
Shauder's Fixed Point Theorem  claims that a compact mapping of a closed bounded convex set $G$
in a Banach space $E$ into itself has a fixed point, see \cite[p.~74]{aub}.
 By $(a)$ and Shauder's Fixed Point Theorem, for any $\tau>0$ the mapping
$S_\tau$ has a fixed point $u_\tau^0\in G$, i.e., $S_\tau u_\tau^0=u_\tau^0$.
In view of the semigroup property $S_t\circ S_\tau=S_{t+\tau}$,
the function $u_\tau(t)=S_tu_\tau^0$  is $\tau$-periodic.
Take a sequence of positive numbers $\{\tau_n\}_{n=1}^\infty$ such that
$\lim\limits_{\,n\rightarrow\infty}\tau_n=0$, and denote $u_n(t)=u_{\tau_n}(t)$, $u_n^0=u_{\tau_n}^0$.
By conditions, the sequence of functions $\{u_n(t)\}_{n=1}^\infty$
satisfies on $[a,\,b]$ assumptions of Lemma~\ref{genArz};
hence, it is precompact in $C([a,\,b],\,E)$.
Thus, it is possible to select from the sequence $u_n(t)$ a subsequence $u_{n_k}(t)$
converging to a function $u(t)$ in the $C([a,\,b],\,E)$-norm.
Applying Lemma~\ref{lmconvperfunc1} to this subsequence, we find
that $u(t)$ is constant, i.e., $u(t)\equiv u_*$ in $[a,\,b]$.
Since $u_{n_k}(t)\in G$ for any fixed $t\in[a,\,b]$ and $G$ is closed in $E$,
we obtain $u_*\in G$. Since $u_{n_k}(t)=S_{t-a}u_{n_k}(a)\;(t\in[a,\,b])$
and the mappings $S_t:\,G\rightarrow G$ are continuous, we get, tending
$k\rightarrow\infty$, that $S_{t-a}u_*=u_*$ for any $t\in[a,\,b]$.
\qed

\subsection{Solutions of the nonlinear heat equation}

In this section we investigate the existence of solutions of a semi-linear elliptic equation.
For this we prove existence of global solutions of associated non-linear parabolic equation and
study its stable stationary solutions.
We reduce this problem to the existence of a common fixed point for the one-parameter semigroup
of mappings, see Theorem~\ref{thcomfixpoint}, corresponding to the non-linear parabolic equation.
Some results of this section may be known, but for convenience of a reader, we give the proofs.

\subsubsection{Global solutions}

Let $(F,g)$ be a closed Riemannian manifold.
Define a bounded, closed and convex set in $C(F)$ by
\begin{equation}\label{dfG}
 G=\{u\in C(F):\ u_-(x)\le u(x)\le u_+(x) \ \ {\rm for\ all}\ \ x\in F\},
\end{equation}
where $u_-,u_+\in C(F)$, $u_-\le u_+$,
and the following compact domain in $\RR\times F$ by:
\begin{equation}\label{dfD}
 D:=\{(u,\,x)\in\RR\times F:\;u_-(x)\le u\le u_+(x)\}.
\end{equation}
Consider the Cauchy's problem for a non-linear heat equation, more general than \eqref{Cauchy},
\begin{equation}\label{nonlinheat1}
 \partial_t u=\Delta u+{f}(u,x),\quad u(x,0)=u_0(x)\in C(F)\,,
\end{equation}
where ${f}\in C(D)$,
and the stationary version of equation \eqref{nonlinheat1}$_1$:
\begin{equation}\label{nonlinstat}
 \Delta u+{f}(u,x)=0.
\end{equation}

\begin{definition}\label{dfclassolnonlinheat}\rm
A function $u(x,t)$ is a \textit{solution} of \eqref{nonlinheat1} in the domain $F\times[0,T]$,
if  it is continuous, satisfies the initial condition \eqref{nonlinheat1}$_2$,
and in $F\times(0,T]$ it is continuously differentiable by $t$,
twice continuously differentiable by $x$ and satisfies \eqref{nonlinheat1}$_1$.
A~function $u(x)$ is a \textit{solution} of \eqref{nonlinstat} in $F$, if it belongs to $C^{2}(F)$ and satisfies this equation in $F$.
\end{definition}

Let $S_t:C(F)\rightarrow C(F)$ be the map which relates to each initial value $u_0\in C(F)$ the value
of the classical solution of \eqref{nonlinheat1} at the moment $t\in[0,T)$ (if this solution exists and is unique).
Since ${f}(u,x)$ does not depend explicitly on $t$, the family $\{S_t\}_{0\le t<T}$ has
the semigroup property, and it is a semigroup (i.e., $T=\infty$) when \eqref{nonlinheat1}
admits a global solution for any $u_0(x)\in C(F)$.

It is known, see \cite[Theorem~B.6.3]{cc2}, that the Cauchy's problem for the heat~equation,
\begin{equation}\label{linheat1}
\partial_tv=\Delta v,\quad v(x,0)=v_0(x)
\end{equation}
admits a unique global solution $v(x,t)$ for any $v_0\in C(F)$. Let
$S_t^0:\,C(F)\rightarrow C(F)\ (t\ge0)$ be the semigroup of linear
mappings corresponding to \eqref{linheat1}. Since $v(x,t)$ is
continuous in $F\times[0,T]$ then $S_t^0v_0\to v_0$ in the
$C(F)$-norm as $t\downarrow 0$. This means that the semigroup
$\{S_t^0\}$ is strongly continuous in $C(F)$. It is known that
\begin{equation}\label{reprsolCauch}
v(x,t)=(S_t^0v_0)(x)=\int_FH(x,y,t)\,v_0(y)\,dy,
\end{equation}
where $H(x,y,t)$ is the fundamental solution of
\eqref{linheat1}$_1$, called the heat kernel, which belongs to
$C^\infty(\Omega)$ in the domain $\Omega:=F\times F\times\{t\in
\RR:\,t>0\}$, see~\cite{cc2}. We shall use the properties
\begin{equation}\label{propheatkern}
 H(x,\xi,t)>0,\quad\int_F H(x,\xi,t)\,{\rm d}\xi=1\quad (x\in F,\ t>0).
\end{equation}
If a solution of \eqref{nonlinheat1} exists then it satisfies the  integral equation (Duhamel's principle):
\begin{equation}\label{nonlininteq}
 u(x,t)=\int_F H(x,y,t)\,u_0(y)\,{\rm d}y +\int_0^t\int_F H(x,y,t-\tau)\,{f}(u(y,\tau),y)\,{\rm d}y\dtau\,.
\end{equation}
Denote by $\mathbb{T}$ the set of all $T>0$ such
that~\eqref{nonlininteq} has a continuous solution in the domain $F\times[0,T]$.
 For $u_0\in\mathrm{Int}(G)$ and $r>0$, let
$B_r(u_0)=\{u\in C(F):\,\Vert u-u_0\Vert_{C(F)}\le r\}$
be a closed ball contained in~$G$.
One may take $r=\min\{\min_{\,F}|u_0-u_-|,\ \min_{\,F}|u_0-u_+|\}$.

\begin{proposition}\label{P-local-sol}
 If $q\in C^\infty(D)$ and $u_0\in\mathrm{Int}(G)$ then

\noindent\
$(i)$
$\mathbb{T}\neq\emptyset$, namely, there are $T>0$ and $r>0$ such that
the integral equation~\eqref{nonlininteq} has a unique continuous solution $u(x,t)$ in
$F\times[0,T]$, with the property $u(\cdot\,,t)\in B_r(u_0)$ for any~$t\in(0,T]$;

\noindent\
$(ii)$ for any $T\in\mathbb{T}$, a continuous solution $u(x,t)$ of \eqref{nonlininteq}
in the domain $F\times[0,T]$ is a solution of \eqref{nonlinheat1}.
Moreover, $u(\cdot,\cdot)\in C^\infty(F\times(0,T])$;

\noindent\
 $(iii)$ for any $T\in\mathbb{T}$, a solution of \eqref{nonlinheat1} is unique in $F\times[0,T]$.
\end{proposition}

\noindent\textbf{Proof}.
satisfies the Lipschitz condition w.r.t. $u$, i.e., there exists $L>0$ such that
\begin{equation}\label{lipcond}
 |{f}(u_2,x)-{f}(u_1,x)|\le L|u_2-u_1|\quad \forall\,(u_1,x)\,,(u_2,x)\in D.
\end{equation}
Hence, the superposition operator $(\mathcal{X}u)(x)={f}(u(x),x)$ satisfies
the Lipschitz condition:
 $\Vert\mathcal{X}u_1-\mathcal{X}u_2\Vert_{C(F)}\le L\Vert u_1-u_2\Vert_{C(F)}\ (\forall\,u_1,u_2\in G)$.
Let $\Theta$ be the operator expressed by the rhs of~\eqref{nonlininteq} and defined on
the set $C([0,T],B_r(u_0))$, which is closed in the Banach space $C([0,T],C(F))$.
By Lemma~\ref{lmproplinsemgr}$(i)$ and
the proof of Proposition~1.1 in \cite[p.~315]{ta},
we can choose $T>0$ such that
$\Theta$ maps the set
$C([0,T],B_r(u_0))$ into itself and it is a contraction there.
Hence, \eqref{nonlininteq} has in $C([0,T],B_r(u_0))$ a unique solution.

$(ii)$ The proof consists of two steps.

\noindent \textbf{Step 1}.
Let us show that $u(\cdot,t)\in C^1(F)$ for any $t\in(0,T]$.
Since $H(\,\cdot\,,\cdot\,,\cdot\,)\in C^\infty(\Omega)$,
the first integral in (\ref{nonlininteq}) belongs to class $C^\infty(F\times(0,T])$.
It remains to prove that the second integral in (\ref{nonlininteq}),
denoted by $I_0(x,t)$, belongs to $C^1(F)$ for any $t\in(0,T]$.
Consider the truncated integral
 $I_\varepsilon(x,t)=\int_0^{t-\varepsilon}\int_F H(x,y,t-\tau)\,{f}(u(y,\tau),y)\,{\rm d}y\dtau$ for $\varepsilon\in(0,t)$.
We have
\begin{eqnarray*}
 |\,I_\varepsilon(x,t)-I_0(x,t)\,|
 \eq\big|\int_{t-\varepsilon}^t\int_FH(x,y,t-\tau)\,{f}(u(y,\tau),y)\,{\rm d}y\dtau\big|\\
 &\le&\Vert{f}(\cdot,\cdot)\Vert_{C(D)}\int_{t-\varepsilon}^t\int_F H(x,y,t-\tau)\,{\rm d}y\dtau
 =\varepsilon\,\Vert{f}(\cdot,\cdot)\Vert_{C(D)}.
\end{eqnarray*}
Hence, for any $t\in(0,T]$, the integral $I_\varepsilon(x,t)$ converges to $I_0(x,t)$ as $\varepsilon\downarrow 0$ uniformly on $F$.

Observe that since $H(\,\cdot\,,\cdot\,,\cdot\,)\in C^\infty(\Omega)$, thus
$I_\varepsilon(\cdot\,,t)\in C^\infty(F)$. Hence, in order to prove that
$I(\cdot\,,t)\in C^1(F)$ for $t>0$, it is sufficient to show that
the first order partial derivatives of $I_\varepsilon(x,t)$ by all variables
converge as $\varepsilon\downarrow 0$ uniformly
for any local coordinates $(x_k)$ with compact support $W$
on~$F$.
 Take $x\in W$ and consider derivatives
\[
 \partial_{x_k} I_\varepsilon(x,t)=\int_0^{t-\varepsilon}\partial_{x_k}\int_F  H(x,y,t-\tau)\,{f}(u(y,\tau),y)\,{\rm d}y\dtau.
\]
Using \eqref{reprsolCauch} and estimate
 $\Vert S_t^0\Vert_{\mathcal{B}(C(F),C^1(F))}\le C\,t^{-1/2}$ with
 $t\in(0,1]$,
see  \cite[(1.11), p.~315]{ta},
we have for $0<\varepsilon_1<\varepsilon_2<t$:
\begin{eqnarray*}
 &&|\partial_{x_k} I_{\varepsilon_1}(x,t) -\partial_{x_k} I_{\varepsilon_2}(x,t)|
 \le \int_{t-\varepsilon_2}^{t-\varepsilon_1} |\partial_{x_k}( S^0_{t-\tau}\mathcal{X}(u))(x)|\,\dtau\\
 &&\le C\,\Vert\mathcal{X}(u)\Vert_{C(F)}\int_{t-\varepsilon_2}^{t-\varepsilon_1}(t-\tau)^{-1/2}\,\dtau
 \le 2\,C\,\sqrt{\varepsilon_2}\,\|{f}(\cdot,\cdot)\|_{C(D)}.
\end{eqnarray*}
This estimate shows us that the following integral exists:
\[
 J_k(x,t):=\int_0^t\int_F \partial_{x_k} H(x,y,t-\tau)\,{f}(u(y,\tau),y)\,{\rm d}y\dtau,
\]
and $\partial_{x_k} I_{\varepsilon}(x,t)\rightarrow J_k(x,t)$
as $\varepsilon\downarrow 0$ uniformly on $F$ for any $t\in(0,T]$.
Hence, $I(\cdot,t)\in C^1(F)$, and, therefore, $u(\cdot,t)\in C^1(F)$ for any $t\in(0,T]$.

\noindent \textbf{Step 2}. Let us show that $u(\cdot,\cdot)\in C^\infty(F\times(0,T])$.
Observe that for any $\sigma\in (0,T)$, the restriction of $u(x,t)$ on $F\times[\sigma,T]$
is a solution of the integral equation
\[
 u(x,t)=\int_F H(x,y,t)\,u_\sigma(y)\,{\rm d}y
 +\int_\sigma^t\int_F H(x,y,t-\tau)\,{f}(u(y,\tau),y)\,{\rm d}y\dtau,
\]
where $u_\sigma(x)=u(x,\sigma)$. By Step 2, $u_\sigma\in C^1(F)$.
Taking into account that ${f}(\cdot,\cdot)\in C^\infty(D)$, and using \cite[Proposition 1.2, p.~316]{ta},
we obtain that $u\in C^\infty(F\times(\sigma, T])$.
Since $\sigma\in (0,T)$ is arbitrary, then $u\in C^\infty(F\times(0, T])$.
Furthermore, one may conclude from (\ref{nonlininteq}), that
$u(x,t)$ satisfies initial condition~\eqref{nonlinheat1}$_2$, and in
the domain $F\times(0, T]$ it obeys \eqref{nonlinheat1}$_1$.

$(iii)$
Assume that~\eqref{nonlinheat1} has two solutions $u_1(x,t)$ and $u_2(x,t)$ in the domain $F\times[0,T]$.
Then, in view of~\eqref{lipcond}, the function $w(x,t)=u_2(x,t)-u_1(x,t)$ satisfies the differential inequalities:
\begin{equation*}
 \Delta w -L|w|\le\partial_t w\le\Delta w +L|w|.
\end{equation*}
By the maximum principle, $w_-(t)\le w(x,t)\le w_-(t)$, where $w_-(t),w_+(t)$ solve the
problems
\[
 {d w_-}/{dt}=-L\,|w_-(t)|,\ \ w_-(0)=0,\qquad
 {d w_+}/{dt}=-L\,|w_+(t)|,\ \ w_+(0)=0.
\]
Hence, $w(x,t)\equiv 0$ in $F\times[0,T]$. \qed

\begin{theorem}\label{prglob}
Suppose that $q\in C^\infty(D)$ and $u_0\in\mathrm{Int}(G)$. If
there exist continuous functions $\tilde u_-(x)$ and $\tilde u_+(x)$
such that $u_-<\tilde u_-<\tilde u_+<u_+$, and for any
$T\in\mathbb{T}$ for the solution $u_T(x,t)$ of Cauchy's problem~\eqref{nonlinheat1}
the estimates $\tilde u_-(x)\le u_T(x,t)\le \tilde u_+(x)$ are valid in the domain $F\times[0,T]$.
Then~\eqref{nonlinheat1} has a global solution $u(x,t)$, i.e., it is
defined in the domain $F\times(0,\infty)$. Furthermore, it is unique
there and satisfies the inequalities $\tilde u_-(x)\le u(x,t)\le \tilde u_+(x)$.
Moreover, $u(\cdot,\cdot)\in F\times(0,\infty)$.
\end{theorem}

\noindent\textbf{Proof}. By Proposition~\ref{P-local-sol}$(i)$,
$\mathbb{T}\neq\emptyset$. Denote $\tilde T=\sup(\mathbb{T})$.
We should prove that $\tilde T=\infty$. Assume on the contrary that $\tilde T<\infty$.
Since by Proposition~\ref{P-local-sol}(iii), for any $T\in\mathbb{T}$,
$u_T(x,t)$ is a unique solution of \eqref{nonlinheat1} in $F\times[0,T]$,
then we can consider the function $u(x,t)$, defined on
$F\times[0,\tilde T)=\bigcup_{\,T\in\mathbb{T}}F\times[0,T]$ such that for any
$T\in\mathbb{T}$ $u_T=u\vert_{F\times[0,T]}$.
It is a unique solution of \eqref{nonlinheat1} in the domain $F\times[0,\tilde T)$;
hence, it satisfies in this domain the integral equation~\eqref{nonlininteq}.
We have for $(x,t_k)\in F\times[0,\tilde T)$ using~\eqref{propheatkern}:
\begin{eqnarray*}
&&\big|\int_0^{t_2}\int_F H(x,y,t-\tau)\,{f}(u(y,\tau),y)\,{\rm d}y\dtau
 -\int_0^{t_1}\int_F H(x,y,t-\tau)\,{f}(u(y,\tau),y)\,{\rm d}y\dtau\big|\\
&&\le |t_2-t_1|\cdot\Vert{f}(\cdot,\cdot)\Vert_{C(D)}.
\end{eqnarray*}
This estimate and~\eqref{nonlininteq} show us that $u(x,t)$ tends to a
continuous function $\tilde u(x)$ as $t\uparrow\tilde T$ in the
$C(F)$-norm. Since $u_-(x)<\tilde u_-(x)\le u(x,t)\le\tilde
u_+(x)<u_+(x)$ in $F\times[0,\tilde T)$, then $\tilde u\in\mathrm{Int}(G)$.
Therefore, by Proposition~\ref{P-local-sol}(i)-(ii) there exists $\delta>0$
such that the Cauchy's problem
\begin{equation*}
 \partial_t v=\Delta u+{f}(v,x),\quad v(x,\tilde T)=\tilde u(x)
\end{equation*}
has a solution $v(x,t)$ in $F\times[\tilde T,\tilde T+\delta]$.
It is easy to check that the function
\begin{displaymath}
 w(x,t)=\left\{\begin{array}{ll} u(x,t),&\mathrm{if}\quad (x,t)\in F\times[0,\tilde T),\\
 v(x,t),&\mathrm{if}\quad (x,t)\in F\times[\tilde T,\tilde T+\delta]
\end{array}\right.
\end{displaymath}
is a continuous solution of the integral equation~\eqref{nonlininteq} in $F\times[0,\tilde T+\delta]$.
This fact contradicts to the definition of the number $\tilde T$.
Thus, $\tilde T=\infty$; hence, $u(x,t)$ is a unique global solution
of Cauchy's problem~\eqref{nonlinheat1} satisfying in $F\times[0,\infty)$
the estimates $\tilde u_-(x)\le u(x,t)\le \tilde u_+(x)$.
Furthermore, by Proposition~\ref{P-local-sol}$(ii)$,
$u(\cdot,\cdot)\in C^\infty(F\times(0,\infty))$.
\qed

\subsubsection{Stationary solutions}

The proof of the following theorem is supported by Lemmas~\ref{lmcontSt}--\ref{lmintcompmaps} in what follows.

\begin{theorem}\label{P-exist}
Let the following conditions are satisfied:

-- ${f}\in C^\infty(D)$, for $D$ in \eqref{dfD},

-- \eqref{nonlinheat1} admits a global solution for any $u_0(x)\in G$, and

-- the set $G$ is invariant w.r.t. the corresponding semigroup $S_t\ (t\ge0)$.

\noindent
Then the set of all solutions of \eqref{nonlinstat} lying in $G$ is nonempty and compact in $C(F)$.
\end{theorem}

\noindent\textbf{Proof}. Take $u_0\in G$. By the Duhamel's principle, we have
\begin{equation}\label{Duhamel}
 S_tu_0=S_t^0u_0+\int_0^tS_{t-\tau}^0 {f}(S_\tau u_0,\,\cdot)\,\dtau\,,
\end{equation}
where $S_t^0$ is the semigroup associated with \eqref{linheat1}$_1$.
Denote by $\|\cdot\|_{\cal{B}(C(F))}$ the operator norm.
For any $t\in[a,\,b]\;(0<a<b)$, $\delta\in(0,\,a)$ and $h\in(0,\,\delta)$, we have
\begin{eqnarray}\label{E-sect-3-2}
\nonumber
&&\hskip-6mm\Vert u(\cdot,t+h)-u(\cdot,t)\Vert_{C(F)}
 \le\delta\max_{\tau\in[t+h-\delta,\,t+h]}
 \Vert S_{t+h-\tau}^0\Vert_{\cal{B}(C(F))}\cdot\Vert {f}(S_\tau u_0,\,\cdot)\Vert_{C(F)}\\
&&\hskip-6mm +\,\delta\max_{\tau\in[t-\delta,\,t]}\big(\Vert S_{t-\tau}^0\Vert_{\cal{B}(C(F))}
 \cdot\Vert {f}(S_\tau u_0,\,\cdot)\Vert_{C(F)}\big)
 +\Vert S_{t+h}^0-S_t^0\Vert_{\cal{B}(C(F))}\cdot\Vert u_0\Vert_{C(F)}\\
\nonumber
&&\hskip-6mm +h\!\max_{\tau\in[t-\delta,t-\delta+h]}
 \Vert S_{t+h-\tau}^0\Vert_{\cal{B}(C(F))}\Vert {f}(S_\tau u_0,\cdot)\Vert_{C(F)}
 {+}\max_{\tau\in[0,t-\delta]}
 \Vert S_{t+h-\tau}^0{-}S_{t-\tau}^0\Vert_{\cal{B}(C(F))}\Vert {f}(S_\tau u_0,\cdot)\Vert_{C(F)}.
\end{eqnarray}
Given any $\varepsilon>0$, by Lemma~\ref{lmproplinsemgr}(i), we can choose
$\delta>0$ such that the sum of the first two terms in the rhs of
\eqref{E-sect-3-2} is less than $\varepsilon/2$ for any $u_0\in G$.
Furthermore, in view of Lemma~\ref{lmproplinsemgr}(iii), the family
$\{S_t^0\}$ is uniformly continuous by $t$ in the operator norm on
each compact interval which does not contain $t=0$, we can choose
$h>0$ such that the sum of the remain terms in the rhs of the last
estimate will be less than $\varepsilon/2$ for any $u_0\in G$. This means
that the semigroup $S_t$ is continuous by $t$ in $[a,\,b]$ uniformly
w.r.t. $u_0\in G$ for any $0<a<b$. Then, in view of the
continuity of ${f}(u,x)$ in $D$ and the invariance of $G$ with respect
to the semigroup $S_t$, the family of mappings $Q_tu_0:={f}(S_t u_0,\,\cdot)$
is continuous by $t$ in $[a,\,b]$ uniformly w.r.t. $u_0\in G$ for any $0<a<b$.
These circumstances, equality (\ref{Duhamel}),
Lemmas~\ref{lmcontSt}, \ref{lmproplinsemgr}(i-ii) and \ref{lmintcompmaps}(ii)
imply that each mapping $S_t$ with $t>0$ is compact on $G$. So,
$S_t$ satisfies all conditions of Theorem~\ref{thcomfixpoint}.
Hence, it has in $G$ a common fixed point $u_*(x)$, i.e.,
$S_tu_*=u_*$ for any $t>0$. On the other hand, it is known
that for any $u_0\in G$ and $t>0$ $S_tu_0\in C^\infty(F)$ (see
Proposition~\ref{P-local-sol}). Hence $u_*(x)$ belongs to
$C^\infty(F)$ and it is a solution of \eqref{nonlinstat}.

By continuity of $S_t u_0$ by $u_0$, the set ${\rm Fix}(G)$ of all common
fixed points of $S_t\,(t>0)$ in $G$ is closed w.r.t. the
$C$-norm. Since $S_t({\rm Fix}(G))={\rm Fix}(G)$ for $t>0$, and
$S_t$ maps any $C$-bounded set on a $C$-precompact set, then ${\rm
Fix}(G)$ is $C$-precompact. Thus, ${\rm Fix}(G)$ is $C$-compact.
\qed

\begin{lemma}\label{lmcontSt}
In conditions of Theorem~\ref{P-exist}, for any $t>0$, the mapping $S_t:\,G\rightarrow G$ is continuous.
\end{lemma}

\noindent\textbf{Proof}.
Take $u_1^0,u_2^0\in G$ and denote
$u_k(x,t)=(S_tu_k^0)(x)\;(k=1,2)$. Then, in view of \eqref{lipcond},
the function $w(x,t)=u_2(x,t)-u_1(x,t)$ satisfies the differential
inequalities:
\begin{equation*}
 \Delta w -L|w|\le\partial_t w\le\Delta w +L|w|.
\end{equation*}
Let $w_-(t),w_+(t)$ be solutions of the following
Cauchy's problems with $w_0=\Vert u_2^0-u_1^0\Vert_{C(F)}$:
\[
 {d w_-}/{dt}=-L\,|w_-(t)|,\ \ w_-(0)=-w_0,\qquad
 {d w_+}/{dt}=-L\,|w_+(t)|,\ \ w_+(0)=w_0.
\]
By the maximum principle, $w_-(t)\le w(x,t)\le w_-(t)$ and $|w(x,t)|\le w_0\,e^{-Lt}$.
\qed

\begin{lemma}\label{lmproplinsemgr}
The semigroup $S_t^0:\,C(F)\rightarrow C(F)$ has the properties:

$(i)$ $\Vert S_t^0\Vert_{\cal{B}(C(F))}\le 1$ for any $t\ge 0$;

$(ii)$ the linear operator $S_t^0$ is compact for any $t>0$;

$(iii)$ the family $S_t^0$ is continuous by $t\in(0,\,\infty)$ in the operator norm.
\end{lemma}

\noindent\textbf{Proof}. For $v_0\in C(F)$ denote $v(x,t)=(S_t^0v_0)(x)$ .
By \eqref{reprsolCauch} and \eqref{propheatkern}, we get for $x,\,y\in F$,~$t>0$:
\begin{equation}\label{estnrmSt0}
 |v(x,t)|\le\Vert v_0\Vert_{C(F)},
\end{equation}
\begin{equation}\label{estnrmdifSt0}
 |v(x,t)-v(y,t)|\le{\rm Vol}(F)\sup\nolimits_{\,\xi\in F}|H(x,\xi,t)-H(y,\xi,t)|\cdot\Vert v_0\Vert_{C(F)}.
\end{equation}
Thus, \eqref{estnrmSt0} implies $(i)$.

Consider the unit ball $B_1=\{f\in C(F): \Vert f\Vert_{C(F)}\le 1\}$ in $C(F)$.
Estimates (\ref{estnrmSt0}), (\ref{estnrmdifSt0})
and continuity of the heat kernel $H(x,\xi,t-\tau)$ on each compact of the form
$K_\delta=F\times F\times\{(t,\,\tau):\;0 \le\tau\le t-\delta\}\;(\delta>0)$
imply that for $t>0$ the set $S_t^0(B_1)$ is bounded in $C(F)$ and it is equicontinuous.
By the Arzela-Ascoli Theorem, it is precompact in $C(F)$. This proves~$(ii)$.

Let us prove $(iii)$. As above, put $v(x,t)=(S_t^0v_0)(x)$.
For $t_1,\,t_2\in(0,\,\infty)$ and $x\in F$ we~get
\begin{equation*}
 |v(x,t_1)-v(x,t_2)|\le{\rm Vol}(F)\sup\nolimits_{\,\xi\in F}
 |H(x,\xi,t_1)-H(x,\xi,t_2)|\cdot\Vert v_0\Vert_{C(F)}.
\end{equation*}
This estimate and the continuity of the heat kernel on each compact
$K_\delta$,
imply~$(iii)$.
\qed

\begin{lemma}\label{lmconvcompmap}
Let $\{T_n\}_{n=1}^\infty$ be a set of compact mappings acting
from a bounded subset $B$ of a Banach space $E_1$ into a Banach
space $E_2$ and converging uniformly to
$T:B\rightarrow E_2$.
Then $T$ is compact.
\end{lemma}

\noindent\textbf{Proof}.
The continuity of $T$ is obvious. Take an arbitrary $\varepsilon>0$ and choose $n\in\NN$
such that $\sup_{\,x\in B}\Vert Tx-T_nx\Vert_{E_2}<\varepsilon$.
This means that the set $T_n(B)$ forms a precompact $\varepsilon$-net for the set $T(B)$ in $E_2$.
Hence the set $T(B)$ is precompact in $E_2$.
\qed

\begin{lemma}\label{lmintcompmaps}
Let $\{T_t\}_{t\in[a,\,b)}$ be a family of compact mappings acting
from a bounded subset $B$ of a Banach space $E_1$ into a Banach
space $E_2$.


$(i)$ If $c\in (a,\,b)$ and $T_t$ is continuous by $t$ on $[a,\,c]$
uniformly w.r.t. $x\in B$ then the mapping
$J^cx:=\int_a^cT_tx\,{\rm d}t\ (x\in B)$ is compact;


$(ii)$ If the condition of $(i)$ is satisfied for any $c\in (a,\,b)$
and the family $J^c$ converges as $c\uparrow b$ to the mapping $Jx=\int_a^bT_tx\,{\rm d}t$
uniformly w.r.t. $x\in B$ then $J$ is compact.
\end{lemma}

\noindent\textbf{Proof}.
($i$) For any $n\in\NN$ consider on $[a,\,c]$ the grid $t_k=a+\frac{k}{n}(c-a)\;(k=1,2,\dots,n)$
and the mapping $J^c_nx=\frac{1}{n}(c-a)\sum_{\,k=1}^nT_{t_k}x\;(x\in B)$.
One may show that each $J^c_n$ is compact and the sequence
$\{J^c_n\}_{n=1}^\infty$ converges to the mapping $T$ uniformly.
By Lemma~\ref{lmconvcompmap}, $T$ is compact.
Thus, $(ii)$ follows from $(i)$ and Lemma~\ref{lmconvcompmap}.
\qed

\subsection{Attractors of the nonlinear heat equation}
\label{sec:attr-all}

This section studies stable stationary solutions of \eqref{Cauchy} for three cases.

\subsubsection{Case of $\Psi_3>0$}
\label{sec:attr2}

Let $\Psi_3>0$, $\Psi_1>0$, $\Psi_2>0$ and $\lambda_0>0$, see Section~\ref{R-burgers-heat}, case (a).
For $y>0$, put
\[
 \phi(y,\theta) = -\lambda_0\,y +\theta_1\,y^{-1}-\theta_2\,y^{-3} +\theta_3\,y^{3} = P_{\,\phi}(y^2)/y^3,
\]
where $P_{\,\phi}(z) = \theta_3 z^{3} -\lambda_0 z^{2} +\theta_1 z-\theta_2$ and
$\theta=(\theta_1,\theta_2,\theta_3)\in\mathcal{P}$.
Then $\phi_-(y)\le\phi(y,\theta)\le\phi_+(y)$~for
\begin{equation*}
 \phi_+(y) = -\lambda_0y +\Psi_1^+ y^{-1} -\Psi_2^- y^{-3} +\Psi_3^+ y^{3},\quad
 \phi_-(y) = -\lambda_0y +\Psi_1^- y^{-1} -\Psi_2^+ y^{-3} +\Psi_3^- y^{3}.
\end{equation*}
The discriminant of $P_{\,\phi}(z)$ is the following cubic polynomial in $\lambda_0$:
\begin{equation}\label{E-discr2}
 D(P_\phi)(\lambda_0) = -4\,\theta_2\lambda_0^3 +\theta_1^2\lambda_0^2 +18\,\theta_1\theta_2\theta_3\lambda_0
 -\theta_3(4\,\theta_1^3 +27\,\theta_2^2\theta_3)\,.
\end{equation}
If $D(P_\phi)>0$ for some $\lambda_0>0$
then $P_\phi(z)$ has 3 real roots $z_3(\theta)<z_2(\theta)<z_1(\theta)$, and $y_k=z^2_k\ (k=1,2,3)$ are roots of $\phi(y,\,\cdot)$.
Since $P_\phi(z)<0$ for $z<0$, all its roots are positive.

 By Maclaurin method, positive $\lambda_0$-roots of $D(P_\phi)$ are bounded above by
\[
 1+\max\{\,\theta_1^2,\ 18\,\theta_1\theta_2\theta_3\,\}/(4\,\theta_2)
 \le K := 1+\max\{\,(\Psi_1^+)^2,\ 18\,\Psi_1^+\Psi_2^+\Psi_3^+\,\}/(4\,\Psi_2^-).
\]
Since $D(P_\phi)(-\infty)=\infty$ and $D(P_\phi)(0)<0$ for any $\theta\in\mathcal{P}$,
there is one negative root.
Indeed, by Vieta's formulas, the sum of $\lambda_0$-roots is $\theta_1^2/(4\theta_2)>0$; hence, three negative roots are impossible.

\begin{figure}
 \begin{tabular*}{\textwidth}{@{}l@{}l@{}r@{}}
\hskip25mm\subfigure{\includegraphics[width=0.34\textwidth]{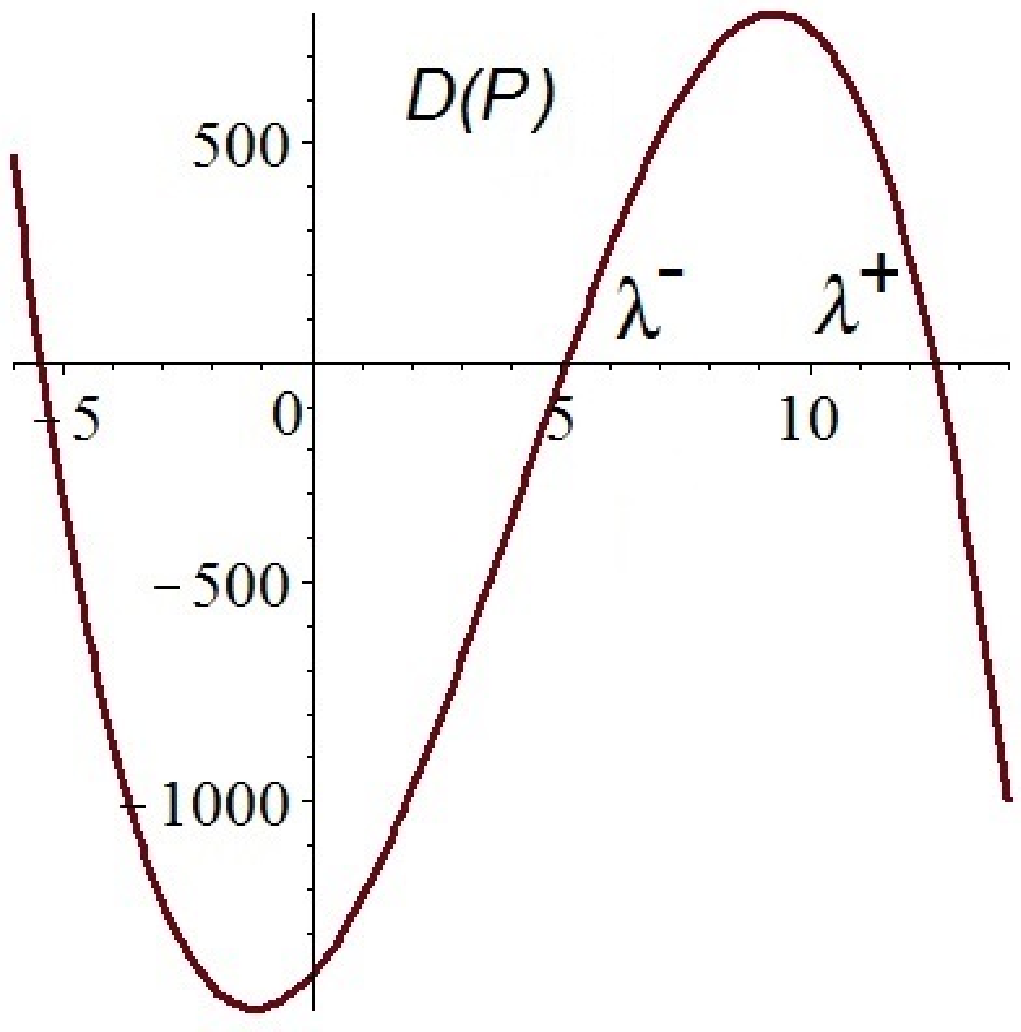}}
\hskip5mm\subfigure{\includegraphics[width=0.28\textwidth]{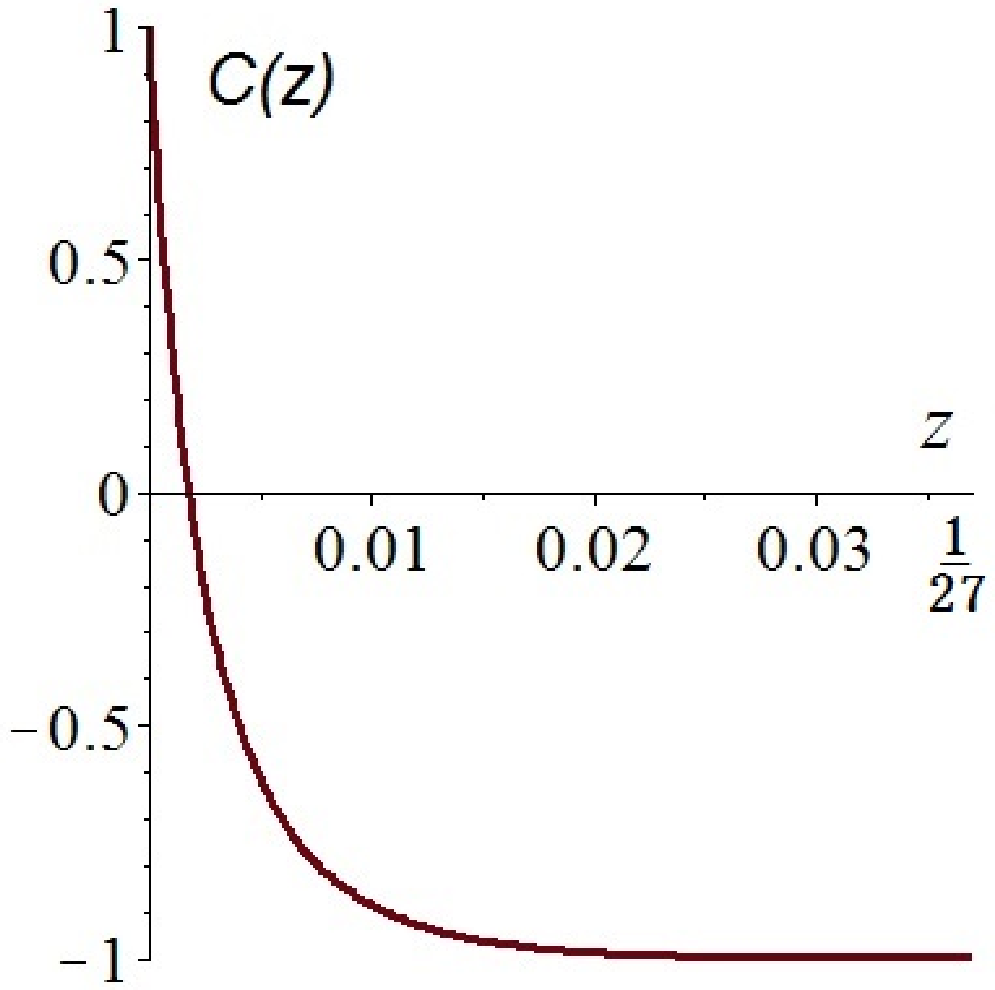}}
\end{tabular*}
\vskip-2mm
\caption{\small
(a)~Roots of $D(P_\phi)$.\quad (b)~Graph of $C(z)$ for $0\le z\le 1/27\approx0.037$.}
\label{fig:3}
\end{figure}

The discriminant by $\lambda_0$ of $D(P_\phi)$ is $16\,\theta_{3}(\theta_{1}^{3}-27\,\theta_{2}^{2}\theta_{3})^{3}$.
If $\theta_{1}^{3}<27\,\theta_{2}^{2}\theta_{3}$ then $D(P_\phi)$ has one real $\lambda_0$-root,
which as was shown is negative; this case in not useful for us, because $D(P_\phi)<0$ for $\lambda_0>0$.
If $\theta_{1}^{3}>27\,\theta_{2}^{2}\theta_{3}$ then $D(P_\phi)$ has three real $\lambda_0$-roots: one negative
and other two positive, $\lambda^+(\theta)>\lambda^-(\theta)$;
moreover, $D(P_\phi)>0$ when $\lambda_0\in I_\lambda(\theta)=(\lambda^-(\theta),\,\lambda^+(\theta))$, Fig.~\ref{fig:3}(a).
In~this case, $\phi(y,\theta)$ has three positive roots $y_1(\theta)>y_2(\theta)>y_3(\theta)$,
$\partial_y\phi(y,\theta)$ has two positive roots $y_4(\theta)\in(y_2(\theta),\,y_1(\theta))$ and $ y_5(\theta)\in(y_3(\theta),\,y_2(\theta))$.
Thus, in what follows we assume
\begin{equation}\label{E-Psi-2}
 (\Psi^-_{1})^{3}>27\,(\Psi^+_{2})^{2}\Psi^+_{3}.
\end{equation}
Since ${\cal P}$ is compact,
there exist $\Lambda^-=\max_{\,\cal P}\lambda^-(\theta)$ and $\Lambda^+=\min_{\,\cal P}\lambda^+(\theta)$.

Denote by
$y_3^+ < y_2^+ < y_1^+$ the positive roots of $\phi_+(y)$,
by
$y_3^- < y_2^- < y_1^-$ the positive roots of $\phi_-(y)$,
and $y_5^- < y_4^-$ and $y_5^+ <  y_4^+$, respectively, the positive roots of functions
\begin{eqnarray*}
 (\partial_y\phi)_-(y) = -\lambda_0 -\Psi_1^+ y^{-2} +3\Psi_2^- y^{-4} +3\Psi_3^- y^2,\\
 (\partial_y\phi)_+(y) = -\lambda_0 -\Psi_1^- y^{-2} +3\Psi_2^+ y^{-4} +3\Psi_3^+ y^2.
\end{eqnarray*}
We calculate
$\partial_y\phi(y,\theta)=-\lambda_0-\theta_1 y^{-2}+3\,\theta_2 y^{-4} +3\,\theta_3 y^2$.
For any $\theta\in\mathcal{P}$ and $y>0$ we have
\begin{equation}\label{estderivphi1}
 (\partial_y\phi)_-(y)\le(\partial_y\phi)(y,\theta)\le(\partial_y\phi)_+(y).
\end{equation}

We need the following condition:
\begin{equation}\label{E-Psi-1}
 3\,\Psi^-_{3} > \Psi^+_{3}.
\end{equation}

\begin{proposition}\label{P-auxprop-2}
If \eqref{E-Psi-2} holds then, for any $\theta\in\mathcal{P}$ and $\lambda_0\in I_\lambda(\theta)$, we have
\begin{eqnarray*}
 && y_3^+ \le y_2(\theta) \le y_3^-,\quad y_2^- \le y_2(\theta) \le y_2^+,\quad y_1^+ \le y_1(\theta) \le y_1^-,\\
 && y_5^- \le y_5(\theta) \le y_5^+,\quad y_4^- \le y_4(\theta) \le y_4^+ .
\end{eqnarray*}
If, in addition, \eqref{E-case1-mu}, \eqref{E-Psi-1}  and
\begin{equation}\label{E-case1-P}
 \delta_3^2 \le \min\!\Big\{1, \frac{8\Psi_3^- D_{P^-}}
 {27(\Psi_2^+)^{2} {+}18(4\Psi_1^-\Psi_2^+\lambda_0 {+}(\Psi_1^-)^{3} {+}9(\Psi_2^+)^{2}\Psi_3^-)},
 \frac{(3\,\Psi_3^- -\Psi_3^+)D_{P^+}}{9(4\Psi_1^+\Psi_2^-\lambda_0 {+}(\Psi_1^+)^{3} {+}9(\Psi_2^-)^{2}\Psi_3^+)}
 \Big\}
\end{equation}
hold for any $\lambda_0\in(\Lambda^- +\varepsilon,\,\Lambda^+ -\varepsilon)$
and some positive $\varepsilon<\frac12\,(\Lambda^+ -\Lambda^-)$ then
\begin{equation}\label{E-auxcor}
 y_3^+ < y_3^- < y_5^+ < y_2^- < y_2^+ < y_4^- < y_1^+ < y_1^-\,.
\end{equation}
\end{proposition}

\noindent\textbf{Proof}.
For implicit derivatives
 $\partial_{\theta_k}y_l = -({\partial_{\theta_k}\phi}/{\partial_y\phi})\,\vert_{\,y=y_l(\theta)}$,
 $\partial_{\theta_k}y_j = -({\partial\,^2_{\theta_k y}\phi}/{\partial\,^2_{yy}\phi})\,\vert_{\,y=y_j(\theta)}$
 where $k,l=1,2,3,\ j=4,5$,
we calculate
\begin{eqnarray*}
 &&\partial_{\theta_1}\phi(y,\theta) = y^{-1},\quad
 \partial_{\theta_2}\phi(y,\theta) = -y^{-3},\quad
 \partial_{\theta_3}\phi(y,\theta) =  y^3,\\
 &&\partial_{y}\phi\vert_{y=y_1(\theta)} > 0,\ \
 \partial_{y}\phi\vert_{y=y_2(\theta)} < 0,\ \
 \partial_{y}\phi\vert_{y=y_3(\theta)} > 0,\ \
 \partial^2_{yy}\phi\vert_{y=y_4(\theta)} > 0,\ \ \partial^2_{yy}\phi\vert_{y=y_5(\theta)} < 0, \\
 &&\partial^2_{\theta_1 y}\phi(y,\theta) = -y^{-2},\quad
 \partial^2_{\theta_2\,y}\phi(y,\theta) = 3\,y^{-4},\quad
 \partial^2_{\theta_3 y}\phi(y,\theta) = 3\,y^2\,,
\end{eqnarray*}
where $\partial^2_{yy}\phi(y,\theta)= 2\,\theta_1 y^{-3}-12\,\theta_2 y^{-5} +6\,\theta_3 y$.
Thus, the following inequalities hold:
\begin{eqnarray*}
 \partial_{\theta_1}y_1(\theta)<0,\quad
 \partial_{\theta_1}y_2(\theta)>0,\quad
 \partial_{\theta_1}y_3(\theta)<0,\quad
 \partial_{\theta_1}y_4(\theta)>0,\quad
 \partial_{\theta_1}y_5(\theta)<0,\\
 \partial_{\theta_2}y_1(\theta)>0,\quad
 \partial_{\theta_2}y_2(\theta)<0,\quad
 \partial_{\theta_2}y_3(\theta)>0,\quad
 \partial_{\theta_2}y_4(\theta)<0,\quad
 \partial_{\theta_2}y_5(\theta)>0,\\
 \partial_{\theta_3}y_1(\theta)<0,\quad
 \partial_{\theta_3}y_2(\theta)>0,\quad
 \partial_{\theta_3}y_3(\theta)<0,\quad
 \partial_{\theta_3}y_4(\theta)<0,\quad
 \partial_{\theta_3}y_5(\theta)>0.
\end{eqnarray*}
The first claim follows from the above, see also Section~\ref{R-burgers-heat}, case (a).
The proof of the second claim is divided into three parts:
1) $\Lambda^- < \Lambda^+$,
2) $y_3^- < y_5^+ < y_2^-\,$,
and
3) $ y_2^+ < y_4^- < y_1^+\,$.


1. Changing variables, $\lambda_0=\mu+\theta_1^2/(12\,\theta_2)$, we reduce $D(P_\phi)$ of \eqref{E-discr2} to depressed form
 $P(\mu)=\mu^3+p(\theta)\,\mu+q(\theta)$,
where
\[
 p(\theta)=-\theta_1(\theta_1^3+216\,\theta_2^2\theta_3)/(48\,\theta_2^2)<0,\quad
 q(\theta)=-(\theta_1^6-540\,\theta_2^2\theta_3\theta_1^3-5832\,\theta_2^4\theta_3^2)/(864\,\theta_2^3).
\]
Due~to trigonometric solution of $P(\mu)=0$, three real roots  are
\begin{equation}\label{E-roots-mu}
 \mu_1(\theta)=A\cos\varphi>0,\quad
 \mu_2(\theta)=A\cos(\varphi-2\pi/3),\quad
 \mu_3(\theta)=A\cos(\varphi+2\pi/3)<0,
\end{equation}
where the amplitude is $A=2\,(-p/3)^{1/2}>0$ and the angle variable is given by
 $\cos(3\,\varphi) =-{4q(\theta)}/{A^3}$.
Introducing $z={\theta_3\theta_2^2}/{\theta_1^3}\in[0,\frac 1{27})$,
we obtain a decreasing (from 1 to -1) function in one variable,
 $\cos(3\,\varphi)=C(z):=-\frac{5832\,z^2+540\,z-1}{(216\,z+1)^{3/2}}$,
see Fig.~\ref{fig:3}(b).
Hence, there is a unique $\varphi=\frac 13\arccos C(z)\in[0,\frac\pi3)$.
Since $\cos(\varphi+2\pi/3)<\cos(\varphi-2\pi/3)<\cos\varphi$,
the~roots \eqref{E-roots-mu} are ordered~as
 $\mu_1(\theta)>\mu_2(\theta)>\mu_3(\theta)$.

Two positive roots $\lambda^-(\theta)<\lambda^+(\theta)$ of $D(P_\phi)$ are given by
$\lambda^-(\theta)=\mu_2(\theta)+\frac{\theta_1^2}{12\,\theta_2}$
and
$\lambda^+(\theta)=\mu_1(\theta)+\frac{\theta_1^2}{12\,\theta_2}$.
 By \eqref{E-Psi-2},
we obtain $0\le z^-\le z\le z^+<\frac 1{27}$ and $0\le\varphi^-\le\varphi\le\varphi^+<\frac\pi3$,~where
\[
 z^+ = \Psi_3^+(\Psi_2^+)^2/(\Psi_1^-)^3,\ \ z^-=\Psi_3^-(\Psi_2^-)^2/(\Psi_1^+)^3,\quad
 3\,\varphi^+ = \arccos C(z^+),\ \ 3\,\varphi^-=\arccos C(z^-).
\]
Thus, $\mu_k^-\le \mu_k(\theta)\le \mu_k^+,\ (k=1,2,3)$, where $A^\pm = 2\,(-p^\mp/3)^{1/2}$ and
$\mu_1^\pm = A^\pm\cos\varphi^\mp$,
\begin{eqnarray*}
 &&\hskip-5mm \mu_2^+ {=}\Big\{\begin{array}{c}
          \!A^+\cos\big(\varphi^+ {-}\frac{2\,\pi}3\big) \ {\rm if} \ \cos\big(\varphi^+ {-}\frac{2\,\pi}3\big)>0,\\
          \!A^-\cos\big(\varphi^+ {-}\frac{2\,\pi}3\big) \ {\rm if} \ \cos\big(\varphi^+ {-}\frac{2\,\pi}3\big)<0,
             \end{array}\
 \mu_2^- {=}\Big\{\begin{array}{c}
          \!A^-\cos\big(\varphi^- {-}\frac{2\,\pi}3\big) \ {\rm if} \ \cos\big(\varphi^- {-}\frac{2\,\pi}3\big)>0,\\
          \!A^+\cos\big(\varphi^- {-}\frac{2\,\pi}3\big) \ {\rm if} \ \cos\big(\varphi^- {-}\frac{2\,\pi}3\big)<0,
             \end{array} \\
 && p^+ = -\frac{\Psi_1^-}{48\,(\Psi_2^+)^2}\,((\Psi_1^-)^3+216\,(\Psi_2^-)^2\Psi_3^-),\quad
 p^-=-\frac{\Psi_1^+}{48\,(\Psi_2^-)^2}\,((\Psi_1^+)^3+216\,(\Psi_2^+)^2\Psi_3^+).
\end{eqnarray*}
Finally,
$\lambda^-(\theta)\le \mu_2^+ +\frac{(\Psi_1^+)^2}{12\,\Psi_2^-}$ and
$\lambda^+(\theta)\ge \mu_1^- +\frac{(\Psi_1^-)^2}{12\,\Psi_2^+}$.
To establish $\Lambda^- < \Lambda^+$, we need to show
\begin{equation}\label{E-case1-mu}
 (\Psi_1^+)^2/(12\,\Psi_2^-) - (\Psi_1^-)^2/(12\,\Psi_2^+) < \mu_1^- - \mu_2^+ .
\end{equation}
The lhs of \eqref{E-case1-mu} tends to 0, when $\delta_i\ge0$ are small enough,
while rhs of \eqref{E-case1-mu} tends to a positive constant
(estimates may be obtained using trigonometric series).
In~this case, $\Lambda^- < \Lambda^+$, and
there exists $K\in(0, (\Lambda^+ -\Lambda^-)/2)$ such that
 $D(P_\phi)$ is positive for all $\lambda_0\in(\Lambda^-+K,\,\Lambda^+-K)$ and $\theta\in{\cal P}$.

2. Consider the functions
\begin{eqnarray*}
 \phi_-(y) \eq P_{\,\phi_-}(y^2)/y^3,\ \ {\rm where}\ \
 P_{\,\phi_-}(z) = \Psi_3^- z^{3} -\lambda_0\,z^{2} +\Psi_1^- z -\Psi_2^+,\\
 \partial_y(\phi_-)(y) \eq P_{\,\partial_y(\phi_-)}(y^2)/y^4,\ \ {\rm where}\ \
 P_{\,\partial_y(\phi_-)}(z) = 3\Psi_3^- z^3 -\lambda_0 z^{2} -\Psi_1^- z +3\Psi_2^+,\\
 (\partial_y\phi)_+(y) \eq P_{\,(\partial_y\phi)_+}(y^2)/y^4,\ \ {\rm where}\ \
 P_{\,(\partial_y\phi)_+}(z) = 3\Psi_3^+ z^3 -\lambda_0 z^{2} -\Psi_1^- z +3\Psi_2^+ .
\end{eqnarray*}
It is sufficient to show that the resultant
 $R_1(t)=-{\rm Res}(P_{\,\phi_-},\,(1-t)P_{\,\partial_y(\phi_-)} +t\,P_{\,(\partial_y\phi)_+})/\Psi_2^+$
of two cubic polynomials does not vanish  for any $t\in[0,1]$ (i.e., they have no common roots).
Computation with a little help of Maple shows that $R_1(t)$ is a cubic polynomial with
 coefficients
\begin{eqnarray*}
 && a_0 = -27\,\delta_3^{3}(\Psi_2^+)^{2},\quad
    a_1 = 18\,\delta_3^{2}\big(4\,\Psi_1^-\Psi_2^+\lambda_{0}-(\Psi_1^-)^{3}-9\,(\Psi_2^+)^{2}\Psi_3^-\big), \\
 && a_2 = 12\,\delta_3 D(P_{\,\phi_-}),\quad a_3 = R_1(0) = 8\,\Psi_3^- D(P_{\,\phi_-})\,,
\end{eqnarray*}
where the discriminant $D(P_{\,\phi_-})>0$ is a cubic polynomial in $\lambda_0\in(\Lambda^-+K,\,\Lambda^+-K)$,
\[
 D(P_{\,\phi_-})= -4\,\Psi_2^+\lambda_0^3 +(\Psi_1^-)^2\lambda_0^2 +18\,\Psi_1^-\Psi_2^+\Psi_3^-\lambda_0
 -4(\Psi_1^-)^3\Psi_3^- -27\,(\Psi_2^+\Psi_3^-)^2\,.
\]
The condition \eqref{E-a3} reads as $ a_3 > |a_0| + |a_1|$, i.e.,
\[
 8\Psi_3^- D(P_{\,\phi_-}) > 27\,\delta_3^{3}(\Psi_2^+)^{2}
 +\,18\,\delta_3^2\left| 4\Psi_1^-\Psi_2^+\lambda_{0} -(\Psi_1^-)^{3} -9(\Psi_2^+)^{2}\Psi_3^-\right|\,.
\]
It is valid for small $\delta_3\ge0$ (since $0<\lambda_0\le K$).
Assuming on the contrary that either $y_2^- \le y_5^+$ or $y_3^-\ge y_5^+$, we get
$R_1(1)\le0$; hence, a contradiction: $R_1(t_0)=0$ for some $t_0\in(0,1]$.


3. Consider the functions
\begin{eqnarray*}
 \phi_+(y) \eq P_{\,\phi_+}(y^2)/y^3,\ \ {\rm where}\ \
 P_{\,\phi_+}(z) = \Psi_3^+ z^{3} -\lambda_0\,z^{2} +\Psi_1^+ z -\Psi_2^-,\\
 \partial_y(\phi_+)(y) \eq P_{\,\partial_y(\phi_+)}(y^2)/y^4,\ \ {\rm where}\ \
 P_{\,\partial_y(\phi_+)}(z) = 3\Psi_3^+ z^3 -\lambda_0 z^{2} -\Psi_1^+ z +3\Psi_2^-, \\
 (\partial_y\phi)_-(y) \eq P_{\,(\partial_y\phi)_-}(y^2)/y^4,\ \ {\rm where}\ \
 P_{\,(\partial_y\phi)_-}(z) = 3\Psi_3^- z^3 -\lambda_0 z^{2} -\Psi_1^+ z +3\Psi_2^- .
\end{eqnarray*}
It is sufficient to show that the resultant
 $R_2(t)=-{\rm Res}(P_{\,\phi_+},\,(1-t)P_{\,\partial_y(\phi_+)} +t\,P_{\,(\partial_y\phi)_-})/\Psi_2^-$
of two cubic polynomials does not vanish for any $t\in[0,1]$ (hence, they have no common roots).
Computation (again with Maple) shows that $R_2(t)$ is a cubic polynomial with coefficients
\begin{eqnarray*}
 && a_0 = 27\,\delta_3^{3}(\Psi_2^-)^{2},\quad
    a_1 = 18\,\delta_3^{2}\big(4\,\Psi_1^+\Psi_2^-\lambda_{0}-(\Psi_1^+)^{3}-9\,(\Psi_2^-)^{2}\Psi_3^+\big), \\
 && a_2 = -12\,\delta_3 D(P_{\,\phi_+}),\quad a_3 = R_2(0) = 8\,\Psi_3^+ D(P_{\,\phi_+})\,,
\end{eqnarray*}
where the discriminant $D(P_{\,\phi_+}) > 0$  is a cubic polynomial in $\lambda_0\in(\Lambda^-+K,\,\Lambda^+-K)$,
\[
 D(P_{\,\phi_+})= -4\,\Psi_2^-\lambda_0^3 +(\Psi_1^+)^2\lambda_0^2 +18\,\Psi_1^+\Psi_2^-\Psi_3^+\lambda_0
 -4(\Psi_1^+)^3\Psi_3^+ -27\,(\Psi_2^-\Psi_3^+)^2\,.
\]
The condition \eqref{E-a3} reads as $ a_3 > |a_1| + |a_2|$, i.e.,
\[
 2(3\,\Psi_3^- -\Psi_3^+)D(P_{\,\phi_+}) >
 9\,\delta_3^2\left| 4\Psi_1^+\Psi_2^-\lambda_{0} -(\Psi_1^+)^{3} -9(\Psi_2^-)^{2}\Psi_3^+\right|\,.
\]
By \eqref{E-Psi-1}, this is valid for small $\delta_3\ge0$ (since $0<\lambda_0\le K$).
Assuming on the contrary that either $y_2^+\ge y_4^-$ or $y_1^+\le y_4^-$, we get
$R_2(1)\le0$; hence, a contradiction: $R_2(t_0)=0$ for some $t_0\in(0,1]$.
\qed

\smallskip

Define closed in $C(F)$ nonempty sets
\begin{eqnarray*}
 && {\cal U}^{\;\varepsilon,\eta}=\{u_0\in C(F):\ y_2^- -\varepsilon\le u_0/e_0 \le y_2^++\eta\},\quad
 \varepsilon\in(0,\,y_2^--y_3^-),\ \eta\in(0,\,y_1^+-y_2^+).
\end{eqnarray*}
We have ${\cal U}^{\;\varepsilon,\eta}\subset\mathcal{U}_{\,1}$,
where the set ${\cal U}_{\;1}\!=\!\{\tilde u\in C(F): \ y_3^-< \tilde u/e_0 < y_1^+\}$ is open.

\begin{proposition}\label{T-new2}
Let \eqref{E-Psi-2} holds. Then

\noindent\
$(i)$ for any $u_0\in{\cal U}^{\;\varepsilon,\eta}$, Cauchy's problem \eqref{Cauchy}
has a~unique global solution of class $C^\infty(F\times(0,\infty))$,
and ${\cal U}^{\;\varepsilon,\eta}$ are invariant sets for associated semigroup
$\,\mathcal{S}_t:u_0\to u(\cdot\,,t)\ (t\ge 0)$ in~$\mathcal{C}_\infty$;

\noindent\
$(ii)$ for any $\sigma\in(0,\varepsilon)$ and $\tau\in(0,\eta)$ there is $t_1>0$ such that
$\;\mathcal {S}_t({\cal U}^{\;\varepsilon,\eta})\subseteq{\cal U}^{\;\sigma,\tau}$ for all~$t\ge t_1$.
\end{proposition}

\noindent\textbf{Proof}.
$(i)$ Let $u(\,\cdot\,,t)\ (t\ge0)$ solve (\ref{Cauchy}) with $u_0\in\mathcal{U}^{\,\varepsilon,\eta}$ for
$\varepsilon\in(0,\,y_2^- -y_3^-)$ and $\eta\in(0,\,y_1^+ -y_2^+)$.
Let $y_-(t,\varepsilon)$ and $y_+(t,\eta)$ solve the following Cauchy's problems for~ODEs, respectively:
\begin{equation*}
 y^\prime=\phi_-(y),\quad y(0)=y_2^- -\varepsilon,\qquad
 y^\prime=\phi_+(y),\quad y(0)=y_2^+ +\eta.
\end{equation*}
Since $\phi_-(y)>0$ in $(y_3^-,y_2^-)$,
the function  $y_-(t,\varepsilon)$ is increasing and $\lim\nolimits_{\,t\to\infty} y_-(t,\varepsilon)= y_2^-$.
Similarly, since $\phi_+(y)<0$ in $(y_2^+,y_1^+)$, the function
$y_+(t,\eta)$ is decreasing and $\lim\limits_{\,t\to\infty}y_+(t,\eta)=y_2^+$.

In order to apply Proposition~\ref{P-local-sol} and Theorem~\ref{prglob} to (\ref{Cauchy}), denote
\[
 {f}(u,x)=\beta(x)\,u  +\Psi_1(x)\,u^{-1}-\Psi_2(x)\,u^{-3} +\Psi_3(x)\,u^{3},
\]
and consider the closed domain $G=\{u_0\in C(F):\ y_2^- -\varepsilon_1\le u_0/e_0 \le y_2^+ +\eta_1\}$,
where $\varepsilon<\varepsilon_1<y_2^-$ and $\eta_1>\eta$, whose interior contains $\mathcal{U}^{\,\varepsilon,\eta}$.
We see that ${f}(\cdot,\cdot)\in C^\infty(G)$.
By $(i)$ and $(ii)$ of Proposition~\ref{P-local-sol}, the set $\mathbb{T}$ of
such numbers $T>0$, for which a solution $u_T(x,t)$ of Cauchy's problem \eqref{Cauchy}
exists in the domain $F\times [0,T]$, is not empty.
By Proposition~\ref{P-weak-max}, applied to
(\ref{difineqR-C}),
for any $T\in\mathbb{T}$ in the domain $F\times [0,T]$, the following inequalities are valid:
\begin{equation}\label{estwabove}
 0<y_2^- -\varepsilon\le y_-(t,\varepsilon)\le w_T(\,\cdot\,,t)\le y_+(t,\eta)\le y_2^++\eta,
\end{equation}
where $w_T(x,t)=u_T(x,t)/e_0(x)$. By Theorem~\ref{prglob},
the solution $u(x,t)$ of (\ref{Cauchy}) exists for all $(x,t)\in{\cal C}_\infty$,
$u(\cdot,\cdot)\in C^\infty(F\times(0,\infty))$ and the set ${\cal U}^{\;\varepsilon,\eta}$
is invariant for operators ${\cal S}_t\ (t\ge 0)$, that proves $(i)$.
Claim $(ii)$ follows immediately from (\ref{estwabove}).
\qed

\smallskip
By \eqref{E-auxcor}, we have $y_2^->y_5^+>y_3^-$ and
$y_2^+<y_4^+<y_1^+$. Define the following quantity:
\begin{equation}\label{E-mu-pm-3}
 \mu^+(\sigma,\tau):=-\max\nolimits_{\,y\in (y_2^-
 -\sigma,\,y_2^++\tau)}\,(\partial_y\phi)_+(y)>0
\end{equation}
for $\sigma\in(0,\,y_2^- -y_5^+)$ and $\tau\in(0,\,y_4^+-y_2^+)$.

\begin{theorem}\label{thattract2}
\ $(i)$ If \eqref{E-Psi-2} holds then \eqref{Cauchy-stat1} has a
solution $u_*\in\mathcal{U}_{\,1}\cap C^\infty(F)$; moreover, the
set $U_*$ of all such solutions is compact is $C(F)$ and
$U_*\subset\{u_0\in C(F): \ y_2^- \le \tilde u/e_0 \le y_2^+\}$.

\ $(ii)$ If,~in addition, \eqref{E-Psi-1} holds and
$\delta_i=\Psi^+_i-\Psi^-_i\ (i=1,2,3)$ are small enough then the
above solution is unique in $\mathcal{U}_{\,1}$, and
$u_*=\lim\nolimits_{\,t\to\infty} u(\cdot\,,t)$, where $u$
solves~\eqref{Cauchy} with $u_0\in\mathcal{U}_{\,1}$; moreover, for
any $\sigma\in(0,\,y_2^- -y_5^+)$ and $\tau\in(0,\,y_4^+-y_2^+)$,
the set $\mathcal{U}^{\,\sigma,\tau}$ is attracted by associated
semigroup exponentially fast to $u_*$ in~$C$-norm:
\begin{equation}\label{rateconv2}
 \| u(\cdot\,,t)-u_*\|_{C(F)} \le \delta^{-1}(e_0)\,e^{-\mu^+(\sigma,\tau)\,t}\| u_0-u_*\|_{C(F)}
 \quad(t>0,\ u_0\in\mathcal{U}^{\,\sigma,\tau}).
\end{equation}
\ $(iii)$ Let $\beta,\Psi_1,\Psi_2,\Psi_3$ be smooth functions
on $F\times\RR^n$ with a smooth metric $g(\cdot,q)$.
If \eqref{E-Psi-2}, \eqref{E-case1-mu}, \eqref{E-Psi-1} and
\eqref{E-case1-P} hold for any $F\times\{q\}\ (q\in\RR^n)$ then the
unique solution $u_*$, see $(ii)$, is smooth on $F\times\RR^n$.
\end{theorem}

\noindent\textbf{Proof}. (i) By Proposition~\ref{T-new2}(i), the set
$\mathcal{U}^{\,\varepsilon,\eta}$ is invariant for the semigroup
$\,\mathcal{S}_t\ (t\ge 0)$ corresponding to (\ref{Cauchy})$_1$,
i.e., ${\cal S}_t\big(\mathcal{U}^{\,\varepsilon,\eta}\big)\subseteq\mathcal{U}^{\,\varepsilon,\eta}\
(t\ge 0)$. By Theorem~\ref{P-exist} with $u_-=y_3^- e_0$, $u_+=y_1^+
e_0$, $G={\cal U}^{\;\varepsilon,\eta}$ and
\[
 f(u,x)=\beta\,u + \Psi_1(x)\,u^{-1} -\Psi_2(x)\,u^{-3} +\Psi_3(x)\,u^{3},
\]
the set $U_*^{\;\varepsilon,\eta}$ of all solutions of \eqref{Cauchy-stat1}
lying in ${\cal U}^{\;\varepsilon,\eta}$ is nonempty and compact in $C(F)$.
Since the intersection of any finite subfamily of the family of
compact sets $\{U_*^{\;\varepsilon,\eta}\}_{\varepsilon,\eta>0}$ is nonempty and
compact in $C(F)$, thus the whole family has nonempty and compact in
$C(F)$ intersection $U_*$.

(ii) To prove the second claim, take initial values
$u_i^0\in\mathcal{U}^{\,\sigma,\tau}\ (i=1,2)$ with
$\sigma\in(0,\,y_2^- -y_5^+)$ and $\tau\in(0,\,y_4^+-y_2^+)$,
and denote~by
\[
 u_i(\cdot\,,t)={\cal S}_t(u_i^0),\quad
 w_i(\,\cdot,t)=u_i(\,\cdot,t)/e_0,\quad
 w_i^0=u_i^0/e_0.
\]
 From (\ref{Cauchw1R-C}), using the equalities
\[
 2\,\bar w\,\Delta\,\bar w=\Delta(\bar w^2)-2\,\|\nabla\,\bar w\|^2,\quad
 \nabla\,(\bar w^2)=2\,\bar w\nabla\,\bar w
\]
with $\bar w=w_2-w_1$, we obtain
\begin{eqnarray*}
 &&\partial_t\big((w_2-w_1)^2\big)=2\,(w_2-w_1)\,\partial_t(w_2-w_1)\le\Delta\big((w_2-w_1)^2\big)\\
 &&+\,\<2\,\nabla\log e_0,\,\nabla (w_2-w_1)^2\> +2\,(f(w_2,\,\cdot\,)-f(w_1,\,\cdot\,))(w_2-w_1).
\end{eqnarray*}
Observe that in view of (\ref{Cauchw1R-C}), and (\ref{estderivphi1}), for all $x\in F$ we have
\begin{equation}\label{estderivf2}
 (\partial_w\phi)_-(w)\le\partial_w f(w,x)\le(\partial_w\phi)_+(w).
\end{equation}
We estimate the last term, using $y_5^+<y_2^- -\sigma\le w_i\le y_2^++\tau<y_4^+\;(i=1,2)$,
\eqref{E-mu-pm-3} and the right inequality of (\ref{estderivf2}):
\begin{eqnarray*}
 &&(f(w_2,\,\cdot\,)-f(w_1,\,\cdot\,))(w_2-w_1)\\
 &&=(w_2-w_1)^2\!\int_0^1\partial_w f(w_1+\tau(w_2-w_1),\,\cdot\,)\dtau \le-\mu^+(\sigma,\tau)(w_2-w_1)^2.
\end{eqnarray*}
Thus, $v=(w_2-w_1)^2$ satisfies the differential inequality
 $\partial_t v \le\Delta\,v +\<2\,\nabla\log e_0,\ \nabla v\> -2\,\mu^+(\sigma,\tau)\,v$.
By Proposition~\ref{P-weak-max}, we obtain
$v(\,\cdot\,,t)\le v_+(t)$, where $v_+(t)$ solves the Cauchy's problem for~ODE:
\[
 v\,_+'=-2\,\mu^+(\sigma,\tau)\,v_+(t),\quad v_+(0)=\| w_2^0-w_1^0\|_{C(F)}^2.
\]
Thus,
\begin{eqnarray}\label{cndcontr}
\nonumber
 &&\|\mathcal{S}_t(u_2^0)-\mathcal{S}_t(u_1^0)\|_{C(F)}
 \le \| w_2(\,\cdot\,,t)-w_1(\,\cdot\,,t)\|_{C(F)} \cdot\max_{\,F}e_0\\
 &&\le e^{-\mu^+(\sigma,\tau)\,t}\,\| w_2^0-w_1^0\|_{C(F)}\cdot\max_{\,F}e_0
 \le \delta^{-1}(e_0) e^{-\mu^+(\sigma,\tau)\,t}\,\| u_2^0-u_1^0\|_{C(F)}\,,
\end{eqnarray}
i.e., the operators $\mathcal{S}_t\ (t\ge 0)$ for (\ref{Cauchy})$_1$ satisfy
in $\mathcal{U}^{\,\sigma,\tau}$, where $\sigma\in(0,\,y_2^--y_5^+)$ and $\tau\in(0,\,y_4^+-y_2^+)$,
the Lipschitz condition for $C$-norm with the Lipschitz constant
$\delta^{\,-1}(e_0) e^{\,-\mu^+(\sigma,\tau)\,t}$.

By~Proposition~\ref{T-new2}(i), each operator ${\cal S}_t\ (t\ge 0)$
for (\ref{Cauchy})$_1$ maps the set $\mathcal{U}^{\,\sigma,\tau}$, which is closed in $C(F)$,
into itself and, by the above arguments, for $t>\frac 1{\mu^+(\sigma,\tau)}\ln \delta^{-1}(e_0)$
it is a contraction there. Since all operators ${\cal S}_t$ commute one with another,
they have a unique common fixed point $u_*$ in $\mathcal{U}^{\,\sigma,\tau}$ and,
in view of (\ref{cndcontr}), the inequality (\ref{rateconv2}) holds
for any $u_0\in\mathcal{U}^{\,\sigma,\tau}$ and $t\ge 0$.

On the other hand, by Proposition~\ref{T-new2}(ii), if
$\varepsilon\in(y_2^- -y_5^+,\,y_2^- -y_3^-)$,
$\eta\in(y_4^+-y_2^+,\,y_1^+-y_2^+)$, $\sigma\in(0,\,y_2^--y_5^+)$
and $\tau\in(0,\,y_4^+-y_2^+)$ then $\sigma<\varepsilon$, $\tau<\eta$ and
$\;{\cal
S}_t\big(\mathcal{U}^{\,\varepsilon,\eta}\big)\subseteq\mathcal{U}^{\,\sigma,\tau}$
for some $t_1>0$ and any $t\ge t_1$. Hence, $u_*$ is a unique
fixed point of the operators ${\cal S}_t$ also in the sets
$\mathcal{U}^{\,\varepsilon,\eta}$ with $\varepsilon\in(0,\,y_2^--y_3^-)$ and
$\eta\in(0,\, y_1^+-y_2^+)$. Since $\varepsilon$ and $\eta$ are arbitrary
in the corresponding intervals, $u_*$ is a unique fixed point of
${\cal S}_t$ in the wider set $\mathcal{U}_1$; moreover, $y_2^-\le
u_*/e_0\le y_2^+$. By the above arguments,
$u_*=\lim\limits_{\,t\to\infty} u(\cdot\,,t)$, where $u$
solves~\eqref{Cauchy} with $u_0\in\mathcal{U}_{\,1}$. Thus, in view
of Proposition~\ref{P-local-sol}, $u_*$ is a solution of
(\ref{Cauchy-stat1}) belonging to $C^\infty(F)$.


(iii) Let $e_0(x,q)>0$ be the normalized eigenfunction for the
minimal eigenvalue $\lambda_0(q)$ of
$\mathcal{H}_q=-\Delta-\beta(x,q)$. By Theorem~\ref{contdeponqeigfunct},
$\lambda_0\in C^\infty(\RR^n)$ and $e_0\in
C^\infty(F\times\RR^n)$, hence $y_3^-$ and $y_1^+$
smoothly depend on $q$. As we have proved in (ii), for any
$q\in\RR^n$ the stationary equation,
\begin{equation}\label{nonlinstatpar1}
 \Delta_q\,u + f(u,x,q)=0,
\end{equation}
see also~(\ref{Cauchy-stat1}),
where $f(u,x,q)=\beta(x,q)u+\Psi_1(x,q)\,u^{-1}
-\Psi_2(x,q)\,u^{-3} +\Psi_3(x,q)\,u^{3}$ has a unique solution
$u_*(x,q)$ in the open set ${\mathcal U}_{\;1}(q)=\{u_0\in
C(F\times\RR^n):\,y_3^-(q)< u_0/e_0(\cdot,q)<y_1^+(q)\}$.

Since $y_3^-(q)$, $y_1^+(q)$  and $e_0(x,q)$ are continuous,
for any $k\in\NN$ and $\alpha\in(0,1)$,
there exist open neighborhoods $U_*\subseteq C^{k+2,\alpha}$ of
$u_*(x,0)$ and $V_0\subset\RR^n$ of $0$ such that
\begin{equation}\label{inclneighb}
  U_*\subseteq{\mathcal U}_{\;1}(q)\quad \forall\, q\in V_0.
\end{equation}
We claim that all eigenvalues of the linear operator
 $\mathcal{H}_*=-\Delta_{0}-\partial_u\,f(u_*(x,0),x,0)$,
acting in $L_2$ with the domain $H^2$, are positive.
To show this, observe that $y_2^-(0)\le
u_*(\cdot,0)/e_0(\cdot,0)\le y_2^+(0)$.
Let $\tilde u(x,t)$ be a solution of Cauchy's problem for the evolution equation
\begin{equation}\label{Cauchstar}
 \partial_t\tilde u=-\mathcal{H}_*(\tilde u),\quad \tilde u(x,0)=\tilde u_0(x)\in C(F).
\end{equation}
Using the same arguments as in the proof of (ii), we obtain that
$v(x,t)={\tilde u^{\,2}(x,t)}\,{e^{-2}_0(x,0)}$ obeys
the differential inequality
 $\partial_t v \le\Delta_{0}\,v +\<2\,\nabla\log e_0(\cdot,0),\ \nabla v\> -2\,\mu^+_0\,v$
with $\mu^+_0>0$. By~Proposition~\ref{P-weak-max},
$v(\,\cdot\,,t)\le v_+(t)$, where $v_+(t)$ solves the Cauchy's problem for~ODE
\[
 v\,_+'=-2\,\mu^-_0\,v_+,\quad v_+(0)=\|{\tilde u_0}/{e_0(\cdot,0)}\|_C^2\,;
\]
moreover, for any $\tilde u_0\in C(F)$ the function $\tilde u(x,t)$
tends to $0$ exponentially fast, as $t\rightarrow\infty$. On the
other hand, if $\tilde\lambda_\nu$ is any eigenvalue of
$\mathcal{H}_*$ and $\tilde e_\nu(x)$ is the corresponding
normalized eigenfunction then $\tilde u=e^{-\tilde\lambda_\nu
t}\tilde e_\nu$ solves (\ref{Cauchstar}) with $\tilde
u_0(x)=\tilde e_\nu(x)$. Thus, $\tilde\lambda_\nu>0$ that
completes the proof of the claim.

By Theorem~\ref{prsmoothparstat}, for any integers $k\ge 0$ and $l\ge 1$ we can restrict the neighborhoods
$U_*$ of $u_*(x,0)$ and $V_0$ of $0$ in such a way that

-- for any $q\in V_0$ there exists in $U_*$ a unique solution $\tilde u(x,q)$
of (\ref{nonlinstatpar1}), and

-- the mapping $q\rightarrow \tilde u(\cdot,q)$ belongs to class $C^l(V_0,U_*)$.

\noindent
In view of (\ref{inclneighb}), $\tilde u(\cdot,q)=u_*(\cdot,q)$ holds for any $q\in V_0$.
\qed


\begin{remark}\label{R-04}\rm
Similarly, for $\Psi_2\equiv0$ when $\Psi_3>0$, $\Psi_1>0$, and $\lambda_0>0$,
we have $P_\phi(z)=\theta_3\,z^{2} -\lambda_0\,z +\theta_1$.
The function $\phi(y)=P_\phi(y^2)/y$ has two positive roots $y_{2}(\theta)<y_{1}(\theta)$,
$y_{1,2}^2(\theta)=\frac{\lambda_0\pm(\lambda_0^2-4\theta_1\theta_3)^{1/2}}{2\theta_3}$,
when
\begin{equation}\label{E-rem4}
 \lambda_0^2 > 4\Psi_1^+ \Psi_3^+.
\end{equation}
Note that $y_{2}$ decreases in $\lambda_0$.
Its derivative $\partial_y\phi(y)$ has one positive root $y_{4}(\theta)\in(y_{2}(\theta),y_{1}(\theta))$;
moreover, $\partial_y\phi\,\vert_{\,y=y_1(\theta)}>0$, $\partial_y\phi\,\vert_{\,y=y_2(\theta)}<0$.
Consider the functions
\begin{eqnarray*}
 \phi_+(y) \eq P_{\,\phi_+}(y^2)/y,\ \ {\rm where}\ \
 P_{\,\phi_+}(z) = \Psi_3^+ z^{2} -\lambda_0\,z +\Psi_1^+,\\
 \phi_-(y) \eq P_{\,\phi_-}(y^2)/y,\ \ {\rm where}\ \
 P_{\,\phi_-}(z) = \Psi_3^- z^{2} -\lambda_0\,z +\Psi_1^-,\\
 \partial_y(\phi_+)(y) \eq P_{\,\partial_y(\phi_+)}(y^2)/y^2,\ \ {\rm where}\ \
 P_{\,\partial_y(\phi_+)}(z) = 3\Psi_3^+ z^2 -\lambda_0 z -\Psi_1^+ ,\\
 (\partial_y\phi)_-(y) \eq P_{\,(\partial_y\phi)_-}(y^2)/y^2,\ \ {\rm where}\ \
 P_{\,(\partial_y\phi)_-}(z) = 3\Psi_3^- z^2 -\lambda_0 z -\Psi_1^+ .
\end{eqnarray*}
Denote by $y_2^+ < y_1^+$ the positive roots of $\phi_+(y)$, by $y_2^- < y_1^-$ the positive roots of $\phi_-(y)$,
and $y_4^-$ the positive root of $(\partial_y\phi)_+(y)$.
Then \eqref{E-auxcor} reduces to
\begin{equation}\label{E-auxcor-0}
 y_2^- < y_2^+ < y_4^- < y_1^+ < y_1^-\,.
\end{equation}
To find sufficient conditions for this, we will show that the resultant of two quadratic polynomials
 $R_2(t)=-{\rm Res}(P_{\,\phi_+},\,(1-t)P_{\,\partial_y(\phi_+)} +t\,P_{\,(\partial_y\phi)_-})/\Psi_3^+$
does not vanish for any $t\in[0,1]$; hence, the polynomials have no common roots.
 Thus, $R_2(t)=a_0t^{2}+a_{1}t+a_2$ is a quadratic polynomial with
\begin{equation*}
  a_0 =  -9\,\delta_3^2\Psi^+_1,\quad
  a_1 = -6\,\delta_3(\lambda_0^2 -4\,\Psi_1^+\Psi_3^+), \quad
  a_2 = R_2(0) = 4\,\Psi^+_3(\lambda_0^2 -4\,\Psi_1^+\Psi_3^+)\,,
 \end{equation*}
Note that $D(P_{\,\phi_+}) = \lambda_0^2 -4\,\Psi_1^+\Psi_3^+ > 0$.
Hence, \eqref{E-a3} reads as $ a_2 > |a_1| + |a_1|$, i.e.,
\begin{equation*}
 2(3\,\Psi^-_3 - \Psi^+_3)(\lambda_0^2- 4\,\Psi_1^+\Psi_3^+) > 9\,\delta_3^2\Psi^+_1.
\end{equation*}
We conclude that \eqref{E-auxcor-0} follows from the inequalities \eqref{E-Psi-1} and
\[
 \lambda_0^2 > \Psi_1^+\big(4\Psi_3^+ + (9/2)\,\delta_3^2/(3\,\Psi_3^- -\Psi_3^+)\big)
 =\Psi_1^+\frac{(3\,\Psi_3^- +\Psi_3^+)^2}{2(3\,\Psi_3^- -\Psi_3^+)}\,.
\]
Note that the last inequality yields \eqref{E-rem4}.
\end{remark}

\subsubsection{Case of $\Psi_3<0$}
\label{sec:attr}

Let $\Psi_3<0$, $\Psi_1<0$ and $\Psi_2>0$ and $\lambda_0<0$, see Section~\ref{R-burgers-heat}, case (b).
Consider the function  for $y>0$
\[
 \phi(y,\theta) = -\lambda_0\,y -\theta_1\,y^{-1}-\theta_2\,y^{-3} -\theta_3\,y^{3} = P_{\,\phi}(y^2)/y^3,
\]
where $P_{\,\phi}(z) = -\theta_3\,z^{3} -\lambda_0\,z^{2} -\theta_1\,z-\theta_2$
and $\theta=(\theta_1,\theta_2,\theta_3)\in\mathcal{P}$.
Then $\phi_-(y)\le\phi(y,\theta)\le\phi_+(y)$,~where
\begin{eqnarray*}
 \phi_+(y) \eq P_{\,\phi_+}(y^2)/y^3,\ \ {\rm where}\ \
 P_{\,\phi_+}(z) = -\Psi_3^- z^{3} -\lambda_0\,z^{2} -\Psi_1^- z -\Psi_2^-,\\
 \phi_-(y) \eq P_{\,\phi_-}(y^2)/y^3,\ \ {\rm where}\ \
 P_{\,\phi_-}(z) = -\Psi_3^+ z^{3} -\lambda_0\,z^{2} -\Psi_1^+ z -\Psi_2^+.
\end{eqnarray*}
 We calculate
\[
 \partial_y\phi=-\lambda_0+\theta_1\,y^{-2}+3\,\theta_2\,y^{-4} -3\,\theta_3\,y^2,\quad
 \partial^2_{yy}\phi=-2\,\theta_1\,y^{-3}-12\,\theta_2\,y^{-5} -6\,\theta_3\,y.
\]
Since $\partial^2_{yy}\phi<0$ for $y>0$ and $\phi(0+,\,\theta)=\phi(\infty,\,\theta)=-\infty$,
the function $\phi$ is concave by $y$ and ``$\cap\,$"-shaped,
and $\partial_y\phi$ is decreasing from $\infty$ to $-\infty$ for $y\in(0,\infty)$.
Note that $\phi_-(y)$ and $\phi_+(y)$ are also concave.
The discriminant of $P_\phi(z)$ is the following cubic polynomial in $-\lambda_0$:
\begin{eqnarray*}
 D(P_\phi) \eq 4\,\theta_2(-\lambda_0)^3 +\theta_1^2(-\lambda_0)^2 -18\,\theta_1\theta_2\theta_3\,(-\lambda_0)
 -(4\,\theta_1^3\theta_3 +27\,\theta_2^2\theta_3^2)\\
 &&\hskip-6mm\ge \bar D := 4\,\Psi^-_2(-\lambda_0)^3 +(\Psi_1^-)^2(-\lambda_0)^2 -18\,\Psi_1^+\Psi_2^+\Psi_3^+(-\lambda_0)
 -4(\Psi_1^+)^3\Psi_3^+ -27(\Psi_2^+\Psi_3^+)^2\,.
\end{eqnarray*}
By Maclaurin method, the following condition is sufficient for $\bar D>0$:
\begin{equation}\label{E-ineqL0}
 \lambda_0 < -\bar K,\quad
 \bar K=1+\big(\max\{18\,\Psi_1^+\Psi_2^+,\ 4(\Psi_1^+)^3 +27(\Psi_2^+)^2\Psi_3^+\,\}\Psi_3^+/(4\Psi_2^-)\big)^{1/2}.
\end{equation}
By the above, if \eqref{E-ineqL0} holds then $\phi(y,\theta)$ for any $\theta\in\mathcal{P}$ has two positive roots $y_2(\theta)>y_1(\theta)$, and $\partial_y\phi$ has a unique positive root $ y_3(\theta)\in(y_1(\theta),\,y_2(\theta))$.
Note that $\partial_y\phi\,\vert_{\,y=y_2(\theta)}<0$ and  $\partial_y\phi\,\vert_{\,y=y_1(\theta)}>0$.

Let
$y_1^+ < y_2^+$ be positive roots of $\phi_+(y)$, $y_1^- < y_2^-$ the positive roots of $\phi_-(y)$,
and $y_3^-,\;y_3^+$ positive roots of decreasing functions
\begin{equation*}
 (\partial_y\phi)_-(y)=-\lambda_0  +\Psi_1^- y^{-2} +3\Psi_2^- y^{-4}-3\Psi_3^+y^2,\quad
 (\partial_y\phi)_+(y)=-\lambda_0 +\Psi_1^+ y^{-2}+3\Psi_2^+ y^{-4} -3\Psi_3^-y^2.
\end{equation*}
Note that $(\partial_y\phi)_-(y)\le\partial_y\phi(y,\theta)\le(\partial_y\phi)_+(y)$ for all $\theta\in{\cal P}$ and $y>0$.

\begin{proposition}\label{P-auxprop}
If \eqref{E-ineqL0} holds then for any $\theta\in\mathcal{P}$,
\[
 y_1^+ \le y_1(\theta) \le y_1^-,\quad
 y_2^- \le y_2(\theta) \le y_2^+,\quad
 y_3^- \le y_3(\theta) \le y_3^+\,.
\]
If, in addition, \eqref{E-Psi-1} and
\begin{eqnarray}\label{E-Kb}
 &&-\lambda_0 > 1+\sqrt K,\quad
 {\rm where}\\
\nonumber
 &&\hskip-6mm K=\frac{\max\,\{\,36\,\Psi_1^+\Psi_2^+\Psi_3^-(\Psi_3^- +\Psi_3^+),\,
 27\,\Psi_3^+(\Psi_2^+)^2((\Psi_3^+)^2 + 3(\Psi_3^-)^2) +(\Psi_1^+)^3(\Psi_3^+ +3\Psi_3^-)^2\,\}}
 {8\,\Psi_2^+(3\,\Psi_3^- -\Psi_3^+)},
\end{eqnarray}
hold then there exist $K>\bar K$ such that for all $\lambda_0<-K$ we have
\begin{equation}\label{E-auxcor-01}
 y_1^+ < y_1^- < y_3^+ < y_2^- < y_2^+\,.
\end{equation}
\end{proposition}

\noindent\textbf{Proof}.
For implicit derivatives
 $\partial_{\theta_k}y_l = -({\partial_{\theta_k}\phi}/{\partial_y\phi})\,\vert_{\,y=y_l(\theta)}$,
 $\partial_{\theta_k}y_3 = -({\partial\,^2_{\theta_k y}\phi}/{\partial\,^2_{yy}\phi})\,\vert_{\,y=y_3(\theta)}$
 where $l=1,2,\;k=1,2,3$, we calculate
\begin{eqnarray*}
 &&\partial_{\theta_1}\phi = -y^{-1},\quad
 \partial_{\theta_2}\phi = -y^{-3},\quad
 \partial_{\theta_3}\phi = -y^3,\quad
 \partial_{y}\phi\vert_{y=y_2(\theta)} <0,\quad \partial_{y}\phi\vert_{y=y_1(\theta)} >0, \qquad\\
 &&\partial^2_{\theta_1 y}\phi = y^{-2},\quad
 \partial^2_{\theta_2 y}\phi = 3y^{-4},\quad
 \partial^2_{\theta_3 y}\phi = -3y^2.
\end{eqnarray*}
Recall that ${\partial\,^2_{yy}\phi}<0\ (y>0)$.
Thus, the following inequalities hold:
\begin{eqnarray*}
 \partial_{\theta_k}y_1(\theta)>0,\quad
  \partial_{\theta_k}y_2(\theta)<0,\quad (k=1,2,3),\quad
 \partial_{\theta_k}y_3(\theta)>0\quad (k=1,2),\quad
 \partial_{\theta_3}y_3(\theta)<0.
\end{eqnarray*}
The first claim follows from the above, see also Section~\ref{R-burgers-heat}, case (b).
 For the second claim, is is sufficient find $K>\bar K$ such that for all $\lambda_0 < -K$ we have
$y_1^- < y_3^+ < y_2^-$. Consider the functions
\begin{eqnarray*}
 \partial_y(\phi_-)(y) \eq P_{\,\partial_y(\phi_-)}(y^2)/y^4,\ \ {\rm where}\ \
 P_{\,\partial_y(\phi_-)}(z) = -3\Psi_3^+ z^3 -\lambda_0 z^{2} +\Psi_1^+ z +3\Psi_2^+,\\
 (\partial_y\phi)_+(y) \eq P_{\,(\partial_y\phi)_+}(y^2)/y^4,\ \ {\rm where}\ \
 P_{\,(\partial_y\phi)_+}(z) = -3\Psi_3^- z^3 -\lambda_0 z^{2} +\Psi_1^+ z +3\Psi_2^+ ,
\end{eqnarray*}
for $y>0$, where $\partial_y(\phi_-)$ and $(\partial_y\phi)_+$ are decreasing.
 Notice that $\phi_-(y)>0$ for $y\in(y_1^-,\,y_2^-)$,
and $\phi_-(y)<0$ for $y\in(0,\,\infty)\setminus[y_1^-,\,y_2^-]$,
and we have  $\phi_-(0+)=-\infty$ and  $\phi_-(\infty)=-\infty$;
moreover, $\phi_-(y)$ increases in $(0,y_3^-)$ and decreases in $(y_3^-,\,\infty)$.
The~function $(\partial_y\phi)_+(y)$ decreases on $(0,\infty)$ from $+\infty$ to $-\infty$;
moreover, $(\partial_y\phi)_+(y)>0$ in $(0,y_3^+)$ and $(\partial_y\phi)_+(y)<0$ in $(y_3^+,\infty)$.

 Since, the positive root of $\partial_y(\phi_-)$ belongs to $(y_1^-, y_2^-)$,
we will show that the resultant of two cubic polynomials
 $R_3(t)=-{\rm Res}(P_{\,\phi_-},\,(1-t)P_{\,\partial_y(\phi_-)} +t\,P_{\,(\partial_y\phi)_+})/\Psi_2^+$
does not vanish (hence, they have no common roots) for any $t\in[0,1]$.
Indeed,
 $R_3(0) = 8\,\Psi_3^+\,D(P_{\,\phi_-}) \ge 8\,\Psi_3^+\,\bar D(-\lambda_0)>0$,
where
\[
 D(P_{\,\phi_-})= -4\,\Psi_2^+\lambda_0^3 +(\Psi_1^+)^2\lambda_0^2 +18\,\Psi_1^+\Psi_2^+\Psi_3^+\lambda_0
 -4(\Psi_1^+)^3\Psi_3^+ -27\,(\Psi_2^+\Psi_3^+)^2.
\]
Assuming on the contrary that either $y_3^+\ge y_2^-$ or $y_3^+\le y_1^-$, we get
$R_3(1)\le0$; hence, a contradiction: $R_3(t_0)=0$ for some $t_0\in(0,1]$.
In our case, $R_3(t)$ is a cubic polynomial with coefficients
\begin{eqnarray*}
 && a_0=27\,\delta_3^{3}(\Psi_2^+)^{2},\quad
 a_1 = -18\,\delta_3^{2}\big( 4\,\Psi_1^+\Psi_2^+(-\lambda_{0}) +(\Psi_1^+)^{3} +9\,(\Psi_2^+)^{2}\Psi_3^+\big), \\
 && a_2 = -12\,\delta_3 D(P_{\,\phi_-}),\quad
 a_3= 8\,\Psi_3^+ D(P_{\,\phi_-})\,.
\end{eqnarray*}
Hence, the condition \eqref{E-a3} reads as $ a_3 > |a_1| + |a_{2}|$ (since $a_0>0$), i.e.,
\begin{equation}\label{cndineq}
 2\,(3\,\Psi_3^- -\Psi_3^+) D(P_{\,\phi_-}) > 9\,\delta_3^2\left( 4\Psi_1^+\Psi_2^+(-\lambda_{0})
 +(\Psi_1^+)^{3} +9(\Psi_2^+)^{2}\Psi_3^+\right)\,.
\end{equation}
By \eqref{E-Psi-1}, this is valid if either $\delta_3\ge0$ is small or $P(-\lambda_0)=\sum_{i=0}^3 b_{3-i}\,(-\lambda_0)^i$
is positive, where
\begin{eqnarray*}
 b_0 \eq 8\,\Psi_2^+(3\,\Psi_3^- -\Psi_3^+),\quad
 b_1 = 2\,(\Psi_1^+)^{2}(3\,\Psi_3^- -\Psi_3^+),\quad
 b_2 =  -36\,\Psi_1^+\Psi_2^+\Psi_3^-(\Psi_3^- +\Psi_3^+)<0,\\
 b_3 \eq -27\,\Psi_3^+(\Psi_2^+)^2((\Psi_3^+)^2 + 3(\Psi_3^-)^2) -(\Psi_1^+)^3(\Psi_3^+ +3\Psi_3^-)^2 <0\,.
\end{eqnarray*}
By Maclaurin method, the inequality $-\lambda_0 > K$, where
 $K=1+(\max\,\{-b_2,\,-b_3\}/\,b_0)^{1/2}$,
yields $P(-\lambda_0)>0$
(if $\delta_3\ge0$ is small enough,
then one may take $K=\bar K$).
\qed

\smallskip

Define closed in $C(F)$ nonempty sets
\begin{eqnarray*}
 && {\cal U}^{\;\varepsilon,\eta}=\{\tilde u\in C(F):\ y_2^- -\varepsilon\le \tilde u/e_0 \le y_2^+ +\eta\},\quad
 \varepsilon\in(0,\,y_2^- -y_1^-),\ \eta\in(0,\,\infty].
\end{eqnarray*}
We have ${\cal U}^{\;\varepsilon,\eta}\subset\mathcal{U}^{\,\varepsilon,\infty}\subset\mathcal{U}_{\,1}$,
where the set ${\cal U}_{\;1}\!=\!\{\tilde u\in C(F): \tilde u/e_0 > y_1^-\}$ is open.
 The proof of the following proposition and theorem is similar to the proof of Proposition~\ref{T-new2}
and Theorem~\ref{thattract2}.

\begin{proposition}\label{T-new1}
Let \eqref{E-ineqL0} holds. Then


$(i)$ for any $u_0\in{\cal U}^{\;\varepsilon,\eta}$, Cauchy's problem
\eqref{Cauchy} admits a~unique global solution. Moreover, ${\cal
U}^{\;\varepsilon,\eta}$ are invariant sets for the associated semigroup
$\,\mathcal{S}_t:u_0\to u(\cdot\,,t)\ (t\ge 0)$ in~${\cal
C}_\infty$;


$(ii)$ for any $\sigma\in(0,\varepsilon)$ there exists $t_1>0$  such that
$\;\mathcal {S}_t({\cal U}^{\;\varepsilon,\infty})\subseteq{\cal U}^{\;\sigma,\infty}$
for all $t\ge t_1$.
\end{proposition}

By \eqref{E-auxcor-01}, we have $y_2^- - y_3^+ > 0$. Define the following quantity for $\sigma\in(0,\,y_2^- -y_3^+)$:
\begin{equation*}
 \mu^+(\sigma):=-\sup\nolimits_{\,y\ge y_2^- -\sigma}\,(\partial_y\phi)_+(y)=-(\partial_y\phi)_+(y_2^- -\sigma)>0\,.
\end{equation*}

\begin{theorem}\label{thattract1}
\ $(i)$ If \eqref{E-ineqL0} holds then \eqref{Cauchy-stat1} has a
solution $u_*\in\mathcal{U}_{\,1}\cap C^\infty(F)$; moreover, the
set $U_*$ of all such solutions is compact is $C(F)$ and
$U_*\subset\{\tilde u\in C(F): \ y_2^- \le \tilde u/e_0 \le y_2^+\}$.

\ $(ii)$ If,~in addition, $\Psi_3^+ < 3\,\Psi_3^-$ then there exists
$K>\bar K$ such that if $\lambda_0 < -K$ then the above solution is
unique in $\mathcal{U}_{\,1}$, and $u_*=\lim\limits_{\,t\to\infty}
u(\cdot\,,t)$ where $u$ solves~\eqref{Cauchy} with
$u_0\in\mathcal{U}_{\,1}$; moreover, for any $\sigma\in(0,\,y_2^-
-y_3^+)$, the set $\,\mathcal{U}^{\,\sigma,\infty}$ is attracted by
the corresponding semigroup exponentially fast to the point
$u_*$ in $C$-norm:
\begin{equation*}
 \| u(\cdot\,,t)-u_*\|_{C(F)} \le \delta^{-1}(e_0)\,e^{-\mu^+(\sigma)\,t}\| u_0-u_*\|_{C(F)}
 \quad(t>0,\ u_0\in\mathcal{U}^{\,\sigma,\infty}).
\end{equation*}
\ $(iii)$ Let $\beta,\Psi_1,\Psi_2,\Psi_3$ be smooth functions
on
$F\times\RR^n$ with a smooth metric $g(\cdot,q)$.
If~\eqref{E-Psi-1}, \eqref{E-ineqL0} and \eqref{E-Kb} hold for any
$F\times\{q\}\ (q\in\RR^n)$ then the solution $u_*$, see $(ii)$, is smooth on~$F\times\RR^n$.
\end{theorem}

\begin{remark}\label{R-03}\rm
Let $\Psi_2\equiv0$ when $\Psi_3<0$, $\Psi_1\le0$ and $\lambda_0<0$.
Due to geometric definition \eqref{E-Psi-geom} of $\Psi_i$ in \eqref{E-Yam1-init}, we are forced to assume $\Psi_1=0$.
Then we have $P_\phi(z)=-\theta_3 z^{2} -\lambda_0 z$, and for $\lambda_0 < 0$
the function $\phi(y)=P_\phi(y^2)/y$ has one positive root $y_{1}(\theta)=({-\lambda_0}/{\theta_3})^{1/2}$,
and its derivative $\partial_y\phi(y)=-\lambda_0-3\,\theta_3\,y^{2}$
has one positive root $y_{3}(\theta)=({-\lambda_0}/({3\,\theta_3}))^{1/2}$;
moreover, $\partial_y\phi\,\vert_{\,y=y_2(\theta)}<0$.
In~aim to find sufficient conditions for \eqref{E-auxcor-01}, consider the following functions:
\begin{eqnarray*}
 \phi_+(y) \eq P_{\,\phi_+}(y^2)/y,\ \ {\rm where}\ \quad P_{\,\phi_+}(z) = -\Psi_3^- z -\lambda_0, \\
 \phi_-(y) \eq P_{\,\phi_-}(y^2)/y,\ \ {\rm where}\ \quad P_{\,\phi_-}(z) = -\Psi_3^+ z -\lambda_0, \\
 \partial_y(\phi_-)(y) \eq P_{\,\partial_y(\phi_-)}(y^2)/y^2,\ \ {\rm where}\ \
 P_{\,\partial_y(\phi_-)}(z) = -3\,\Psi_3^+ z -\lambda_0,\\
 (\partial_y\phi)_+(y) \eq P_{\,(\partial_y\phi)_+}(y^2)/y^2,\ \ {\rm where}\ \
 P_{\,(\partial_y\phi)_+}(z) = -3\,\Psi_3^- z -\lambda_0\,.
\end{eqnarray*}
Then $y_2^+=({-\lambda_0}/{\Psi_3^-})^{1/2}$ and $y_2^-=({-\lambda_0}/{\Psi_3^+})^{1/2}$ are positive roots
of $\phi_+(y)$ and $\phi_-(y)$, and $y_3^+=({-\lambda_0}/{(3\Psi_3^-)})^{1/2}$ is the positive root of $(\partial_y\phi)_+(y)$.
 We need to examine when the resultant
\[
 R_3(t)=-{\rm Res}(P_{\,\phi_-},\,(1-t)P_{\,\partial_y(\phi_-)} +t\,P_{\,(\partial_y\phi)_+})/\lambda_0
 =-3(\Psi^+_3-\Psi^-_3)\,t + 2\,\Psi^+_3
\]
does not vanish for any $t\in[0,1]$; hence, the polynomials have no common roots.
We have $R_3(0) = 2\,\Psi_3^+ > 0$. Hence,
$3\,\Psi_3^- >\Psi_3^+$, see \eqref{E-a3}, provides $y_3^+ < y_2^- < y_2^+$, see \eqref{E-auxcor-01}.
\end{remark}

\subsubsection{Case of $\Psi_3=0$}
\label{sec:attr3}

Let $\Psi_3=0$, $\Psi_1>0$ and $\Psi_2>0$, see Section~\ref{R-burgers-heat}, point (c$_1$).
Then \eqref{Cauchy-stat1} becomes
\begin{equation}\label{Cauchy-stat2}
 \mathcal{H}(u) = \Psi_1(x)\,u^{-1} -\Psi_2(x)\,u^{-3}\,,
\end{equation}
where $\mathcal{H}(u):=-\Delta u -\beta\,u$. Certainly, Cauchy's problem \eqref{Cauchy} reads
\begin{eqnarray}\label{Cauchy3}
 \partial_t u + \mathcal{H}(u) = \Psi_1(x)\,u^{-1}-\Psi_2(x)\,u^{-3},\quad
 \quad u(x,0) = u_0(x)>0\,.
\end{eqnarray}
Then functions $\phi_-$ and $\phi_+$ in \eqref{difineqR-C}  become
\begin{eqnarray*}
 \phi_+(y) \eq P_{\,\phi_+}(y^2)/y^3,\ \ {\rm where}\ \quad P_{\,\phi_+}(z) = -\lambda_0 z^2 +\Psi_1^- z -\Psi_2^+ , \\
 \phi_-(y) \eq P_{\,\phi_-}(y^2)/y^3,\ \ {\rm where}\ \quad P_{\,\phi_-}(z) = -\lambda_0 z^2 +\Psi_1^+ z -\Psi_2^-,
\end{eqnarray*}
and $f(w,\,\cdot\,)=-\lambda_0\,w+(\Psi_1 e_0^{-2})\,{w^{-1}}-(\Psi_2 e_0^{-4})\,w^{-3}$.
It is easy to see that
 $\partial_w f(w,\,x\,)\le\partial_w\phi_-(w)$.
Assume that
\begin{equation}\label{E-lambda0-cond1}
 0<\lambda_0<{(\Psi_1^-)^2}/({4\Psi_2^+}).
\end{equation}
Each of functions $\phi_-(y)$ and $\phi_+(y)$ has two positive roots; moreover, $y_1^-<y_2^-$ and $y_1^+<y_2^+$.
Since $\phi_-(y)<\phi_+(y)$ for $y>0$, we also have $y_2^- < y_2^+$ and $y_1^- > y_1^+$.
 Denote by $y_3^- \in (y_1^-,\,y_2^-)$
a uni\-que positive root of $\partial_y\phi_-(y)=-\lambda_0 -\Psi_1^- y^{-2} +3\,\Psi_2^+ y^{-4}$.
Notice that $\phi_-(y)>0$ for $y\in(y_1^-,\,y_2^-)$ and $\phi_-(y)<0$ for $y\in(0,\,\infty)\setminus[y_1^-,\,y_2^-]$;
moreover, $\phi_-(y)$ increases in $(0,y_3^-)$ and decreases in $(y_3^-,\,\infty)$.
The line $z=-\lambda_0\,y$ is asymptotic for the graph of $\phi_-(y)$ when $y\to\infty$,
and $\lim_{\,y\downarrow 0}\phi_-(y)=-\infty$.
The function $\partial_y\phi_-(y)$ decreases in $(0,y_4^-)$ and increases in
$(y_4^-,\,\infty)$, where $y_4^-:=(6\,\Psi_2^+/\Psi_1^-)^{1/2}>y_3^-$,
and $\lim\nolimits_{\,y\rightarrow\infty}\partial_y\phi_-(y)=-\lambda_0$, see Fig.~\ref{fig:2}.
We conclude that
 $y_1^+ < y_1^- < y_3^- < y_2^- < y_2^+$.
Hence, the following function positive for $\sigma\in(0,\,y_2^- -y_3^-)$:
\begin{equation*}
  \mu^+(\sigma):=-\sup\nolimits_{\,y\ge y_2^- -\sigma}\,\partial_y\phi_-(y)
 =\min\{|\partial_y\phi_-(y_2^- -\sigma)|,\,\lambda_0\}.
\end{equation*}
Define closed in $C(F)$ nonempty sets
\begin{eqnarray*}
 && {\cal U}^{\;\varepsilon,\eta}=\{\tilde u\in C(F):\ y_2^- -\varepsilon\le \tilde u/e_0 \le y_2^+ +\eta\},\quad
 \varepsilon\in(0,\,y_2^- -y_1^-),\ \eta\in(0,\infty].
\end{eqnarray*}
We have $\mathcal{U}_{\,0}\subset{\cal U}^{\;\varepsilon,\eta}\subset\mathcal{U}^{\,\varepsilon,\infty}\subset\mathcal{U}_{\,1}$,
where the set ${\cal U}_{\;1}\!=\!\{\tilde u\in C(F): \tilde u/e_0 > y_1^-\}$ is open,
and $\,\mathcal{U}_{\,0}=\!\{\tilde u\in C(F): y_2^- \le \tilde u/e_0 \le y_2^+ \}$.

 The proof of next results is similar to the proof of Proposition~\ref{T-new2} and Theorem~\ref{thattract2}.

\begin{proposition}\label{T-new3}
Let \eqref{E-lambda0-cond1} holds. Then

$(i)$ for any $u_0\in{\cal U}^{\;\varepsilon,\eta}$, Cauchy's problem
\eqref{Cauchy3} admits a~unique global solution. Moreover, ${\cal
U}^{\;\varepsilon,\eta}$ are invariant sets for associated semigroup
$\,\mathcal{S}_t:u_0\to u(\cdot\,,t)\ (t\ge 0)$
in~$\mathcal{C}_\infty$;

$(ii)$ for any $\sigma\in(0,\varepsilon)$ there exists $t_1>0$ such that
$\;\mathcal {S}_t({\cal U}^{\;\varepsilon,\infty})\subseteq{\cal U}^{\;\sigma,\infty}$
for all $t\ge t_1$.
\end{proposition}

\begin{theorem}\label{thattract3}
\ $(i)$ If \eqref{E-lambda0-cond1} holds then \eqref{Cauchy-stat2}
has in $\,\mathcal{U}_{\,1}\cap C^\infty(F)$ a unique solution $u_*$, which obeys
$y_1^-\le u_*/e_0\le y_1^+$;
moreover, $u_*=\lim\nolimits_{\,t\to\infty} u(\cdot,t)$,
where $u$ solves~\eqref{Cauchy3} with $u_0\in\mathcal{U}_{\,1}$,
and for any $\sigma\in(0,\,y_2^- -y_3^-)$, the set $\mathcal{U}^{\,\sigma,\infty}$
is attracted by associated semigroup exponentially fast to $u_*$ in $C$-norm:
\begin{equation*}
 \|u(\cdot\,,t)-u_*\|_{C(F)}\le\delta^{-1}(e_0)\,e^{-\mu^+(\sigma)\,t}\| u_0-u_*\|_{C(F)}
 \quad(t>0,\ u_0\in\mathcal{U}^{\,\sigma,\infty}).
\end{equation*}

$(ii)$ Let $\beta,\Psi_1,\Psi_2$ be smooth functions on
$F\times\RR^n$ with a smooth metric $g(\cdot,q)$.
If \eqref{E-lambda0-cond1} holds for any
$F\times\{q\}\ (q\in\RR^n)$ then the solution $u_*$, see $(i)$, is smooth on $F\times\RR^n$.
\end{theorem}

\begin{remark}\label{R-05}\rm
Similarly, for $\Psi_2\equiv0$ when $\Psi_3=0$, $\Psi_1>0$, condition \eqref{E-lambda0-cond1} reduces to $\lambda_0>0$.
Each of the functions $\phi_-(y)=-\lambda_0\,y +\Psi_1^- y^{-1}$ and $\phi_-(y)=-\lambda_0\,y +\Psi_1^+ y^{-1}$
has one positive root $y_2^-=(\Psi_1^-/\lambda_0)^{1/2}$ and $y_2^+=(\Psi_1^+/\lambda_0)^{1/2}$;
moreover, $\partial_y\phi_-(y)<0$ for $y>0$.
\end{remark}

\end{document}